\documentclass[12pt]{article}
 \usepackage{amsmath,amsfonts,theorem}
\input epsf
 \newcommand{\ntag}[1]{} 
\numberwithin{equation}{section}
 
 \newtheorem{prop}{Proposition}

 \newtheorem{cor}{Corollary}

 \newtheorem{thm}{Theorem}

 \theorembodyfont{\rmfamily}
 \newtheorem{dfn}{Definition}

 \newtheorem{rmk}{Remark}

 \newcommand{\qed}{\ifhmode\unskip\nobreak\fi\quad\ensuremath\square}

 \newcommand{\Span}[1]{\left< #1 \right>}
 \newcommand{\half}{\textstyle\frac12}
 
 \newcommand{\id}{\operatorname{id}}

 \newcommand{\hF}{\operatorname{hF}}

 \newcommand{\CLRep}{\mathrm{CLRep}}

 \newcommand{\res}{\mathrm{res}}
 \newcommand{\rest}[1]{_{{\textstyle{|}}#1}} 

 \newcommand{\op}{\overline{\partial}}

  \newcommand{\p}{\partial}
 \newcommand{\w}{\widetilde}
 \newcommand{\tensor}{\otimes}
 \newcommand{\ov}{\overline}

 \newcommand{\sA}{\mathcal A} 
 \newcommand{\sB}{\mathcal B} 
 \newcommand{\sG}{\mathcal G} 
 \newcommand{\sD}{\mathcal D}
 
 \newcommand{\sC}{\mathcal C}
\newcommand{\sU}{\mathcal U}
 \newcommand{\sH}{\mathcal H}
 
 \newcommand{\sL}{\mathcal L}
 \newcommand{\sM}{\mathcal M}
 \newcommand{\Oh}{\mathcal O}
 
 \newcommand{\sS}{\mathcal S}
 
 \newcommand{\sE}{\mathcal E}
 \newcommand{\sP}{\mathcal P}

 \newcommand{\al}{\alpha}
 \newcommand{\be}{\beta}
 \newcommand{\de}{\delta}
 \newcommand{\ep}{\varepsilon}
 \newcommand{\fie}{\varphi}
 \newcommand{\ga}{\gamma}
 
 \newcommand{\om}{\omega}
 \newcommand{\si}{\sigma}
 \newcommand{\De}{\Delta}
 \newcommand{\Ga}{\Gamma}
 \newcommand{\La}{\Lambda}
 \newcommand{\Om}{\Omega}
 \newcommand{\la}{\lambda}
 \newcommand{\Si}{\Sigma}
 \newcommand{\na}{\nabla}

 \newcommand{\PP}{\mathbb P}
 \newcommand{\C}{\mathbb C}
 
 \newcommand{\Q}{\mathbb Q}
 \newcommand{\R}{\mathbb R}
 \newcommand{\Z}{\mathbb Z}

 \newcommand{\ad}{\operatorname{Ad}} 

\newcommand{\Diff}{\operatorname{Diff}}
\newcommand{\Tr}{\operatorname{Tr}}
 \newcommand{\End}{\operatorname{End}}

 \newcommand{\Hom}{\operatorname{Hom}}
 \newcommand{\Conj}{\operatorname{Conj}}
\newcommand{\CL}{\operatorname{CL}}
 \newcommand{\SL}{\operatorname{SL}}

 \newcommand{\SU}{\operatorname{SU}}
 \newcommand{\fsl}{\operatorname{\mathcal{\frak{sl}}}}
 \newcommand{\fsu}{\operatorname{{\frak{su}}}}
 \newcommand{\pt}{\mathrm{pt}}

\begin{document}

 \title{Quantization, Classical and Quantum Field Theory and Theta-Functions.}

  \author{\fbox{Andrei Tyurin}}
  \date{24.02.1940 -- 27.10.2002}

\maketitle

\begin{center}
 {\em To Igor Rostislavovich Shafarevich on his 80th birthday}
 \end{center}
 \bigskip

\section*{Introduction}

 Arnaud Beauville's  survey "Vector bundles on Curves and Generalized
Theta functions: Recent Results and Open Problems"
\cite{Be} appeared 10 years ago. This elegant survey is short (16 pages)  but
provides a complete introduction to a specific part of algebraic geometry.
To repeat his succes now we need more pages, even though we  assume that
the reader is already acquainted with the material presented
there. Moreover, in  Beauville's survey the  relation
between generalized theta functions and conformal field theories (classical and
quantum) was presented already.

Following  Beauville's  strategy we do not provide any proof or motivation.
 But we would like  to propose {\it all constructions} of
this large domain of mathematics in
such a way that the proofs    can be guessed
  from the geometric picture. Thus this text is not a mathematical monograph
   yet, but rather a digest of a field of mathematical investigations.

In the abelian case ( the
  subject of several beautiful classical books (\cite{B},
  \cite{C}, \cite{Wi}, \cite{F1} and  many others ) fixing some
  combinatorial structure
   (a so called
  theta structure of level $k$)   one obtains  special basis in the space
  of  sections of powers of the canonical polarization powers on jacobians.
  These sections can be presented  as holomorphic functions on the "abelian
  Schottky" space $(\C^*)^g$.
  This fact provides various applications of these  concrete analytic formulas
   to  integrable
  systems, classical mechanics and PDE's  (see the references in  \cite{DKN} ).

   Our practical goal
 is to do the same in the non-abelian case, that is, to give the answer to
 the final question of the Beauville's survey (Question 9 in
 \cite{Be}).

   It has been observed many time
   that the {\it construction} of theta functions with characteristics is
   intricately
   related to the paradigm of the quantization procedure (which is
   a quantum field theory in  dimension 1).
New features   came from Conformal Field Theory (which is a field theory
in dimension 2 = 1+1).
 This new stream brings the
standard physical paradigm of "symmetries, fields, equations etc ... and
 gluing properties corresponding to local
Lagrangians". New CFT methods provided powerful  computational tools while
"algebraic
geometers would have never dreamed of being able to perform such
computations" (A. Beauville).

In  future we hope to extend this digest to a mathematical monograph
with the title "VBAC".

{\it Acknowledgments.} These notes were written for a series of lectures given at the Centre de Recherches Mathematiques at the Universite de Montreal in September 2002. I am grateful to J. McKay, J. Hurtubise, J. Harnad, D. Korotkin
and A. Kokotov who made my stay in Montreal very agreeable and productive.  I would like to express my gratitude to
my collaborators C.~Florentino, J.~Mourao, J.P.~Nunes and
participants of Quantum Gravity seminar in IST (Lisboa, 2000). I would like to thank
 Jean Le Tourneux and Andre Montpetit for their assistance in the preparation of these notes. Special thanks  go to my daughter, Yulia Tiourina,
 for  catching  numerous mistakes and misprints in this huge text.

\tableofcontents

\section{Quantization procedure}

\subsection{The Framework}

The first question: what we want  to quantize? The
framework of geometric quantization is usually described as
follows ( see, for example \cite{GS1}). Let $(M, \om)$ be a
symplectic manifold which represents a classical mechanical
system with finite number of degrees of freedom or the bialgebra
$C^\infty(M)$ with the usual commutative multiplication  $f_1
\cdot f_2$ and Poisson structure given by well known construction:
 every $f \in C^\infty(M)$ defines a Hamiltonian vector field
\begin{equation}
H_f = \om^{-1}(df)
\end{equation}
and a pair $(f_1, f_2)$ defines the function
\begin{equation}
\{f_1, f_2\} = \om(H_{f_1}, H_{f_2})
\end{equation}
-Poisson brackets, such that
\begin{equation}
d \om = 0 \implies \text{ Jacobi equality}.
\end{equation}

If $M$ is compact, then even abelian multiplication is enough to
reconstruct $M$ and then the Lie algebra structure reconstructs
$\om$.

 So in this case the pair $(M, \om)$ and the (observables) algebra
$C^\infty(M)$ are equivalent objects.

 (Pre BRST)-quantization rules is a map $Q$ sending classical
objects to corresponding quantum objects:
\begin{enumerate}
\item $Q(M) = \sH$ is a Hilbert space (of wave functions);
\item $Q(f) \in Op(\sH)$ where $Op(\sH)$ is the space of
 self adjoint operators
\item which
should satisfy the {\it correspondence principle}
\begin{equation}
[Q(f_1), Q(f_2)] = i \cdot h \cdot Q(\{ f_1, f_2 \})
\end{equation}
(Dirac correspondence) and the representation of our Poisson
Lie-algebra  be irreducuble.
\end{enumerate}

Unfortunately, such construction couldn't exist at all due to van
Hove theorem. Nevertheless  one usually uses two basic examples:
Souriau - Kostant quantization doesn't satisfy the irreducibility
condition while  Berezin quantization doesn't satisfy the
correspondence principle.

But in both approaches we  extend classical mechanical
 date by
 a {\it prequantization date} $(L, \nabla)$, where $L$ is a
complex line bundle (with a Hermitian structure $h$) and $\nabla$
is a unitary connection on $L$ with the curvature form
\begin{equation}
F_\nabla = 2 \pi \cdot i \cdot \om.
\end{equation}

 This equality implies  very strong constraint to the
  symplectic structure:
\begin{equation}
[\om] \in H^2(M, \Z)
\end{equation}
that is the cohomology class of the symplectic form has to be
integer.
 Such quadruple
\begin{equation}
(M, \om, L, \nabla)
\end{equation}
 is ready to be quantized.
First of all the space of wave functions is the space of sections
\begin{equation}
\sH = \Ga^\infty (L)
\end{equation}
with the inner product given by the formula
\begin{equation}
<s, s'> = 1/n! \int_M (s, s')_h \om^{1/2 dim M},
\end{equation}
(usually one  extends this space to $L^2$ - sections). For every
function (a classical observable) we use the operator
\begin{equation}
  Q(f) = \nabla_{H_f} + 2 \pi i \cdot f,
\end{equation}
that is
$$
Q(f)(s) = \nabla_{H_f} s + 2 \pi i \cdot f \cdot s.
$$
Two simple exercises:
\begin{enumerate}
\item Dirac equality holds (2.4);
\item  the representation of
our Poisson Lie-algebra is very far  from being irreducible (even
for the simplest case $M = \R^2 = T^{*} \R)$
\end{enumerate}


Well known remedy is a choice of a {\it polarization}. A
polarization is a distribution that is a sub bundle of the
complexified tangent space
\begin{equation}
F \subset T^\C M , \quad rk F = 1/2 dim M;
\end{equation}

Fixing such polarization  the wave function space can be obtained
as
\begin{equation}
Q(M)= \sH^F = \{ s \in \Ga^\infty (L) \vert \nabla_F s = 0 \}.
\end{equation}
 There are two natural choices:
\begin{enumerate}
\item $\overline{F} = F$  is a {\it real polarization};
\item $\overline{F} = F^\perp $  is a
{\it complex polarization}, that is, a choice of an almost complex
structure $I$ such that $(F = T^{0,1})$.
\end{enumerate}

Of course for the general choice of a polarization we don't know
so much but there is a couple of good choices:

in the real case: choose $F$  integrable and moreover completely
integrable that is defined by a Lagrangian fibration:
\begin{equation}
\pi \colon M \to B , \quad \pi^{-1}(b_{gen}) = T^n
\end{equation}
that is, the generic fiber here is  a Lagrangian torus of
dimension $n = 1/2 dim M$ (see Fig.1);

 and  in the complex case: choose  an
almost complex structure which is integrable and compatible with
$\om $ such that the pair $(\om, I)$ gives a Kahler metric with
condition (1.6) that is a Hodge metric.

\begin{figure}[tbn]
\centerline{\epsfxsize=3in\epsfbox{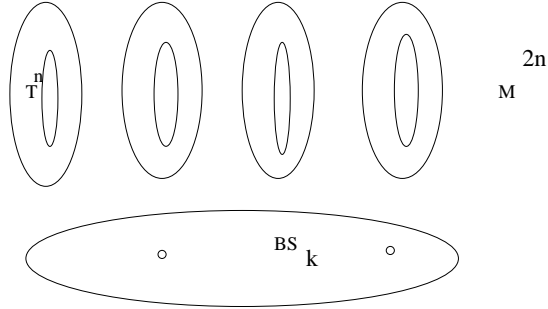}}
\caption{\sl Completely integrable system}
\label{F2}
\end{figure}

There is a couple of equivalent conditions for both of these
cases: for a real polarization, for any two vector fields
\begin{equation}
\p_1, \p_2 \in F \subset T^\C M
\end{equation}
the commutator
\begin{equation}
[\p_1, \p_2] \in F \subset T^\C M
\end{equation}
i.e.  vector fields from our real distribution form a Lie
subalgebra.
 (Otherwise their commutators form a quadratic form on $ F$ with
 coefficients in $T^\C M / F$ just like the second quadratic
 form in the standard differential geometry of a surface
 in $\R^3$ ).

For the complex case we have to add the equality of the
 same nature
\begin{equation}
[I \p_1, I \p_2] = [\p_1, \p_2] + I [I \p_1, \p_2]
+ I [\p_1, I \p_2]
\end{equation}
where we consider $I$ as a section of $End TM$ with
\begin{equation}
 I^2 = -1
\end{equation}
$$
\om(\p_1, I \p_2) = - \om (I \p_1,  \p_2)
$$
$$
\om(\p_1, I \p_2) > 0.
$$

These cases are lying in different domains of mathematics: the
real polarizations  are contained in Symplectic Geometry and the
complex polarization belongs to Algebraic Geometry. A
quantization is {\it perfect} if it  admits both of these
polarizations. So the perfect quantization belongs to the union of
Algebraic and Symplectic Geometries.

\subsection{The real polarization.(Symplectic Geometry)}

Recall that a  subcycle $T \subset M$ is Lagrangian if it is a
middle dimensional subcycle such that
\begin{equation}
\om \vert_T = 0.
\end{equation}
Thus  the restriction of the pair $(L, \nabla)\vert_T$ is  the
trivial
 line bundle with a flat connection
\begin{equation}
F_\nabla \vert_T = 2 \pi i \om \vert_T = 0.
\end{equation}
Recall that a gauge class of a flat connection is given by
 a character
of the fundamental group of $T$
\begin{equation}
\chi \colon \pi_1(T)     \to U(1).
\end{equation}

\begin{dfn} A Lagrangian cycle $T$ is Bohr-Sommerfeld
if $\chi = 1$, that is, if there exists a covariant constant
section of $(L, \nabla)\vert_T$.
\end{dfn}
 To kill a character,  we need $rk H_1(T, \Z)$ conditions
 so we may expect that
 in a family
of Lagrangian cycles the subfamily of Bohr-Sommerfeld cycles has
this codimension. Moreover, it was shown
  that in any small
 Darboux-Weinstein
 neighborhood of a smooth
 Bohr-Sommerfeld cycle $T$ every other Bohr-Sommerfeld cycle $T'$
 intersects
 $T$.

For example,  our real polarization (1.13) defines a family
 of Lagrangian
cycles parameterized by a space $B$ of dimension $n = 1/2 \dim M$.
But  a general fiber of any real polarization fibration $\pi$
(1.13) is a n-dimensional torus, that is
$$
 rk H_1(T^n) = n = dim B.
$$
 Thus the set of Bohr-Sommerfeld fibers is discrete and if $B$ is
 compact it forms
  a finite set $BS \subset B$ of Bohr-Sommerfeld fibers of $\pi$.

 Then the Hilbert space of wave functions is given as
 the direct sum of lines
\begin{equation}
 Q_F(M) = \sH_\pi = \oplus_{b \in BS \subset B} \C \cdot
s_b
\end{equation}
where $s_b$ is a covariant constant section of $(L,
\nabla)\vert_{\pi^{-1}(b)}$).


Recall that every covariant constant section is defined up
to $U(1)$-scaling,
so we can define naturally an Hermitian scalar product on this space.

We can extend  the notion of
 Bohr-Sommerfeld Lagrangian cycles:

\begin{dfn} A Lagrangian cycle $T$ is Bohr-Sommerfeld of level $k$ (or $BS_k$ - cycle)
 if $\chi^k = 1$, that is, if there
exists a covariant constant section of $(L^k, \nabla_k)\vert_T$,
where $\nabla_k $ is the connection generated by $\nabla$ on the
tensor power $L^k$ of the line bundle $L$.
\end{dfn}

Reasoning  the same way as before,  we can prove that the subset
$BS_k \subset B$ is  finite and that the Hilbert space of wave
functions of level $k$ is given by the direct sum of lines
\begin{equation}
 Q_F^k(M) = \sH_\pi^k = \oplus_{b \in BS_k \subset B} \C \cdot
s_b
\end{equation}
where $s_b$ is a covariant constant section of $(L^k,
\nabla_k)\vert_{\pi^{-1}(b)}$).

 But we are quantizing a classical dynamical system $(M, \om)$!
So the result does not have to  depend on additional data  like
$(L, \nabla, \pi, I)$ and so on.

In our construction the wave space $\sH_\pi$ depends on a gauge
class of a connection $\nabla$ only and this class is given by
the curvature form if $M$ is simply connected (see \cite{AB}).
(Recall that in non-simple connected case
 we have to fix additionally  a point of the
 "jacobian" $H^1(M, \R) / H^1(M, \Z)$ of $M$ (see below).

  In any case
 our wave function spaces depend up to a shift of characters on the
curvature tensor which is $\om$ and on a choice of a real
polarization $\pi$.

\begin{figure}[tbn]
\centerline{\epsfxsize=3in\epsfbox{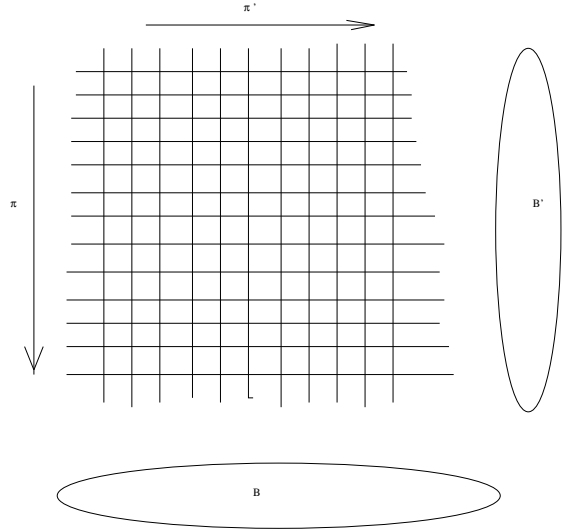}}
\caption{\sl twice integrable system}
\label{F4}
\end{figure}

But the most important thing is that the projectivization of the
wave functions
 space doesn't depend on a real polarization $\pi$ too.
  For any other Lagrangian fibration $\pi' \colon M \to B'$,
Kostant defined the canonical (up to constant )
 pairing between
$\sH_\pi$ and $\sH_{\pi'}$ (see for example \cite{JW2}). So we
get  the identification
\begin{equation}
 \PP \sH_\pi = \PP \sH_{\pi'}
\end{equation}
of projective spaces.

 Thus $\sH_\pi$  is the wave space of a quantization of $(M, \om)$
 and of nothing else.

\subsection{Kahler quantization. (Algebraic Geometry) }

In this case we  have a quadruple $(M, \om, \nabla, I)$ where a
complex structure $I$ is compatible with $\om$. Thus the pair
$(\om, I)$ defines a Kahler metric with the Kahler form $\om$.
Moreover, $\om$ is of Hodge type (1,1). Thus the connection
$\nabla $ on $L$ has the curvature form of type (1,1) and defines
a holomorphic structure on $L$ (see for, example, \cite{GH}).
 As
Kodaira (1958) proved, a Kahler metric  with  integer cohomology
class $[\om] \in H^2(M, \Z)$ of the Kahler form is a Hodge
metric, hence $M_I$ is a projective algebraic variety and $L =
\Oh_M (1)$ is the class of  hyperplane section of $M_I$.

 So, the complex polarizations bring us to Algebraic Geometry!

Naturally, the wave function space
\begin{equation}
Q_F(M) = \sH_I = H^0(M_I, L)
\end{equation}
is the space of holomorphic sections of the polarization $L$ and
 the projectivization $\PP H^0(M_I, L)$ is the classical
 complete linear system on $M_I$.

The Hermitian structure on $L$ gives an Hermitian structure on
$H^0(M_I, L)$ and provides the identification
\begin{equation}
 \sH_I = H^0(M_I, L) = H^0(M_I, L)^*
\end{equation}
   with changing of the complex structure to the complex
conjugate.

The standard construction of the Algebraic Geometry gives the
following map
\begin{equation}
 \phi \colon M_I \to \PP H^0(M_I, L)^*,
\end{equation}
that is, under the previous identification
\begin{equation}
 \ov{\phi} \colon M_I \to \PP H^0(M_I, L).
\end{equation}

Images of points of $M_I$ are called {\it coherent states}.  A
priori such image is defined up to a transformation
 from $PGL( \PP H^0(M_I, L))$. But its every
 {\it stable} orbit contains  unique
unitary orbit preserving our form  $\om$. If our form $\om$ is the
restriction of
 the Fubini-Study form
from $ \PP H^0(M_I, L)$, then we can switch on to  Geometric
Invariant Theory and reduce the theory of coherent states to it.
So in algebraic geometrical set up we can assume that our form is
a restriction of the Fubini-Study form. In the general case we
have to be more careful and to add so called J.Rawnsley
$\ep$-function of
 (see  \cite{R}).

Of course we lose the Dirac correspondence and have to use the
Berezin-Toeplitz quantization rules (see for example \cite{R}).
Thus we  introduce the level parameter $k \in \Z^+$. The space of
wave function of level $k$ is the space of holomorphic sections
\begin{equation}
 \sH_I^k = H^0(M_I, L^k).
\end{equation}
The Dirac correspondence can be restored in the
{\it quasi classical limit} $k \to \infty$ or $1/k \to 0$.
 So,  our complex polarization becomes a polarization in
 algebraic geometrical sence. Recall that if $k$ is large enough,
 then $\ov{\phi}$ (1.27) is an embedding  and the dimension of the
 corresponding projective space is given by the Riemann-Roch
 formula as a pure topological invariant.

So, every underlying symplectic manifold $(M, \om)$ of
 a polarized algebraic   variety $M_I, \Oh(1)$  with a
 natural complex polarization $\Oh(1)$ gives  the collection
  of wave functions
 spaces
$H^0(M_I, L^k) = H^0(M_{I}, \Oh(k))$

 But the quantization set up disturbs
  the Algebraic Geometry Peace.
The main question is:

do these spaces appear as a result of a successful quantization
procedure of a classical system given by the  underlying
symplectic manifold $(M, \om)$?

 That is,

is the projectivization of this space independent  on the complex
structure $I$ ?

(In physical symbols, $\frac{\p \PP H^0(M_I, L)}{\p I} = 0$ .)

 Mathematically correct question is the following:
let  $ \sM $ be the moduli space of such complex structures, that
is,  the moduli space of polarized  (in the algebraic geometryical
sense) varieties and
\begin{equation}
p \colon V_k \to \sM, \quad p^{-1}(M_I) = H^0(M_I, L^k)
\end{equation}
- a vector bundle of holomorphic sections of polarization line
bundles. Note that by the Kodaira vanishing theorem it is a vector
bundle indeed! Then
   there exists a holomorphic projective flat connection on $V$.

Existence of such connection constraints the topology type of the
vector bundle: if a vector bundle $V$ admits a projective flat
connection, then
 \begin{equation}
  c_t(V) = (1 - c \cdot t)^{rk V}.
\end{equation}
In particular,
\begin{equation}
  c_1(V) = rk V \cdot c, \quad c \in H^2(\sM, \Z).
\end{equation}
This is  very strong topological constraint!

 Thus   problems of new type are
 given inside of Algebraic Geometry:

 When a pair $M_I, \Oh(1)$ of an algebraic variety with
 a polarization is a result
  of a quantization procedure
 for the underlying symplectic manifold ?

Before  considering the simplest example when $\dim M = 2$ we give
some   general remarks:
\begin{enumerate}
\item of course our vector bundles $V_k$ (1.29) are
holomorphic but the  required projective flat connection could be
non holomorphic;
\item so, it is quite reasonable to construct may be not
holomorphic projective Hermitian connection;
\item in all classical examples
holomorphic projective flat connections are Hermitian by the
construction;
\item but in the last non-abelian case (which is of the main
  interest here) this fact isn't known
 up to now;
\item if such new connection exists automatically we get a Higgs field
on every vector bundle (1.29) on moduli spaces of complex
structures.
\item In the physical set up it may be  that we do not obtain  vector
 bundles (1.29) but have the collection of projective bundles
 which {\it aren't projectivizations of vector bundles};
\item even in the infinite dimensional case (when Hilbert  vector
 bundles are always trivial) the projective bundles can be non trivial
  and are given by a cohomology class from $H^3(\sM, \Z)$. This
Brauer group gives us non trivial analog of Atiyah's K-functor.
\end{enumerate}

In the conformal field theory context this connection always {\it
has to be holomorphic} but we do not have space and time to
discuss it here.

Because of the time-space problem this very interesting subject
is out of our considerations here.

\subsection{Extended Kodaira-Spencer theory}

The theory of deformations of complex structures was developed by
Kodaira and Spencer many years ago. Trying  "theoretically" to
solve  the
 problem of  successful
Kahler quantization we have to lift this theory to the problem
of
 holomorphic deformations of a triple $(I, L_I, (s)_0)$,
  where $I$ is an integrable complex structure, $L_I$ is
  our holomorphic line bundle and $s \in H^0(L_I)$ is
   zero divisor of a holomorphic section of this line bundle.
   The best reference for such lifting is   Welters paper
   \cite{W}. Recall that the holomorphic structure of $L_I$
    codes the connection $\nabla$ from the quadruple $(1.7)$ by the decomposition
\begin{equation}
 \nabla = \nabla^{1, 0} + \nabla^{0, 1}
\end{equation}
where
\begin{equation}
  \nabla^{1, 0} = 1/2 \cdot(1 - i I) \nabla, \quad \nabla^{0, 1} =
  1/2 \cdot (1 + i I) \nabla.
\end{equation}
Our history starts with an infinitesimal deformation of the
complex structure $\dot I$. Since it is an anti involution $I^2 =
-1$ we have the equality
\begin{equation}
\dot I I + I \dot I = 0
\end{equation}
and one can consider $\dot I$ as a projection of $-i$ eigenspace
of $I$ to the $+i$-eigenspace. That is
\begin{equation}
\dot I \in \Om^{0, 1} (M_I,T^{1, 0})
\end{equation}
(just as  usual Beltrami differential).

Linearization of the integrability condition (1.16) gives us
\begin{equation}
(\op \dot I)^{2, 0} = 0,
\end{equation}
that is, the (0, 2)-component of this tensor is trivial.

Suppose that we have the following list of properties.
\begin{enumerate} \item
Our polarization $L$ is a rational part of the canonical class $K$
of $M_I$. That is
\begin{equation}
- [K] = n [c_1(L)].
\end{equation}
\item Further on we suppose that the Levi-Civita connection of
 the Kahler
 metric $\om(\quad , I \quad )$ (1.17) on $K $ is proportional to
  the unitary  connection on $L$.
\item We consider a sufficiently large level $k$ such that
\begin{equation}
i >o \implies H^i(M_I, \Oh(k \cdot c_1(L) + 1/2 K)) = 0.
\end{equation}
(For example, $M_I$ is a Fano-variety of index $n$ and
 $L = \Oh(k)$).
\item Then we have the symplectic Kahler form
\begin{equation}
( k + n/2) \om.
\end{equation}
\end{enumerate}

From this it is easy to see (for example, by  local
representation of our tensor) that
 there exists a  section
\begin{equation}
W_0 \in \Ga^\infty (S^2 T^{1, 0}M_I) \in (S^2)^{2, 0}
\end{equation}
where $S^2$ is the second symmetrical power of the
holomorphic tangent bundle such that the convolution of
 this tensor with our symplectic form gives us $\dot I$:
\begin{equation}
W_0 * (k + n/2) \om = \dot I.
\end{equation}
(Recall that $\om$ is of a type $(1.1)$ with respect to the
complex structure $I$.)

 To see this very important tensor we have to represent tensors
 as homomorphisms of tensors vector bundles!
  By this notation we
   emphasize
  that this construction is due to Welters \cite{W}.
   So our symplectic form
\begin{equation}
(k+n/2)\om \colon T^{1, 0} \oplus T^{0, 1} \to \Om^{1, 0} \oplus \Om^{0, 1}
\end{equation}
is of the Hodge type $(1,1)$ and it means that
\begin{equation}
\om(T^{1, 0}) \in \Om^{0,1} \text{ and } \om(T^{0, 1}) \in \Om^{1,0}.
\end{equation}
On the other hand,
\begin{equation}
\dot I \colon T^{0,1} \to T^{1,0}.
\end{equation}
Thus the composition of these homomorphisms
\begin{equation}
\dot I \circ \frac{1}{k+n/2} \om^{-1} = W \colon  \Om^{1,0} \to T^{1,0}
\end{equation}
is the very one symmetric tensor.

Now consider the cup product map
\begin{equation}
H^0(M_I, T^{1, 0}) \otimes H^1(M_I, \Om^{0, 1}) \to H^1(M_I, \Oh)
\end{equation}
and the restriction of it to
 $H^0(M_I, T^{1, 0}) \otimes [(k+n/2)\om]$:
\begin{equation}
[(k+n/2)\om] \colon H^0(M_I, T^{1, 0}) \to H^1(M_I, \Oh).
\end{equation}
Suppose that this in an {\it isomorphism}. For example
\begin{equation}
H^0(M_I, T^{1, 0}) = H^1(M_I, \Oh) = 0
\end{equation}
or $M_I$ is a polarized  abelian variety.

In this case our holomorphic line bundle $L_I$ can't be deformed
continuously
 under deformation of the complex structure.
  So we have to investigat  deformations of
 holomorphic sections only.

Having such symmetric form and following Nigel Hitchin let us
infinitesimally deform a section $s \in H^0(M_I, L_I)$ to a $\dot
I$-holomorphic section. Differentiating along a "time" $t$ the
condition $(1 + i I_t) \nabla s_t = 0$ we get the equality
\begin{equation}
i \dot I \nabla^{1,0} s = \nabla^{0,1} a
\end{equation}
where $-a$ is  a section of $L$. Such section depends linearly on
$s$ and $\dot I$ and defines a connection on the bundle $V$
(1.29).
 The main observation of
the extended Kodaira-Spencer theory is the fact that
 this equation admits
a cohomological interpretation. Namely LHS and RHS of (1.49) are
results of applying  first order
 differential
operators on $L$. So now we have to consider the vector
bundle $D^1(L_I)$ of
holomorphic differential operators on $L_I$ applied to
a holomorphic section
 in LHS and to the potential of a connection.

Let the exact sequence
\begin{equation}
0 \to \Om^{1,0} \otimes L \to J^1(L) \to L \to 0
\end{equation}
be the extension defining 1-jet bundle of our line bundle $L$. By
the Atiyah (see \cite{A2}) theorem, the cocycle of this extension
is
\begin{equation}
[c_1(L)] = [\om] \in H^1(M_I, \Om^{1,0}).
\end{equation}

Now
\begin{equation}
D^1(L) = Hom (J^1(L), L )
\end{equation}
is the bundle of the sheaf of germs of differential operators
on $L$ of order
$\leq 1$. So this sheaf is the extension
\begin{equation}
0 \to \Oh \to D^1(L) \to T^{1, 0} \to 0
\end{equation}
(given by the class $- [\om]$).

The long exact cohomology sequence of this exact triple
 decays in two parts
\begin{equation}
0 \to H^0(M_I, \Oh) \to H^0(M_I,D^1(L)) \to 0
\end{equation}
and
\begin{equation}
0 \to H^1(M_I, D^1(L)) \to H^1(M_I, T^{1, 0}) \to H^2(M_I, \Oh)
\dots
\end{equation}
because they are divided by the homomorphism (1.44).
 The first exact sequence states
 that an every globally defined first order holomorphic operator
 on $L$ is
 a {\it multiplication by a constant} and the second states
 that the {\it
 symbol map} $H^1(M_I, D^1(L)) \to H^1(M_I, T^{1,0})$ is a
 { \it monomorphism}.

Now for a holomorphic section  $s \in H^0(M_I, L_I)$ the
 valuation morphism
\begin{equation}
d_s \colon D^1(L) \to L , \quad d_s (D) = D s,
\end{equation}
gives the complex
\begin{equation}
0 \to D^1(L) \to L  \to 0
\end{equation}
which hypercohomology $H_h^1(d_s)$  defines linear infinitesimal
deformations of a triple $(\dot I, L_I, s)$.

Each class of such hypercohomology can be done as a pair
\begin{equation}
(\{ \dot s_i\}, \{\al_{i,j} \}) \in C^0(\sU, L) \oplus C^1(\sU, D^1(L))
\end{equation}
where $\sU$ is a Cech covering. Indeed it is a 1-cocycle of the
total complex, associated with the double complex
\begin{equation}
\pm d_s \colon C^i(\sU, D^1(L)) \to C^i(\sU, L).
\end{equation}
The spectral sequence of hypercohomology (of a double complex)
yields  the exact sequence
\begin{equation}
H^0(M_I, D^1(L)) = \C \to H^0(L_I) \to H_h^1(d_s) \to H^1(D^1(L)) \to H^1(L_I)
\end{equation}
corresponding to the pairs (1.56).

On the other hand,    $d_s$-complex (1.57) defines the exact
quadruple
\begin{equation}
0 \to ker_s \to D^1(L) \to L \to coker_s \to 0
\end{equation}
and the second spectral sequence gives
\begin{equation}
0 \to H^1(ker) \to H_h^1(d_s) \to H^0(coker) \to H^2(ker).
\end{equation}
It is easy to see (this observation due to Welters \cite{W})
 that the support
of the sheaf $coker_s$ can be very small:
\begin{equation}
Supp F = Sing (s)_0
\end{equation}
-the singularity locus of the zero-divisor of a section $s$. In
particular, if a zero divisor
 is smooth (this is a general case)
then $F=0$  and the exact quadruple (1.60) is in fact a triple.

Now let us return to our main equation (1.49).
The RHS admits a solution if the
\begin{equation}
(\op LHS )^{0,2} = 0
\end{equation}
(by the  Dolbeault lemma). So we have to compute this RHS Hodge
type. But by ( 1.41 ) it is easy to see that
\begin{equation}
(\op (i \dot I \nabla^{1,0}))^{0,2} = 0.
\end{equation}
Here we are using
\begin{equation}
i \dot I \nabla^{1,0} \in \Om^{0,1} (M_I, D^1(L))
\end{equation}
and the infinitesimal integrability condition (1.34).
Now we have the interpretations of $i \dot I \nabla^{1,0}$ and
 $-a$ from (1.49)
as (0, 1) - forms with coefficients in $D^1(L)$ and we have
the complex
\begin{equation}
C^i = \Om^{0,i} (M_I, D^1(L)) \oplus \Om^{0, i-1}(M_I, L_I)
\end{equation}
with the differential
\begin{equation}
d(s) (D, a) = (\op D, \op a + (-1)^{i-1} D s)
\end{equation}
Now
\begin{equation}
\op s = 0 \implies d(s)^2 = 0
\end{equation}
and our pair
\begin{equation}
(i \dot I \nabla^{1,0}, -a) \in C^1.
\end{equation}
Moreover, equations (1.65) and (1.49) give
\begin{equation}
d(s) (i \dot I \nabla^{1,0}, -a) = 0.
\end{equation}
Thus, every solution to (1.49) defines a  hypercohomology class
 in $H_h^1(d_s)$!

Now suppose we can lift every pair
\begin{equation}
(\dot I, s) \in \Om^{0,1}(M_I, T^{1,0}) \oplus H^0(M_I, L_I)
\end{equation}
to a cocycle from $H_h^1(d_s)$ such that the natural "symbol" map
\begin{equation}
\si \colon H_h^1(d_s) \to H^1(M_I, T^{1,0})
\end{equation}
sends our pair to the cohomology class $[\dot I] \in H^1(M_I, T^{1,0})$.
That is there exists a map
\begin{equation}
A \colon  \Om^{0,1}(M_I, T^{1,0}) \oplus H^0(M_I, L_I) \to H_h^1(d_s)
\end{equation}
such that the composition $\si A$ is the map to the cohomology class.
Then such map defines a projective  connection on the vector bundle
of holomorphic sections (1.29).

Indeed, we can presents a hypercohomology class as a pair (1.58)
\begin{equation}
(D, s') \in \Om^{0,1} (M_I, D^1(L)) \oplus \Om^{0, 0}(M_I, L_I).
\end{equation}
But now the symbols of operators $-i \si (D)$ and $i \dot I
\nabla^{1,0}$ are cohomological  and, therefore the corresponding
operators are cohomological in $\Om^{0,1}(M_I, D^1(L))$    because
map to symbols is a monomorphism (1.55). Thus there exists a
global operator
 $\sD \Om^0(M_I, D^1(L))$ such that
\begin{equation}
D - i \dot I \nabla^{1,0} = \op \sD.
\end{equation}
Since $(D, s')$ is a cocycle ($d(s)(D, s') = 0 $, (see (1.71))
 one gets
\begin{equation}
Ds + \op s' = i \dot I \nabla^{1,0} s + \op (\sD s + s').
\end{equation}
Thus
\begin{equation}
-a = \sD s + s'
\end{equation}
is a solution for (1.49).

Now for two such solutions hypercocycles
$(i \dot I \nabla^{1,0}, a_1)$
 and $(i \dot I \nabla^{1,0}, a_2)$ are cohomological and
 the hypercocycle
\begin{equation}
(0, a_1 - a_2) = d(s) \sD'.
\end{equation}
 But from the definition of the complex (1.59) we have
 \begin{equation}
\op \sD' = 0 \implies a_1 - a_2 = \la s
\end{equation}
since every first order holomorphic operator is the multiplication
by a constant. Geometrically this means that we get a holomorphic
connection on
 the projectivization of our vector bundle.

Thus our last task is to construct the lifting $A$ (1.74).
But there is a canonical way to associate such lifting
 with a holomorphic symmetrical tensor
$W$ (1.40).

Consider the sheaf $D^2(L)$ of second order  holomorphic
operators given by the
standart extension
 \begin{equation}
0 \to D^1(L) \to D^2(L) \to S^2(T^{1,0}) \to 0.
\end{equation}
Evaluating differential operators on a given section $s \in
H^0(M_I, L_I)$ we get two homomorphisms
\begin{equation}
d_s \colon D^1(L) \to L
\end{equation}
$$
d'_s \colon D^2(L) \to L
$$
giving two complexes: our complex (1.57) and its analog
\begin{equation}
0 \to D^2(L) \to L \to 0.
\end{equation}
Considering the "trivial" third complex
\begin{equation}
0 \to S^2(T^{1,0}) \to 0 \to 0
\end{equation}
we get the homomorphism of the vertical triples  to the triple
$$
0 \to L \to L \to 0 \to 0
$$
and the exact triple of these complexes. The corresponding exact
sequence of hypercohomology has the form
\begin{equation}
0 \to H_h^0(d_s) \to H_h^0(d'_s) \to
\end{equation}
$$
\to H^0(M_I, S^2(T^{1,0})) \to H_h^1(d_s) \to H_h^1(d'_s) \to \dots
$$
and every holomorphic quadratic form on the tangent bundle defines
a  1-hypercohomology class (which may be zero). Of course our
required map depends on a choice of a section $s$, so we denote
this map by the symbol
\begin{equation}
\phi_s \colon H^0(M_I, S^2(T^{1,0})) \to H_h^1(d_s).
\end{equation}
The composition of this map with the forgetful homomorphism $\si$
(1.73) gives the homomorphism
\begin{equation}
\mu_L \colon H^0(M_I, S^2(T^{1,0})) \to H^1(M_I, T^{1,0})
\end{equation}
which doesn't depend on the choice of $s$. More precisely
\begin{equation}
\mu_L (W_0) = -W * [c_1(L)] + \mu_{\Oh}
\end{equation}
It is quite natural to call it the "cohomological heat equation".
 For example
for polarized abelian varieties the map $\mu_\Oh = 0$. Moreover,
all extensions $D^n(\Oh)$ are trivial hence all
 coboundary homomorphisms are trivial.

In particular, from this equality and (1.41) we can see that
\begin{equation}
\mu_\Oh (W_0) = - W_0 * 1/2 K
\end{equation}
because we can say more about sheaves of differential operators
on the trivial line bundle. Namely the extension class for all
jet-bundles is equal to
\begin{equation}
1/2[K_{M_I}] \in \Om^{1,1}.
\end{equation}
From this it is easy to see that
\begin{equation}
\mu_L (W_0) = -W_0 * ([c_1(L)] + 1/2 [K_{M_I}])
\end{equation}
where $K_{M_I}$ is the canonical class of $M_I$ and the
 cohomology class
$[K_{M_I}]$ depends on $\om$ only (and not on the choice of an
admissible complex
 structure).

 Of course we can get more precise formula using representatives
 of classes by
 forms and using the Levi-Civita connections for tensor bundles
  (see \cite{H1})

So the extended Kodaira-Spencer theory gives us (canonically) a
holomorphic projective connection on the vector bundles (1.29) but
of course to compute a curvature of one we have to present it as
 a Dolbeault or Cech class.

An existence of Welters tensor $W_0$ (1.40) gives strong geometrical
constraints. In particular $M_I$ can't be a variety of general
type. More precisely the canonical class $K$ of $M_I$ can't be
positive.

\subsection{Faithful functors}

      Suppose our $M_I$ which we consider as a point $I$ of
       the moduli space $\sM$
of all its deformation defines a new polarized algebraic variety
$M'_{I'}$  which also we consider as a point $I'$ of the
moduli space $\sM'$
such a way, that
 we can reconstruct $M_I$ (from the geometry of $M'_{I'}$).
 Geometrically, mapping $M_I$ to $M'_{I'}$ defines a holomorphic
embedding of moduli spaces
\begin{equation}
J \colon \sM \to \sM',
\end{equation}
and $J(I) = I'$. If a morphism of polarized varities from $\sM $
induces a morphism of corresponding varieties from $\sM'$ then
the map $M_I \to M'_{I'}$ is called a {\it faithful
 functor}. From geometrical point of view these geometrical
 images are  undistinguishable in spite of the fact
 that these varieties can have different dimensions and so on.
  An algebraic geometer is skillful enough if he can recognize
  such geometrically undistinguishable varieties of different
   dimensions and different shapes.

For example, a smooth intersection
\begin{equation}
X = Q_1 \cap Q_2 \subset \PP^5
\end{equation}
 of two quadrics in five dimensional projective space is
  geometrically equivalent
 to an algebraic curve $\Si_I$ of genus 2.
 To reconstruct $\Si_I$ from $X$ it is enough to remark that
\begin{enumerate} \item By the Lefshetz theorem $X$ is
simply connected;
\item $\Oh_X(-2) = K_X$ - the canonical class of $X$.
 Hence the embedding of $X$ into $\PP^5$ is canonical.
\item The space of all quadrics in $\PP^5$ through $X$ is
 a pencil
\begin{equation}
 \vert 2 H - X \vert = \PP^1
\end{equation}
that is parameterized by the projective line.
\item The subset of singular quadrics in this pencil is
 six different points
\begin{equation}
 W =p_1 \cup ... \cup p_6 \subset \PP^1.
\end{equation}
\item The double cover
\begin{equation}
 \phi \colon \Si_I \to \PP^1
\end{equation}
 with  ramification in these six points $W$ is the required curve
$\Si_I$ of genus 2.
\end{enumerate}

Note that the previous double cover is given by the canonical
 linear system of $\Si_I$.
Now a simple exercise in  linear algebra shows that there exists
just one pencil of quadrics in $\PP^5 $ up to a linear
transformation with the subset of singular
 quadrics which is equal to the  ramification locus of the
  double cover.

If $M_I \to M'_{I'}$ is a faithful functor and geometries of both
varietes are undistinguishable it is quite reasonable to consider
a quantization of $\sM'$ as a quantization of $\sM$, that is,
 to say
\begin{equation}
Q(M_I) = Q(M'_{I'})
\end{equation}
and so on ... We will say that such quantization is successful if
the
 quantization of $\sM$ is successful, that is, corresponding
  vector bundles (1.29)
on $\sM'$ admit projective flat holomorphic connections.

\subsection{Perfect quantization}

Suppose our algebraic variety $(M_I, \Oh(1))$ with a Hodge form
$\om$ admits as a symplectic manifold a real polarization that is
a Lagrangian fibration:
\begin{equation}
\pi \colon M \to B , \quad \pi^{-1}(b_{gen}) = T^n
\end{equation}
that is the generic fiber is  a Lagrangian torus of dimension
$n = 1/2 dim M$;

Again we may suppose that the subset $BS_k \subset B$ is
 finite and the Hilbert space of wave functions of level $k$
 is given by the direct sum of lines
\begin{equation}
 Q_F^k(M) = \sH_\pi^k = \oplus_{b \in BS_k \subset B} \C \cdot
s_b
\end{equation}
where $s_b$ is a covariant constant section of $(L^k,
\nabla_k)\vert_{\pi^{-1}(b)}$).

Suppose we can construct a natural isomorphism
\begin{equation}
  \PP Q_F^k(M) = \PP H^0(M_I, \Oh(k)).
\end{equation}
Then we have a perfect quantization providing a lot of
 beautiful and very important
properties.
\begin{enumerate}
\item First of all, both kinds of quantization are successful
 : indeed, the LHS of the equality (1.100) doesn't depend on a complex
  polarization at all. On the other hand, RHS doesn't depend on
  a real polarization.
\item LHS of the equality is decomposed into a sum of lines
 so the RHS as a space of sections admits a special basis
 (may be after fixing some additional structure like  so called
theta - structure for abelian varieties \cite{Mum}).

\item A basis of such type is called the Bohr-Sommerfeld basis.
\item Such situation can be called  {\it perfect quantization.}
\end{enumerate}

Moreover, the first example of a perfect quantization is given by
the classical theory of theta-functions.

 Obviously it is not the simplest example.
The simplest example is the following: consider the standard two
sphere $S^2$ realized as the complex Riemann sphere with the
standard $U(1) = S^1$-action (rotations around North and South
Poles). Now  $U(1)$-invariant symplectic form $\om$ of volume 1
gives the phase space of the classical mechanical system $S^2,
\om$. To quantize this one we have to fix a complex structure, but
it is unique and the real polarization
\begin{equation}
  \pi \colon S^2 \to [-1, 1]
\end{equation}
is given by the projection to the rotation axis.

As a prequantization date we get the line bundle $\Oh(1)$ of
degree 1 with the standard $U(1)$-invariant unitary  connection
$\nabla$. Then
 the wave function spaces
\begin{equation}
  \sH^k_I = H^0(\Oh(k))
\end{equation}
and $U(1)$- action on this space defines the special eigenbasis of
 Fourier polynomials. Fourier monomials
  are in 1-1 correspondence with
 Bohr-Sommerfeld fibers of the projection (1.101) in angle
 coordinates.
 Thus in this case we have the perfect quantization.

 It is easy to see that in this case all "good" theoretical
 conditions
 as (1.48) are
 disturbed but the result of quantization is successful.

 If we put out the word "natural" from the statement (1.100)
  we get the statement
 about ranks of both wave functions spaces. In this case we say that
 the quantization
 is {\it numerically} perfect. The geometry behind this
  numerical coincidence
 will be described in (2.25) - (2.30) and \cite{T3}.
 Such rank's equalities
 were proved for K3 surfaces cases and many others. But historically first case
 was the Gelfand-Cetlin system (see \cite{GS2}).

 \section{Algebraic curves = Riemann surfaces}

\subsection{Direct approach}

 Let $\Si$ be a compact smooth oriented Riemann surface and
$\Si_I$ be a complex structure on it, thus, one has an algebraic
curve of genus $g > 1$. Such curves admit the canonical
 polarization by the cotangent bundle $T^* \Si_I = \Oh(K)$
  where $K$ is the canonical class of $\Si_I$
  (see for example \cite{DSS}). Moreover, if we fix a metric
  on $\Si$ with the conformal
 class given by the complex structure $I$ then
 the Levi-Civita connection on the cotangent
 bundle gives the prequantization data $(T^*\Si, \nabla_{LC})$.
 This pair is equivalent to
the holomorphic pair $(\Si_I, \Oh(K))$.

The union of all spaces $H^0(\Oh(k K))$ of holomorphic sections
gives the collection of vector bundles
\begin{equation}
p_k \colon V_k \to \sM_g, \quad p^{-1}(M_I) = H^0(M_I, \Oh(k K))
\end{equation}
on the moduli space $\sM_g$ of curves of genus $g$.
(Recall that $\dim \sM_g= 3g-3$).


In particular, for $k=2$ the vector bundle $V_{2} = T^* \sM_g$ is
the cotangent bundle of the moduli space, that is, $V_{2}$
 is the bundle with spaces of quadratic differentials as fibers.
In particular the canonical class $K_{\sM_{g}}$ is the first Chern
class of the vector bundle $V_{2}$:
\begin{equation}
 c_1 (V_2) = K_{\sM_g}.
\end{equation}
Now according to Mumford  this canonical class has
the decomposition
\begin{equation}
K_{\sM_g} = 13 \Theta + \text{boundary divisors $\dots$}
\end{equation}
where $\Theta$ is the diviser of zero theta constants.

If the vector bundle $V_{2}$ admits a holomorphic projective flat
connection then by
 (1.31) the number 13 has to be divided by $3g-3$.
 This isn't the case thus $V_{2}$
 doesn't admit
any flat connection and an algebraic curve as a polarized
algebraic variety isn't a
 result of
a successful quantization procedure of the classical dynamical
system  $(\Si, \om)$ where $\om$ is the restriction of the
Fubini-Study form.

But the history isn't finished! We can find a "faithful functor".

\subsection{Jacobians}

May be the first faithful functor for algebraic curve is mapping
of
 the curve to the
polarized jacobian
$(J(\Si_I), \Theta )$. There exists a lot of constructions but
 for us it will be very convinient
to consider $J(\Si_I)$ as the {\it moduli space of holomophic
topological trivial line bundles}. Fixing a point $p_0 \in \Si_I$
we consider for every other point $p\in \Si_I$ the topological
trivial line bundle $\Oh(p - p_0)$, where we consider points as
divisors of degree 1. So we get the map
\begin{equation}
a_{p_0} \colon \Si_I \to J(\Si_I).
\end{equation}
By the Riemann theorem this map is an embedding (if genus $g>0$)
and induces the isomorphism
\begin{equation}
a \colon H_1(\Si, \Z) \to H_1(J(\Si_I), \Z).
\end{equation}
The moduli space of line bundles is a group (with respect to
the tensor product of line bundles). So, we can consider the Pontrjagin map
\begin{equation}
a_{p_0}^{g-1} \colon \Si_I^{[g-1]} \to J(\Si_I)
\end{equation}
sending a finite set of points $\{p_1, ... , p_{g-1}\}$ to
 the line bundle
\begin{equation}
\Oh(p_1 + ... + p_{g-1} - (g-1)p_0) \in J(\Si_I).
\end{equation}
The image of this map is a divisor $\Theta$ (so called
 {\it theta diviser}) and the
corresponding line bundle admits only one (up to $\C^*$)
holomorphic section (with
the zero set $\Theta$).

To reconstruct the curve $\Si_I$ from the pair
 $(\Theta \subset J(\Si_I))$ we have to consider
 the Gauss map
\begin{equation}
G \colon \Theta \to \PP T^* J(\Si_I)
\end{equation}
to the projectivization of the cotangent bundle sending
 non singular point of $\Theta$ to  the tangent hyperplane
  to $\Theta$ at this point.
  But the cotangent bundle (of a group) is
  trivial $\PP T^* J(\Si_I) = \PP^{g-1} \times  J(\Si_I)$, and
  the composition of the projection to $\PP^{g-1}$ and
   the Gauss map gives the rational cover
\begin{equation}
G \circ p \colon \Theta \to \PP^{g-1}.
\end{equation}
Now the  ramification diviser of this cover is the dual divisor
 to
the canonical embedding of $\Si_I$ into $(\PP^{g-1})^*$ (for
hyperelliptic curves we use other arguments). By the classical
projective duality theory we can reconstruct a curve by its dual
divisor. Thus the sending of a curve $\Si_I \to (J(\Si_I),
\Theta)$ is a faithful functor.

Now using our rule (1.97) we may quantize the polarized variety
 $(J(\Si_I), \Theta)$.

The role of wave function spaces now are played by the spaces of
 holomorphic sections:
\begin{equation}
\sH_{I}^k = H^0(J(\Si_I), \Oh (k \Theta)).
\end{equation}
These spaces are called spaces of theta function of
 the Riemann surface $\Si_I$ of level $k$.

This quantization is successful! Corresponding vector bundles
 on
the moduli space admit projective flat connections (see
\cite{W}). But the situation is much better than just a successful
quantization.

\subsection{Algebro geometrical theory of theta functions}.

Consider a principal polarized abelian variety $(A, \Theta)$ of
 dimension $g$. Then the line bundle $\Oh(k \Theta)$ admits
 the finite group $K$ of translations preserving this line bundle.
It is easy to see that
\begin{equation}
  K = A_k \subset A
\end{equation}
is the subgroup of points of order $k$ in the abelian group $A$.
 Recall that we have a symplectic form on $K$ given by the wedge
  product of 1-cycles coupled by the polarization $\Theta$.
The decomposition
\begin{equation}
  K = K_+ \times K_-
\end{equation}
by two isotropic subgroups of this group is called
theta-structure of level $k$.

The beautiful classical result reproduced by Mumford is the following
\begin{enumerate}
\item In the complete linear system $\vert  k \Theta \vert$ there
exists only one divisor $D_+$ which is invariant with respect to
the subgroup $K_+$ - action.
\item Every translation $a \in K_-$ defines a divisor $D_+ + a$
and the rational function $\theta_a$ with pole at the divisor
$D_+$ and zeros at the divisor $D_+ + a$. This function is called
 theta function with  characteristic $a \in K_-$.
\item The space of holomorphic sections
\begin{equation}
 H^0(A, \Oh (k \Theta)) = H^0(A, \Oh(D_+)) = \oplus_{a \in K_-} \C \theta_a.
\end{equation}
\item This is the theta functions with characteristics decomposition.
\end{enumerate}

Note, that the cardinalities of our finite abelian groups are
given by formulas
\begin{equation}
 \vert K \vert = k^{2g} , \quad \vert K_{\pm} \vert = k^{g}.
\end{equation}

Now it is quite useful to apply the extended Kodaira-Spenser
theory to this situation. First of all for every abelian variety
\begin{equation}
  H^0(A, T) \to H^2(A, \Oh)
\end{equation}
is an isomorphism thus the condition (1.48) holds. Moreover, from
cohomology sequence of (1.50) we have
\begin{enumerate}
\item an isomorphism
\begin{equation}
 H^0(A, S^2(T)) = H^1(A,D^1(k \Theta)));
\end{equation}
\item
\begin{equation}
  H^0(A, \Oh) = H^0(A, D^1(\Oh(k \Theta)));
\end{equation}
\item and from (1.60)
we have
\begin{equation}
 0 \to H^0(A, \Oh(k \Theta) )/ \C \cdot s \to H_h^1(d_s) \to
 H^1(A, D^1(\Oh(k \Theta))) \to 0;
\end{equation}
\item so under a deformation of the pair $(A, \Oh(k \Theta))$
for every section $s$ there exits  unique deformation of this
section which is {\it defined by the heat equation} (see the
standard presentation (2.49) below).
\end{enumerate}

\subsection{Symplectic-combinatorial theory of theta functions}

A principal polarization $\Theta$ on an abelian variety $A$ can
be done by
 some unimodular skew symmetrical integer form
\begin{equation}
\om  \in H^2(A, \Z)
\end{equation}
because  $A $ is a $2g$-dimensional real torus:
\begin{equation}
\om \in \wedge^2 H^1(A, \Z) = H^2(A,\Z).
\end{equation}
For a jacobian $J(\Si)$ this form is induced by the intersection
form on 1-cycles on the Riemann surface $\Si$.

The theory of unimodular integer skew symmetrical forms predicts
 that there are two isotropic $\Z$-sublattices in $H^1(A, \Z)$
 such that our form $\om$ has the standard symplectic shape.

Recall that the underlying smooth manifold of our abelian variety
is nothing else but a real torus $T^{2g}$ and the decomposition
 of  $H^1(A, \Z)$ by two isotropic $\Z$-submodules  induces
 the corresponding decomposition of this torus into two families of Lagrangian
 subtori
\begin{equation}
T^{2g} = T^g_+ \times T^g_-.
\end{equation}
Considering the projection to the second $g$-torus
\begin{equation}
\pi \colon T^{2g} \to T^g_-
\end{equation}
we obtain an integrable system that is a real polarization
of the torus $T^{2g} = A$.


\begin{figure}[tbn]
\centerline{\epsfxsize=3in\epsfbox{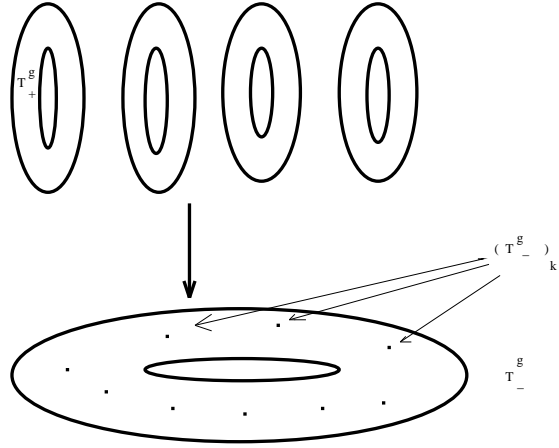}}
\caption{\sl Bohr-Sommerfeld fibers}
\label{Fig 2}
\end{figure}

The direct interpretation of the prequantization line bundle
shows that the subset of $k$-Bohr-Sommerfeld fibers
\begin{equation}
BS_k = (T^g_-)_k \subset T^g_-
\end{equation}
is the subgroup of points of order $k$ on our second torus
$T^g_-$  (see Fig.5). So, the wave function space of the
Bohr-Sommerfeld
 quantization is the sum of lines
\begin{equation}
Q_\pi^k (T^{2g}, \om) = \oplus_{a \in (T^g_-)_k} \C \cdot s_a
\end{equation}
where $s_a$ are covariant constant sections of restrictions of
 the prequantization line bundle.

We can see that the results of the Kahler quantization and the
Bohr-Sommerfeld quantization give wave function spaces of the
same rank so the quantization is numericaly perfect (see the end
of subsection 1.5).

This is the first amazing fact in this story.
We would like to say couple of words about the geometry behind
such coincidence.

Returning to the fibration (2.22) we can construct so called
{\it dual fibration}
\begin{equation}
\pi' \colon (T')^{2g} \to T^g_-
\end{equation}
changing every torus-fiber $T^g_+$ by the dual torus $(T^g_+)^*$.
This dual torus is nothing else but the space of gauge classes of
all $U(1)$-flat connections on the trivial line bundle on
$T^g_+$. We have got  a priori a new $2g$-torus with  new
symplectic form $\om'$ of the same type as before.

Now consider any holomorphic line bundle  $L$ on our old
complex torus $A$ with $U(1)$-connection which  curvature
form has to be proportional to $\om$. (Just like our
prequantization bundle $\Oh(k)$). Then for every point
$a \in T^g_+$ the restriction of $L$ to the Lagrangian fiber
 $\pi^{-1}(a)$ gives a flat line bundle, that is the point of the
 dual fibers
\begin{equation}
s_L \in (T^g_+)^* = (\pi^{-1}(a))^*.
\end{equation}
Thus every such line bundle defines a section
\begin{equation}
s_L \colon  T^g_+ \to (T')^{2g},
\end{equation}
that is a middle dimensional submanifold $s_L$ (more precisely
$s_L(T^g_+)$) in the manifold $(T')^{2g}$.

Of course  trivial connection defines  zero section of the
 new fibration (2.25) that is we have a middle dimensional
  submanifold
\begin{equation}
s_{\Oh(0)} \subset (T')^{2g}.
\end{equation}

Obviously, the intersection set
\begin{equation}
s_{\Oh(0)} \cap s_{\Oh(k)} = BS_k
\end{equation}
is the set of Borh-Sommerfeld fibers of level $k$. Moreover, it
is easy to see that both  our cycles are Lagrangian, oriented
 and admit positive index intersections only. Thus
the number of Bohr-Sommerfeld fibers is given by the
 intersection index of cohomology classes
\begin{equation}
[s_{\Oh(0)}] \cdot [ s_{\Oh(k)}] = \# BS_k.
\end{equation}
Now using the first Chern class of the polarization we can easily
compute these numbers.

The moral is the following: an abelian variety $A = T^{2g}$ has
 a symplectic partner $(T')^{2g}$
such that coherent sheaves geometry on $A$ can be encoded by the
geometry of Lagrangian cycles (or super cycles) on the partner
 $(T')^{2g}$. Usually this partner is called
 {\it mirror partner}.
  Mirror symmetry now is  far developed domain so we
  can't  dive to this subject (see \cite{CK}).
  In spite of the coincidence of our previous construction with
  the Strominger-Yau-Zaslow mirror symmetry construction we would
  like to stop here discussion of this subject but formally it can be extended
  along this line (see \cite{T3}).

 Of course the dimensional coicidence can be also extended. Moreover there are
  theta bases in both of wave function spaces.
  To compare them we have to construct
  a holomorphic object, namely, a section
  of the theta bundle of level $k$ using only
   Bohr-Sommerfeld fiber of the projection $\pi$ (2.22) (see Figure 1).

  We will do this using so called {\it coherent state
transform} or the slightly generalized
Segal-Bargmann isomorphism introduced in the
context of the quantum theory as a transform from the Hilbert
space of  square integrable functions on the configuration
space to the space of holomorphic functions on the phase space.

Considering  zero fiber $T^g_+$ of the polarization map $\pi$
(2.22) we can interpret our complex torus $A$ as a
complexification of this real torus and input our geometrical
situation to the situation of the classical coherent state
transform.

 In the finite dimensional context a configuration space is
 just $\R^g$ as the real part of the phase space $\C^g$.
  The  coherent state transform
\begin{equation}
CST_t \colon L^2(\R^g, d^g x) \to L^2(\C^g, d\mu_t) \cap \sH(\C^g)
\end{equation}
is a unitary isomorphism.(Here $d\mu_t$ is the Gaussian measure
and $\sH$ is the space of holomorphic functions.)

 For our abelian case we have to replace $\R^g$ by $T^g = U(1)^g$ and $\C^n$
 by our complex
torus $A$. Remark
 that $T^g$ in our case is a special Lagrangian
subtorus (see, for example \cite{T4}).

 B. C. Hall proposed a generalization of the CST where $\R^n$ is replaced
by $\SU(2)$ and $\C^n$ by $SL(2, \C)$ (see
 \cite{Ha}) and we will use it later for non-abelian case. But here we use this
 construction for  abelian group $U(1)^g \subset (\C^*)^g$.

 First of all, the decomposition (2.21) induces the decomposition
 \begin{equation}
H^1(A, \Z) = H^1(T^g_+, \Z) \oplus  H^1(T^g_-, \Z).
\end{equation}
Let us fix a basis $(a_1, \dots, a_g)$ in the lattice $H^1(T^g_+, \Z)$ and
a basis $(b_1, \dots, b_g)$ in the lattice
 $ H^1(T^g_-, \Z)$ such that the full system $(a_1, \dots, a_g ,b_1, \dots, b_g)$
is a standard basis of 1-homology of $A$. For a complex Riemann
surface $\Si$
 consider the period matrix
\begin{equation}
\Om = \Vert \Om_{ij}\Vert = Re \Om + i \cdot Im \Om
\end{equation}
  as a point of the Siegel space $H$.
  Let $\vec x =(x_1, \dots, x_g) \in [0,1]^g$
be periodic coordinates on $U(1)^{g} = T^g_+$. Then we have
 $U(1)$-invariant
Laplacian
\begin{equation}
\De^{Im \Om} = \sum_{i,j = 1}^g \frac{Im \Om_{ij}}{2 \pi} \frac{\p^2}{\p x_i
\p x_j}.
\end{equation}
The complexification of $T^g_+$ is $(\C^*)^g$ with coordinates $(e^{2\pi i z_1},
\dots,
e^{2\pi i z_g})$ for $\vec z = \vec x + i \vec y$. Then the Haar mesure on
 $(\C^*)^g$ is $d \vec x d \vec y$ and we have  $U(1)^g$-invariant
  complex
 Laplacian
\begin{equation}
\Delta^{Im \Om}_\C = \sum_{i,j = 1}^g \frac{Im \Om_{ij}}{2 \pi} (\frac{\p^2}{\p x_i
\p x_j} + \frac{\p^2}{\p y_i
\p y_j}).
\end{equation}

As the first step consider the fundamental solution $\mu_t$ at the
identity of the
 heat equation on $(\C^*)^g$
\begin{equation}
4 \frac{\p u}{\p t} = \Delta^{Im \Om}_\C  u
\end{equation}
and average the measure $\mu_t d \vec x d \vec y$ with respect to the action of
$U(1)^g$ :
\begin{equation}
\nu_t(\vec z) d \vec x d \vec y =
(\int_{U(1)^g} \mu_t(\vec z + \vec x) d \vec x) \quad d \vec x d \vec y.
\end{equation}
The explicit formula is
\begin{equation}
\nu_t(\vec z) = (\frac{2}{t})^{g/2} (det (Im \Om))^{-1/2} e^{\frac{\pi}{2t}
\sum_{i,j}
(z_i - \ov{z_i}) (Im \Om)^{-1}_{ij}(z_j - \ov{z_j})}.
\end{equation}
Now consider  non-self adjoint Laplace operator
\begin{equation}
\Delta^{-i \Om} = \sum_{i,j = 1}^g \frac{i}{2 \pi} \Om_{ij}\frac{\p^2}{\p x_i
\p x_j}.
\end{equation}
C. Florentino, J. Mourao and J. Nunes proved in \cite{FMN} the
following
\begin{prop} Let $\sC$ be the analytic continuation from $U(1)^g$ to
$(\C^*)^g$. Then the transform
\begin{equation}
C_t^{(-i \Om)} = \sC \circ e^{t/2 \Delta^{(-i \Om)}} \colon L^2(U(1)^g, d \vec x)
\to
\end{equation}
$$
\to L^2((\C^*)^g, d\nu_t) \cap \sH((\C^*)^g)
$$
is unitary.
\end{prop}
So for any function $f \in L^2(U(1)^g, d \vec x)$ given by its Fourier
 decomposition
\begin{equation}
f(\vec x) = \sum_{\vec n \in \Z^g} a_n e^{2\pi i \vec x \cdot \vec n}
\end{equation}
we have
\begin{equation}
(C_t^{(-i \Om)} f)(\vec z) =
\sum_{\vec n \in \Z^g} a_n e^{ti\pi \vec n \Om \vec n}\cdot
 e^{2\pi i \vec x \cdot \vec n}.
\end{equation}
This function is the analytic continuation to $(\C^*)^g$ of the solution of
the complex heat equation
\begin{equation}
2 \frac{\p u}{\p t} = \Delta^{(-i \Om)} u
\end{equation}
on $(\C^*)^g$ with the initial condition given by f.

But the coherent states transform (2.40) can be extended from
$L^2(U(1)^g, d \vec x)$ to the space of distributions $(C^\infty(U(1)^g))' $
given by Fourier series of the form (2.41) such that there exists an integer
$N > 0$ such that
\begin{equation}
lim_{\vec n \cdot \vec n \to \infty} \frac{| a_n |}{(1 + \vec n \cdot \vec n)^N}
= 0.
\end{equation}
The Laplace operator and its powers act as continuous linear operators on this space of distributions
(by duality from the corresponding action on $C^\infty(U(1)^g)$)
and define for ($t > 0$) the action of the operator $e^{t/2} \Delta^{(-i \Om)}$ on
distributions $f$ of the form (2.41) as
\begin{equation}
e^{t/2} \Delta^{(-i \Om)} (\sum_{\vec n \in \Z^g} a_n e^{2\pi i \vec x \vec n}) =
\end{equation}
$$
= \sum_{\vec n \in \Z^g} a_n e^{t i \pi \vec n \Om \vec n}
e^{2\pi i \vec z \vec n}.
$$

In the same paper \cite{FMN} the authors proved the following
\begin{prop} If a series $f \in (C^\infty(U(1)^g))' $ then the RHS of (2.45)
defines a holomorphic function on
$(\C^*)^g$.
\end{prop}
Now we are ready to map every Bohr-Sommerfeld fiber of the
projection (2.22) to an analytic function on $(\C^*)^g$. It is
enough to do it for  zero torus $T^g_+ = U(1)^g$ of this
projection.

For this let us consider the distribution
\begin{equation}
\theta^\R_0 (x) = \sum_{\vec n \in \Z^g} e^{2\pi i k \vec x \vec n}.
\end{equation}
It is nothing else but the delta-function at the identity of level
$k$. Let us apply  CST (2.40) to this distribution:
\begin{equation}
C_{1/k}^{(-i \Om)} (\theta_k^\R) =
\end{equation}
$$
= \sum_{\vec n \in \Z^g} e^{\pi i k \vec n \cdot \Om / k  \vec n k } \cdot e^{
2 \pi i k\vec n \vec z}.
$$
For every $l \in \Z^g / k \Z^g$ we can start with the function
\begin{equation}
\theta^\R_l (x) = \sum_{\vec n \in \Z^g} e^{2\pi i  \vec x (l+ k\vec n)}
\end{equation}
to get the analytic function
\begin{equation}
C_{1/k}^{(-i \Om)} (\theta_l^\R) = \theta_l(z, \Om) =
\end{equation}
$$
= \sum_{\vec n \in \Z^g} e^{\pi i (l+ k \vec n )\cdot \Om / k (l
+ k\vec n )} \cdot e^{ 2 \pi i (l + k\vec n ) \vec z}.
$$

 Remark that we are substitute  our continious
positive parameter $t$ coming from
 heat kernels by its discrete conterpart $1/k$.

We would like to emphasize that in spite of 150 years of
development of the classical theory of theta functions  this
"comparing quantization" approach was realized first quite
recently by C. Florentino, J. Mourao and J. Nunes in \cite{FMN}.

What we have to do now is just to compare our holomorphic functions on
$(\C^*)^g$ with holomorphic sections of the $\Theta$-line bundle on
$A$.

\subsection{Abelian holomorphic flat connections}

The usual way to define a line bundle $L$ is
  to consider some finite cover $\{U_i\}$ on the base and
  a collection of {\it transition functions} $f_{ij}$ which are
  regular
   and regular inversed on $U_i \cap U_j$.
On the intersection $U_i \cap U_j \cap U_k$ we have the cocycle
equality
\begin{equation}
f_{ij} \cdot f_{jk} \cdot  f_{ki} = 1.
\end{equation}
A collection $f_{ij}$ is equivalent to a collection $f'_{ij}$
iff there exists a collection $\{\beta_i\}$ of functions, where each
$\beta_i$ is regular and regular inversed on $U_i $ such that
\begin{equation}
f'_{ij} = \beta_i f_{ij} \beta_j^{-1}.
\end{equation}

The idea is to make these functions (2.51) as simple as possible using
this equivalence: for example take all $f'_{ij}$
 constant.

Starting with any collection of functions (2.50) let us
 consider the collection of  differentials
forms
\begin{equation}
\{  f_{ij}^{-1} d f_{ij} \} \in Z^{1}(\{U_{i} \} )
\end{equation}
Obviously, this cochain is a cocycle from $H^1(\Om) = H^{1,1}$
whose cohomology class
\begin{equation}
c(L) \in H^{1,1}
\end{equation}
is the first Chern class of $L$.

If this class is zero  then the cocycle is trivial and there are
matrices of differential forms $h_i$ such that
\begin{equation}
f_{ij}^{-1} d f_{ij} = h_i - h_j.
\end{equation}
 \begin{dfn} A collection of  differential forms $h_i$
  on $U_i$ is called  flat holomorphic connection of the
   vector bundle $L$ given by the collection of transition
    functions
$f_{ij}$.
\end{dfn}

Having such collection we can solve  the system of
linear equations on $U_i$
\begin{equation}
\beta_i^{-1} \cdot d \beta_i = h_i
\end{equation}
and obtain a new collection of transition functions (2.51).
Now using all equations we can see that functions
$f'_{ij} = \beta_i \al_{ij} \beta_j^{-1}$  are constant, that is
\begin{equation}
 d f'_ij = 0.
\end{equation}
If transition functions are constant then this vector bundle is
called
 {\it local system of coefficients} and we can trivialize it over
 any simply connected open set. Thus  we have a  character of
 the fundamental group of our base:
\begin{equation}
\rho \colon \pi_1 \to  \C^*.
\end{equation}

For a curve we've proved the classical result of Poincare:

\begin{prop}
Every $L \in J(\Si)$ admits  holomorphic flat
 connection
 given by a character $\chi \colon \pi_1(\Si) \to \C^*$.
\end{prop}

How many characters  give the same vector bundle?
 Just  as many as flat connections admitted by this bundle there are.
 The difference of any two holomorphic connections can be
 identified over
$U_i \cap U_j$, thus
\begin{equation}
\{h_i\} - \{ h'_i\} \in H^0(\Oh_\Si (K_\Si)).
\end{equation}

By Serre duality
\begin{equation}
H^1(\Oh)^* = H^0( \Oh( K_\Si))
\end{equation}
 is a fiber of
the cotangent bundle $T^* J(\Si)$.

The full space $A$ of all
 holomorphic flat connections is the space of all characters
\begin{equation}
A = (\C^*)^{2g}, \quad dim_\C A = 2g
\end{equation}

Every character defines a holomorphic vector bundle
 by the standard construction
\begin{equation}
L = U \times \C / (\pi_1,\rho)
\end{equation}
 where $U$ is the universal cover of our base with the natural
  action of the fundamental group of a base on $U$ and the $\rho$-action
  on $\C$ .

Thus we have  the forgetful map
\begin{equation}
f \colon A \to J(\Si).
\end{equation}
Every fiber of this map is the set of holomorphic flat connections
on a fixed line bundle. Thus the map $f$ sending each character
 to the corresponding line bundle provides on $A$ the structure of an
  {\it affine bundle over the cotangent bundle}.

So fibers of $f$ are affine spaces over $H^0(J(\Si),\Om)$ (but in
this case the cotangent bundle is the trivial $g$-dimensional
vector bundle). Over every small open set $U_i$ our affine bundle
admits a section that is can be identified with the restriction
of the cotangent bundle. Any descrepance  exists only  on the
intersections $U_i \cap U_j$. The collection of differences of
sections gives a 1-cocycle
\begin{equation}
\ep \in H^1(J(\Si), \Om).
\end{equation}
If this cocycle is trivial, our affine bundle coincides with its
vector bundle while non trivial affine bundle is defined by this
cocycle uniquely. By the Dolbeault theorem in our case the cocycle
$\ep_A \in H^{1,1} = H^2(J(\Si), \C)$. Poincare proved that
\begin{equation}
\ep_A = [\Theta] = [\om].
\end{equation}
To prove this equality we have to use the functional equation
for theta functions.

The affine bundle (2.62) admits {\it non holomorphic sections}:
let $U(1) \subset \C^*$ be the subgroup of unit norm numbers.
There is a map
\begin{equation}
mod \colon \C^* \to U(1), \quad z \to \frac{z}{|z|}.
\end{equation}
Then we have the subspace
\begin{equation}
U(1) (\C^*)^{2g}
\end{equation}
of unitary characters of $\pi_1(\Si)$ which is a {\it section} of
 the projection $f$ (2.62) of the  affine vector
 bundle $(\C^*)^{2g}$. It  follows immediately from the maximum
principle for a compact complex manifold.

So we have the smooth  identification
\begin{equation}
 J(\Si) = U(1)^{2g}.
\end{equation}

 Let $H_1(\Si, \Z) = <a_1, ... ,a_g, b_1,... ,b_g> $ be as before
  (see (2.32)).Such
  presentation defines
 very important subspace of the character space, so
 called  {\it abelian Schottky subspace}:
\begin{equation}
S^a_g = (\C^*)^g = \{ \chi \in (\C^*)^{2g} \vert \chi(a_i) = 1 \}.
\end{equation}
and
\begin{equation}
uS^a_g = (U(1))^g = \{ \chi \in U(1)^{2g} \vert \chi(a_i) = 1 \}
\end{equation}
 is
 {\it abelian unitary Schottky space}.

 The restrictions of the forgetful map $f$ (2.62) onto this space has the
 following properties
\begin{enumerate}
 \item $f \colon S^a_g = (\C^*)^g \to J(\Si)$ is infinite ($\Z^g$) cover,
\item $ f^* (\Oh (\Theta))$ is trivial, so holomorphic sections =
holomorphic functions = theta functions;
\item this line bundle is given by the {\it automorphy factors}
\begin{equation}
 e_{a_i}(z) = 1
\end{equation}
$$
e_{b_i}(z) = e^{-2\pi z_i - \pi i \Om_{ii}}
$$
where $\vec z$ are complex coordinates of the universal cover $U$ of $A$
described after formula (2.34), $\Om$ is
 the period matrix of $\Si$ and so on (see the set up around (2.32) - (2.35));
let $\{ u_i\}$ be the basis in $U$ dual to $\{ a_i\}$ and
 \begin{equation}
 u_{g + i} = \sum_1^g \Om_{ij} u_i;
\end{equation}
\item thus using these automorphy factors we see that the space $H^0(J(\Si),
\Oh(k \Theta))$ is naturaly identifyed with the space of holomorphic functions
$\theta$
 on $(\C^*)^g$ such that $\theta(z + u_i) = \theta$ and
\begin{equation}
 \theta(z + u_{g+i}) = e^{-2\pi k  z_i - \pi i k
 \Om_{ii} };
\end{equation}
\item thus any $\theta(z)$ admits the following decomposition
\begin{equation}
 \theta (z) = \sum_{l \in \Z^g / k \Z^g} \theta_l(z, \Om)
\end{equation}
where theta functions with characteristics are functions (2.49).
\end{enumerate}

\subsection{Perfect quantization}

Summarizing all results and constructions, we have

\begin{thm} Principal polarized abelian varieties
(and hence algebraic curves) admit the perfect
quantization by theta functions.
\end{thm}

In particular, we proved  very important fact (which we want to
generalize for the non-abelian case):

\begin{prop}
The  Bohr-Sommefeld wave space (2.24) doesn't depend on the
decomposition (2.21) determining  projection $\pi$ (2.22). This
decomposition defines a basis in this space only and the result
of Bohr-Sommerfeld quantization doesn't depend on the choice of a
real polarization of an "arithmetical" type (2.22).
\end{prop}
Of course this observation is obvious in the set-up of the classical theory
of theta functions, but in the set-up of non-abelian theta functions it is giving
an important "independence condition" in CQFT and low dimensional topology.

Let us list the sequence of tasks necessary to get a perfect quantization:
\begin{enumerate}
\item a $g$-torus $U(1)^g$ with the marked point $0 \in U(1)^g$ we identify with
the zero fiber of the fibration $\pi : A \to T^g_-$ (2.22) containing zero point;
\item for zero point $0 \in U(1)^g$ we construct $\de$-function $\theta_0^\R$
 as the Fourier
serie (2.46);
\item using the period matrix $\Om$ of $A$ we construct  CST
$C^{(-i \Om)}_{1/k}$ (2.40) , (2.41);
\item applying this transform to the distribution $\theta_0^\R$
we get the collection of holomorphic theta functions with characteristics $\{
\theta_l(z, \Om)\}$  (2.49) on the complexification  $(\C^*)^g$ of $U(1)^g$;
\item using the periods matrix $\Om$ we construct the cover $(\C^*)^g \to A$
 and get sections of the line bundle $\Oh(k \Theta)$ as functions subjecting to
 automorphity factors conditions (2.70);
 \item to check that  functions $\{
\theta_l(z, \Om)\}$  (2.49) correspond to sections of  $\Oh(k \Theta)$.
\end{enumerate}

This is the end of story for the moduli spaces of topological trivial
 vector bundles
 of rank one on  algebraic curves=Riemann surfaces. But any curve $\Si$ also
  has the other
 "moduli space", namely the moduli space of topological trivial
 vector bundles of rank 2. This is a faithful functor too. To implement
 this theory we have to plug  Conformal Field Theory in dimension
 $D = 2$. We will use the previous list of tasks as a pattern of
 much more sophisticated  efforts in this "non-abelian" case.

 Note that  $H^0(J(\Si), \Oh ( k \Theta))$ is the space of the irreducible
representation of the Heisenberg
 group of level $k$.

This statement is in the very heart of the classical theory of theta
functions. Actually the group $K$ (2.12) is the quotient of the
Heisenberg group $\Ga_k$ with the kernal  $\Z_k$ (with respect to the
multiplicative action). More precisely, the Heisenberg group acts
on the line bundle $\Oh(k \Theta)$ and on the space of its
holomorphic sections. This is an
irreducible representation of $\Ga_k$ and such representation is
unique up to projectiviazation. The existence of the special basis is provided
  by this
action. ( In particular, $rk H^0(J(\Si), \Oh ( k \Theta)) = k^g$).
This property (to be the space of   unique irreducible
representation of some group or algebra) is true for non-abelian
case too. But in this case the space of theta functions is the
space
 of irreducible representation of the gauge algebra
of WZW CQFT.

\section{Non-abelian theta functions}

The direct generalization of  jacobians of  algebraic curves
 as a faithful functor is the sending of an algebraic curve $\Si_I$
 to the {\it moduli space} $M^{ss}$ of semi-stable
  topologicaly trivial rk 2 vector bundles on this curve.
  Of course, we can consider much more complicated situation,
  but this case is the first {\it non-commutative} case which
   is expressive enough to see all new features of the
   geometrical situation.

Following  our pattern  we can quantize a Riemann surface $\Si$
using the faithful functor
\begin{equation}
 \Si \to M^{ss}(\Si)
\end{equation}
sending a Riemann surface to the moduli space of  semi-stable
vector bundles on it.

For this we have to repeat all steps of the abelian (classical) theory of
theta functions.

\subsection{Algebraic geometry of moduli spaces of vector bundles}

 The algebro-geometric part of the
non-abelian theory of theta functions is well developed by
Narashimhan, Beauville,
 Laszlo, Pauly, Oxbury, Ramanan, Sorger and many others.
In non-abelian context important new notion of  semi-stability
appears. Recall that for our case  semi-stability of $E$ means
that $E$ doesn't contain  linear subbundles of positive degree
(that is, with positive $c_1$). To work with rk 2 vector bundles
we need some information about the structure of coherent sheaves
on algebraic curves and more from the homological algebra. We
have to know that every coherent sheaf on a smooth algebraic
curve is a direct sum of a local free sheaf (= vector bundle) and
a torsion sheaf with support in a finite set of points. From the
homological algebra we need  Riemann-Roch theorem and first
cohomology of sheaves. All  these facts can be found in any
survey on this subject, for example in \cite{DSS}

Let $M^{ss} (\Si)$ be the moduli space of topologically trivial
 semi-stable holomorphic bundles of rank 2 on $\Si$ where every
non-stable (but semi stable) bundle is presented by the direct sum
of two opposite
 topological trivial line bundles.
Then  $ dim M^{ss}(\Si) = 3g-3$.

 The tangent space to the moduli space at a point $E$
has well known  shape
 \begin{equation}
TM^{ss}_E = H^1(\ad E)
\end{equation}
where $\ad E$ is the traceless part of the vector bundle of
endomorphisms of $E$. For any  stable bundle $E$ one has
$H^0(\Si_I, \ad E) = 0$ and thus by the Riemann-Roch theorem
\begin{equation}
dim M^{ss} = rk TM^{ss}_E = rk H^1(\ad E) = 3g-3.
\end{equation}

Definitions of theta divisors are absolutely parallel to (2.6):
let $\si \in Pic_{g-1}(\Si)$ be a line bundle on $\Si$ of degree
$g-1$; then we have

\begin{prop}
\begin{enumerate}
\item $\Theta = \{ E \in M^{ss}(\Si) \vert H^0(\Si, E(\si))
\neq 0\}$;
\item as a variety $\Theta$ is birationally equivalent
 to the projective space  $\PP^{3g-4}$;
\end{enumerate}
\end{prop}

Indeed, a general vector bundle from the divisor $\Theta$ twisted by $\si \in
Pic_{g-1}(\Si)$
 admits a section. (Recall that in VBAC-slang "twisted" means
  tensor product with the line bundle corresponding to $\si$.)
 We may expect that this section has no zeros, thus our general
  vector bundle $E$ can be represented as an extension
\begin{equation}
0 \to \Oh(-\si) \to E \to \Oh(\si) \to 0.
\end{equation}
Of course there exists  trivial extension $E = \Oh(-\si) \oplus
\Oh(\si)$ but this bundle isn't stable at all.

Any other extension is given by non zero cocycle of the space
\begin{equation}
Ext^1(\Oh(\si), \Oh(-\si)) = H^1(\Si_I, \Oh(-2\si)) = \C^{3g-3}
\end{equation}
by the Riemann-Roch theorem. Automorphisms of a subbundle act
 on these cocycles by  constant multiplication.
  Thus if $E(\si)$ has no other section (what happens in   general case),
  our vector bundle $E$ is defined uniquely by a point of
  the projective space $\PP H^1(\Si_I, \Oh(-2\si)) = \PP^{3g-4}$.
It is a simple exercise to estimate the dimension of spaces of
bundles admitting a section with  non trivial zero set.

The moduli space $M^{ss}$ is close to be rational itself. Let us
consider the space $Pic_g$ of classes of divisors of
 degree $g$ (instead of $g-1$ as in Proposition 5) and twist
  all bundles from $M^{ss}$ by some divisor $\si \in Pic_g(\Si)$.
  Then for every vector bundle $E \in M^{ss}$
we have
\begin{equation}
rk H^0(\Si_I, E (\si)) \geq 2.
\end{equation}
Now we have the canonical (evaluating) homomorphism
 \begin{equation}
 can \colon H^0(\Si_I, E(\si)) \otimes \Oh \to E(\si)
\end{equation}
Again by dimensional computations it is easy to see that this
moduli space contains a Zariski open set $M_0$ of bundles
 for which this homomorphism
can be extended to the exact sequence
\begin{equation}
0 \to  \C^2 \otimes \Oh \to E(\si) \to \oplus_{i=1}^{2g}
\Oh_{p_i} \to 0
\end{equation}
where $2g$ points $\{ p_i \}$ form the sum which is an effective
divisor equivalent to $2\si$:
\begin{equation}
p_1 + p_2 + .... + p_{2g} \in \vert 2\si \vert = \PP H^0(\Si_I,
\Oh(2\si)) = \PP^g.
\end{equation}
Why does this sum of points give an effective divisor from $\vert
2\si\vert$? Well, the homomorphism
\begin{equation}
 \wedge^2 can = det \colon \Oh \to \Oh(2 \si) = c_1(E(\si))
\end{equation}
and the effective divisor (3.9) is just  zero set of this
homomorphism.

Sending a vector bundle $E \in M_0$ to such effective divisor
 we
obtain  the map (of course algebraic)
\begin{equation}
det \colon M_0 \to \PP^g.
\end{equation}
The fiber can be obtained by the canonical homomorphism data:
\begin{equation}
 ker can_{p_i} \colon \C^2 \to E_{p_i}.
\end{equation}
So for every point $p_i$ we have a point on the projective line
\begin{equation}
\ga_i = \PP ker  can_{p_i} \in \PP \C^2 = \PP^1.
\end{equation}
Thus the fiber of the map (3.11) is given by a collection of
points
\begin{equation}
(\ga_1, .... , \ga_{2g}) \subset \PP^1
\end{equation}
up to a projective linear transformation.

We can see that any fiber of the map $det$ is rational too
 and $M^{ss}$ is an algebraic fibration with rational
 $g$-dimensional base and rational $2g-3$-dimensional fiber.

The same trick we can use for vector bundles of higher ranks. But
we would like to recall the old problem of the birational
geometry:

{\it Is the moduli space $M^{ss}$ rational?}

In spite of many attempts to solve  this rationality problem,
the answer isn't known up to now.

\begin{dfn}
The space $H^0(M^{ss}(\Si), \Oh ( k \Theta))$ is called the space
of non-abelian theta functions of level $k$.
\end{dfn}
Continuing the parallel description of the geometry of moduli
spaces we get
\begin{prop}
\begin{enumerate}
\item  Theta divisors are ample. Moreover
\item  $H^2(M^{ss}, \Z) = \Z$.
\item (Faltings 92) $H^0(M^{ss}(\Si), \Oh ( k \Theta))$ is the space
  of irreducible representation of the gauge algebra
of WZW CQFT.
\end{enumerate}
\end{prop}
Such spaces have played a noted role in Conformal Field Theory.
The Proposition  will be proved below.

\subsection{Holomorphic flat connections}

Up to now we have considered any vector bundle as a local free
sheaf.
 As far as for line bundles   any vector bundle is defined by
 some finite cover $\{U_i\}$ of the base and by
  a collection of {\it transition matrices} $\al_{ij}$ regular
   and regular inversed on $U_i \cap U_j$.
On the intersection $U_i \cap U_j \cap U_k$ we have again the
cocycle equality
\begin{equation}
\al_{ij} \cdot \al_{jk} \cdot  \al_{ki} = 1.
\end{equation}
A collection $\al_{ij}$ is equivalent to a collection $\al'_{ij}$
iff there exists a collection $\beta_i$ matrices, each $\{\beta_i\}$ is
regular and regular inversed on $U_i $ such that
\begin{equation}
\al'_{ij} = \beta_i \al_{ij} \beta_j^{-1}.
\end{equation}

For the matrix case "as simple as possible" representatives can be
\begin{enumerate}
\item triangle matrices (the theory of extensions, see (3.4));
\item constant matrices (the theory of flat holomorphic
connections).
\end{enumerate}
Again let us
 consider the collection of matrices of differential
forms
\begin{equation}
\{  \al_{ij}^{-1} d \al_{ij} \} \in Z^{1}(\{U_{i} \} )
\end{equation}
Obviously this cochain is a cocycle from $H^1(\Om (End E))$ which
cohomology class
\begin{equation}
c(E) \in H^{1}(\Om(End E))
\end{equation}
is the full Chern class of $E$.

The automorphisms group $Aut E$ acts on the space $H^{1}(\Om(End E
))$ preserving this class.

For example, on a curve by the Serre duality
\begin{equation}
 H^{1}(\Om(End E))^{*} = H^{0}(End E)
\end{equation}
and it is easy to see that

\begin{prop}
\begin{equation}
<c(E), \si> \neq 0 \implies \quad \si \quad \text{is an idempotent
in } \ H^0(End E)
\end{equation}
\end{prop}
This simple observation belonging to Atiyah gives plenty  results
about the structure of multidimensional bundles: the vector bundle
$End E$ splits as
\begin{equation}
End E = \Oh \oplus ad E
\end{equation}
where the first component $\Oh$ is the trivial line bundle
corresponding to homotheties  and $ad E$ corresponds to traceless
endomorphisms. Thus we have the decomposition
\begin{equation}
\Om(End E) = \Om \oplus ad E \otimes \Om
\end{equation}
where the first component sends  endomorphisms to their  traces.
Thus if $H^0(End E)$ contains one indempotent then it contains a
second too and such vector bundle is the direct sum of line
bundles.

 Vector bundle $E$ is called {\it indecomposable} if $E \neq E_{1}
\oplus E_{2}$.

For such bundles
\begin{equation}
c(E) \in H^{1}(\Om)
\end{equation}
since $H^0(End E)$ doesn't contain any idempotents.

If $\{ \al_{ij}^{-1} \cdot d \al_{ij} \}$ is the cochain (3.15)
defining $c(E)$ then
\begin{equation}
 tr (\al_{ij}^{-1} d \al_{ij}) \in H^{1} (\Om).
 \end{equation}
 By the Dolbeault isomorphism
\begin{equation}
H^1(\Om) = H^{1,1} \in H^2(\C).
\end{equation}
The formal matrix equality
\begin{equation}
tr (\al_{ij}^{-1} \cdot d \al_{ij}) = (det \Vert\al_{ij}
\Vert)^{-1} d (det \Vert \al_{ij} \Vert)
\end{equation}
shows that
\begin{equation}
c(E) = c(det E)
\end{equation}
for any indecomposable vector bundle on  curve.

It is easy to see (just by the direct interpretation) that this
cohomology class is integer and moreover
\begin{equation}
\{\al_{ij}^{-1} d \al_{ij}\} = c_1(E)
\end{equation}
where $E$ is the indecomposable vector bundles given by
$\{\al_{ij}\}$. If this class is zero as in our case of topologically
trivial bundles then the cocycle is trivial and there are
matrices of differential forms $h_i$ such that
\begin{equation}
\al_{ij}^{-1} d \al_{ij} = h_i - h_j.
\end{equation}
 \begin{dfn} Collection of matrix differential forms $h_i$
  on $U_i$ is called  flat holomorphic connection on the
   vector bundle $E$ given by the collection of transition
    functions
$\al_{ij}$.
\end{dfn}

Having such collection we can solve  the system of linear
equations on $U_i$
\begin{equation}
\beta_i^{-1} \cdot d \beta_i = h_i
\end{equation}
and  get a new collection of transition matrices. Now using all
equations we can see that matrices $\al'_{ij} = \beta_i \al_{ij}
\beta_j^{-1}$  are constant, that is
\begin{equation}
 d \al'_ij = 0.
\end{equation}
If matrices are constant,  this vector bundle is  called
 {\it local system of coefficients} and we can trivialize it over
 any simply connected open set and obtain the  representation of
 the fundamental group of the base:
\begin{equation}
\rho \colon \pi_1 \to SL(2, \C).
\end{equation}
Actually, because of the equivalence relations  we obtain a class of
 representations only. That is the representation $\rho$ up to
 the adjoint action of $SL(2, \C)$.

For  curve we have the classical result:

\begin{prop} (A. Weil) Every $E \in M^{ss}(\Si)$ admits a
holomorphic flat connection
 given by a class of representations $\rho \colon \pi_1(\Si) \to SL(2, \C)$.
\end{prop}
The author emphasizes   again that he've proved this statement
already.

\begin{rmk}
In the influential paper \cite{We} A. Weil proposed
\begin{enumerate}
\item the proof of the previous Proposition;
\item the idea that high dimensional vector bundles should play
the role of non-abelian analogue of Jacobians;
\item the notion of  {\it class of matrix divisor} which is
equivalent to the notion of  vector bundle.
\end{enumerate}
\end{rmk}

Many classes of representations give the same vector
bundle:
 the difference of two holomorphic connections can be
 identified over
$U_i \cap U_j$, thus
\begin{equation}
\{h_i\} - \{ h'_i\} \in H^0(End E \otimes \Om).
\end{equation}
 An element of the last space is a homomorphism
\begin{equation}
\phi \colon E \to E \otimes \Om
\end{equation}
and for   case of curves (when $K_\Si = \Om_\Si$) it is called
{Higgs field}.

Note that by Serre duality
\begin{equation}
H^1(\ad E)^* = H^0( \ad E \otimes K)
\end{equation}
that is the vector bundle over the moduli space $M^{ss}$ with
fibers equal to spaces of Higgs fields is nothing else but the
cotangent bundle $\Om$ of $M^{ss}$.

For rank 2 vector bundles with fixed determinant we have to
consider traceless Higgs fields only:
\begin{equation}
\phi \in H^0(\Si, \ad E \otimes K).
\end{equation}

The full space $A^{na}$ of all
 holomorphic flat connections
\begin{equation}
A^{na} = CLRep(\pi_1(\Si), SL(2, \C))
\end{equation}
is the space of all classes of representations and
 every representation defines a holomorphic vector bundle
just by the standard construction
\begin{equation}
E  = U \times \C^2 / (\pi_1,\rho)
\end{equation}
 where $U$ is the universal cover of $\Si$ with the natural
  action of the fundamental group of athe base.

Of course in the set of such bundles there are non-stable but
indecomposable
 vector bundles and even semi-stable but indecomposable. Let us remove
 such representations and get the space $A_{na}^{ss}$ of classes
 of representations which gives stable and decomposable  semi-stable bundles.

Consider  the forgetful  map
\begin{equation}
 f \colon A^{na}_{ss} = CLRep(\pi_1(\Si), \SU(2)), SL(2, \C)) \to M^{ss}.
\end{equation}
This is the affine bundle over the cotangent bundle (see (2.62
)). Again such bundle is given by a cocycle
\begin{equation}
\ep_{na} \in H^1(M^{ss}, \Om) = H^{1,1}(M^{ss}, \C)
\end{equation}

\begin{prop}
\begin{equation}
\ep_{na} = [\Theta_{na}] = [\om_{na}]
\end{equation}
\end{prop}
That is, the answer is precisely the same as in the abelian case.
We will prove this fact a little bit latter using special
subspaces.

Both of affine bundles (abelian and non-abelian) (2.62) and (3.39)
admit {\it non holomorphic sections}. Again constructions of these
sections are quite parallel (see (2.67) for the abelian case) :
the including $\SU(2) \subset SL(2, \C)$ of complexification
defines
 the subspace
\begin{equation}
CLRep (\pi_1(\Si), \SU(2)) \subset CLRep (\pi_1(\Si), SL(2, \C)
\end{equation}
which is the section of the projection $f$ (3.39). (Again it
follows from the maximum principle for the compact curve.) But the
new feature of non abelian
 case is the following:

 we don't know if the restriction of forgetfull map to
 $CLRep (\pi_1(\Si), \SU(2))$ is
 onto. Does every semi-stable bundle admits Hermitian flat connection?

 We will prove this statement  later, and note that the restriction map defines
the Narasimhan-Sesadri isomorphism
\begin{equation}
 NS  \colon  CLRep(\pi_1(\Si), \SU(2)) \to M^{ss}.
\end{equation}

So we have differential identification
$$
  CLRep(\pi_1(\Si), \SU(2)) = M^{ss}.
$$

 Let $\pi_1 = <a_1, ... ,a_g, b_1,... ,b_g \vert \prod_{i=1}^g [a_i, b_i] = 1>$
 be a standard presentation
of the fundamental group. Such presentation defines a couple of
 very important subspaces of the spaces of representations, so
 called Schottky subspaces:
\begin{enumerate}
\item $S_g = \{ \rho \in A^{na} \vert \rho (a_i) = 1 \}$ is
non-abelian complex Schottky space.
\item $S_g = SL(2, \C)^g / Ad_{diag}SL(2, \C)$
\item $dim S_g = dim M^{ss}$
\item $uS_g = S_g \cap CLRep(\pi_1(\Si), \SU(2))$ is  unitary Schottki space;
\item $dim \quad uS_g = 1/2 dim M^{ss}$.
\end{enumerate}

We saw in the abelian case that the restriction of the forgetful
map to the Schottki space is a cover giving theta sections as
functions.

For the non-abelian case the situation is much more complicated:
\begin{enumerate}
 \item the forgetful map
 \begin{equation}
 f \colon S_g \to M^{ss}(\Si)
 \end{equation}
  is meromorphic only
 (because of existence of non stable flat bundles).
  For example, the Schottky uniformization of  curve  gives the
  non stable bundles.
\item Properties of this map are unknown yet.
C. Florentino proved only that the differential of this map at the
unitary Schottky space is non degenerated and in section 9 we show that for
 general curve $\Si$ the image of $f$ (3.44) is Zariski dense in $M^{ss}(\Si)$.
\end{enumerate}
 To avoid this gap, we have to use Symplectic and Lagrangian
 geometry and sophisticated analysis.

\subsection{Moduli of stable pairs
and the desingularization of moduli spaces}

The construction of the moduli spaces was achived in the 1960's,
mainly by D.Mumford and the mathematicians of  Tata institute.
But after appearing of the gauge theory approach to vector bundles
these moduli spaces were input in the collection of moduli spaces
of {\it stable pairs} related by a chain of {\it flips} (see
\cite{Re}). This chain of flips was used by Thaddeus, Bertram and
Bredlow-Daskalopoulos for the computations of the cohomology ring
of moduli spaces and by Seshadri \cite{Se} for the construction
of a desingularisation of moduli spaces.

New object to study will be a pair $(E, s)$ consisting of a
 vector bundle $E$ on a curve $\Si$ and a nonzero section
$s \in H^{0}(E)$. S. Bredlow defined  new stability condition for
such pairs and proved  theorem relating stable pairs to vortices
on the Riemann surface. The vortex equation depends on  positive
real parameter $t$ and so the stability condition also depends on
$t$.
\begin{rmk}
The vortex equation requests a Kahler metric, not a conformal
class only. We consider the algebro-geometrical interpretation
 here. We fix the normalization of such metric such that the
 volume $vol \Si = 4 \pi$.
\end{rmk}

Our vector bundles admit nonzero sections so we have to consider
some positive line bundle $D$ of degree $d$ such that $det E =
D$. We obtain this case by twisting vector bundles from $M^{ss}$
by an appropriate line bundle. The stability condition for  pair
$(E, s)$ with rk 2
 vector bundle is the
following: a pair is semi-stable if for any line bundle $ L
\subset E$
\begin{equation}
deg L \leq 1/2 d - t \quad \text{ if } s \in H^{0}(L)
\end{equation}
$$
deg L \leq 1/2 d +t \quad \text{ otherwise. }
$$
The fondation of this theory can be found in \cite{BD},
\cite{Ber}, \cite{Th2}.

Let $M(D, t)$ be the moduli space of stable pairs. It's easy to
see that
\begin{enumerate}
\item $M(D, t)$ is nonempty if and only if $t  \leq d/2$;
\item if $t$ is irrational then $M(D, t)$ is compact;
\item if for $(E_{i}, s_{i}) \in M(D, t), \quad i=1,2$ and if there
exists a homomorphism $\phi \colon E_{1} \to E_{2}$ such that
$\phi(s_{1}) = s_{2}$ then $\phi$ is an isomorphism;
\item for $(E, s) \in M(D, t)$ there are no endomorphisms of $E$
 annihilating the section $s$ exept 0 and no endomorphisms preserving $s$
exept  identity;
\item if $(E, s) \in M(D, t)$ then $(E(\xi), s(\xi)) \in M(D(\xi), t)$
for any effective divisor $\xi$.
\end{enumerate}

The deformation theory of stable pairs is very simple: let $t$ be
irrational, then  for $(E, s) \in M(D, t)$
\begin{enumerate}
\item there is a natural exact sequence
\begin{equation}
o \to H^{0}(ad E) \to H^{0}(E) \to TM(D,t)_{(E,s)} \to H^{1}(ad
E) \to H^{1}(E);
\end{equation}
where homomorphisms are defined by the section $s$ .
\item $rk
TM(D,t)_{(E,s)} = d - g -2.$
\end{enumerate}
So for irrational $t$ the space $M(D, t)$ is smooth and we
concentrate our attention on these smooth moduli spaces. More
precisely we have the case if
\begin{equation}
t_{i} \in (max(0, d/2 - i - 1, d/2 - i)).
\end{equation}
For all $t_{i}$ from these intervals we have the same compact
smooth moduli spaces $M(D, t_{i})$ related by the chain of flips
(see \cite{Re}). We will use them later.

In general case the stability condition for  pair $(E, s)$ is the
following: a pair is semi-stable if for any proper subbundle $
E_{1} \subset E$
\begin{equation}
\mu (E_{1}) \leq \mu (E) - t \quad \text{ if } s \in H^{0}(E_{1})
\end{equation}
$$
\mu E_{1}  \leq \mu(E) +t \quad \text{ otherwise. }
$$
where $\mu$ as usual is  slope of the bundle. The fondations of
this theory can be find in \cite{BD}.

For  moduli spaces of stable pairs we have the same collections
of statements as before.

To apply this method to desingularization of $M^{ss}$ we have to
send a vector bundle $E \in M^{ss}$  to the pair $(End E, s)$
where $s \in H^{0}(End E)$. For a stable $E$ the choice of
section is unique: it has to be the line subbundle  of
homotheties. For
 semi-stable bundles there are two cases: $H^{0}(End E) = 2, 4$.
 Considering the parameter $t$ which is very close to $0$ we get
 the moduli space of stable pairs $\widetilde{M^{ss}}$ with the
 birational projection to $M^{ss}$ which is a smooth variety.
 Always below {\it we will consider this canonical desingularization
will not mention  that specially}.

\begin{rmk}
This desingularization was constructed by Seshadri \cite{Se}. He
considered a {\it parabolic structure} on  vector bundle over
fixed point of a curve and choose  special parameters of such
structure to get the same non singular variety with a birational
map to the $M^{ss}$.
\end{rmk}

\subsection{Holomorpic symplectic geometry of Higgs fields.}

Returning to Higgs fields (3.36),
 we can send every pair $(E, \phi)$ where $\phi$ is a Higgs
 field
to the space of holomorphic quadratic differentials:
\begin{equation}
\wedge^2 \phi \colon \Oh \to \wedge^2 (E \otimes K_{\Si}).
\end{equation}
if $E$ is stable our pair is stable too and  $T^* M^{ss}$ is a
moduli space of such pairs. Then
 we get
the holomorphic map
\begin{equation}
\pi \colon T^* M^{ss} \to H^0(\Si, K_{\Si_I}^2).
\end{equation}

As usual the cotangent bundle admits the holomorphic "action
1-form" $\al$ such that the   holomorphic differential of it $d
\al = \Om$ is the holomorphic symplectic form defining the
holomorphic skew symmetrical isomorphism
\begin{equation}
\Om \colon T T^* \to T^* T^*.
\end{equation}
One can see that we just reproduce notions and constructions of
the classical symplectic geometry in the holomorphic set up. In
particular, we can expect  existence of "completely integrable
systems" that is a holomorphic map of  holomorphic symplectic
variety $\sM$
\begin{equation}
\pi \colon \sM \to B
\end{equation}
such that every fiber is Lagrangian with respect to $\Om$ and the
generic fiber is a middle dimensional complex torus or an abelian
variety.We call such map as  holomorphic moment map.

The brilliant observation of Hitchin is that the map (3.50) gives
such  integrable system. Moreover, the variety $T^* M^{ss}$ can be
slightly extended as a holomphic symplectic variety such a way
that fiber of $\pi$ become compact.

\begin{dfn}  Higgs pair $(E, \phi)$ is called stable if it doesn't
 exist a line
subbundle $L \in E$ of positive degree invariant with respect to
$\phi$. That is, the restriction of $\phi$ to $L$ hasn't the form
\begin{equation}
\phi \vert_L \colon L \to L \otimes K_\Si \subset E \otimes K_\Si.
\end{equation}
\end{dfn}

It can be shown (see \cite{H2}) that the moduli space $MH_\Si$ of
stable Higgs pairs
\begin{enumerate}
\item contains the cotangent bundle $T^* M^{ss}$ as a Zariski
 dense subset;
\item $T^* M^{ss}$ is different from $MH_\Si$ in codimension 2
and we can use  Hartog's principle;
\item the holomorphic symplectic form can be extended to such
form (with the same notation) $\Om$;
\item the map $\pi$ (3.50) can be extended on $MH_\Si$ to
 a holomorphic
 moment map
 \begin{equation}
\pi \colon MH_\Si \to H^0(\Si, \Oh(2K_\Si))
 \end{equation}
 with compact fibers.
 \end{enumerate}

All  these statements can be proved directly and we just
 describe  general
fiber of the map $\pi$. General Higgs field we can consider as an
endomorphism of the
 projectivization
 of our vector bundle $E$
 $$
\phi \colon \PP E \to \PP E = \PP E \otimes K_\Si.
$$
If such map is birational then over every point $p \in \Si$
\begin{enumerate}
\item either $\phi_p$ is a linear automorphism of $\PP E_p$ and
has two fixed points $p_1, p_2 \in \PP E_p$,
\item or the endomorphism $\phi_p$ is degenerated and the
image gives us one point $p \in \PP E_p$.
\end{enumerate}.
Thus we have  double cover
\begin{equation}
\phi \colon \Si_\phi \to \Si
\end{equation}
with  ramification  divisor
\begin{equation}
\xi_\phi = (\wedge^2 \phi)_0 \in \vert 2K_\Si \vert
\end{equation}
(recall that all our endomorphisms are traceless).

The curve $\Si_\phi$ is called  spectral curve. It lies on the
ruled surface $\PP E$ and we have the standard adjoint sequence of
sheaves on this surface
\begin{equation}
0 \to \Oh(- \Si_{\phi}) \to \Oh \to \Oh_{\Si_{\phi}} \to 0.
\end{equation}

Multiplying this sequence by the Grothendieck line bundle $H$ on
$\PP E$ (recall
 that for the projection to base $R^0pr(H) = E, R^1pr (H) = 0 $) we get the
  exact triple
  \begin{equation}
0 \to \Oh(- \Si_\phi) \otimes H \to H  \to \Oh_{\Si_\phi}(H) \to
0.
\end{equation}
The exact sequence of this triple for direct images gives us the
 isomorphism
\begin{equation}
0 \to E \to R^0(\Oh_{\Si_{\phi}} (H)) \to 0
\end{equation}
because of the restriction
$$
\Oh(-\Si_phi) (H )\vert_{\PP E_p} = \Oh(-2) \otimes \Oh(1) =
\Oh(-1)
$$
on the projective line.

So  on the spectral curve $\Si_\phi$ our vector bundle $E$ defines
the line bundle $\Oh_{\Si_{phi}} (H) \in Pic (\Si_\phi)$.

For  double cover $\phi$ (3.55) the norm map
\begin{equation}
Nm \colon J(\Si_\phi) \to J(\Si)
\end{equation}
has connected preimage of zero
\begin{equation}
Nm^{-1}(0) = Pr_\phi
\end{equation}
which is an abelian variety which is called  {\it Prym variety}.

The isomorphism (3.59) shows that
\begin{equation}
Nm (\Oh_{\Si_\phi} (H)) = \wedge^2 E = \Oh.
\end{equation}
Thus
\begin{equation}
\Oh_{\Si_\phi} (H) \in Pr_\phi.
\end{equation}

So,  general quadratic differential
 $w \in H^0(\Si, \Oh(2 K_\Si)$ has
 zero-set
\begin{equation}
(w)_0 = p_1, \dots , p_{4g-4}
\end{equation}
and defines the double cover
\begin{equation}
\phi \colon \Si_w \to \Si
\end{equation}
with the  ramification locus (3.64).

It is easy to see that every $L \in Pr_\phi$ defines rk2 vector
 bundle on $\Si$:
\begin{equation}
E_L = R^0 L
\end{equation}
with $\wedge^2 E_L = \Oh.$ This bundle is stable as a Higgs
bundle if the cover is irreducible.

Thus the fiber of the map $\pi$ (3.54) is
\begin{equation}
\pi^{-1} (w) = Pr_\phi.
\end{equation}

\subsection{Gauge theory on  Riemann surface. Higgs fields.}

 The new features of the gauge theory is a necessity to
 consider objects  of different dimensions.
 For example,
Riemann surface bounds  threefold (handlebody) and we come to
3-dimensional theory to get a result in 2-dimensional case. In
the same vein sometimes it is reasonable to consider graphs
(1-dimensional complexes)
 instead of  Riemann
surfaces and so on.

A classical field theory on a manifold $M$ has three
ingredients:
 \begin{enumerate}
 \item a collection $\sA$ of {\it fields} on $M$, which are
 geometric
objects, such as sections of vector bundles, connections on
 vector bundles,
maps from $M$ to some auxiliary manifold (the target space) and
so on;
\item an {\it action functional}
$ S \colon \sA\to\C$ which is an integral of a function $L$ (the
Lagrangian) of fields;
\item a collection of observable functionals on the space of
 fields,
$ w \colon \sA\to\C.$
 \end{enumerate}
For example: {\it  Chern--Simons functional}. Here
\begin{enumerate}
\item
  $M$ is a
3-manifold, \item the collection of fields
 \[
 \sA=\Om^1(M) \tensor \fsu (2)
 \]
is the space of 1-form with coefficients in the Lie algebra,
\item the action functional
 \begin{equation}
S(a)=\frac{1}{8} \pi^2\int_M tr(a \cdot d a + \frac{2}{3} a^3),
 \end{equation}
\item as an observable  we can consider a {\it Wilson loop},
given by some knot $K\subset M$:
 \[
 W_K(a)= tr(Hol_K(a))
 \]
- the trace of  holonomy of  connection $a$ around  knot $K$.
\end{enumerate}
\begin{rmk} These functions on the space of fields $\sA$ form
an algebra because of the Mandelstam constraint:
$$
tr(A ) \cdot tr(B) = tr(A \cdot B) + tr( A^{-1} \cdot B).
$$
\end{rmk}
Our space $CLRep(\pi_1(\Si), \SU(2))$ consists  of classes of
$\SU(2)$-representations of the fundamental group of a Riemann
surface of genus $g$. This space contains  the subspace of
reducible representations
 \begin{equation}
 CLRep(\pi_1(\Si), \SU(2))^{red} \subset CLRep(\pi_1(\Si),\SU(2)).
 \end{equation}
To get this space as the phase space of some mechanical system,
 consider a
compact smooth Riemann surface $\Si$ of genus $g > 1$ and trivial
Hermitian vector bundle $E_h$ of rank 2 on it.
 As usual, let $\sA_h$ be the
affine space (over the vector space $\Om^1(\End E_h)$)
 of Hermitian
connections and $\sG_h$ is the corresponding Hermitian gauge
group. This space has the subspace:
 \[
 \sA_h^{red} \subset \sA_h
 \]
where the left-hand side is the subset of reducible connections.
 As usual,
let
 \[
\sA_h^*=\sA_h-\sA_h^{red}.
 \]

The sending of  connection to its curvature tensor defines a
$\sG_h$-equivariant map
 \begin{equation}
F \colon \sA(E_h)\to\Om^2(End E_h)=Lie(\sG_h)^*
 \end{equation}
to the coalgebra Lie of the gauge group.

We can consider this map as the moment map with respect to the
action of $\sG_h$. The subset
 \begin{equation}
F^{-1}(0)=\sA_F
 \end{equation}
is the subset of flat connections.

For a connection $a\in \sA_F$  a tangent vector to $\sA_h$ at $a$ is
 \begin{equation}
 \om\in\Om^1(\End E_h)=T \sA_h,
 \end{equation}
and for it the following condition holds
 \begin{equation}
  \na_a (\om)=0.
 \end{equation}
The trivial vector bundle $E_h$ admits  zero connection
 $\theta $,
which is interesting and important from many points of view. In
particular it provides possibility to identify $\sA_h$ with
$\Om^1(End E_h)$  sending  connection $a$ to the
 form $a - \theta$. We
will identify forms and connections  this way.

 The space $\sA_h=\Om^1(End E_h)$ is the {\it collection of fields} of
YM-QFT with the {\it Yang--Mills functional}
 \[
 S(a)=\int_{\Si} |F_a|^2.
 \]
Thus $\sA_F$ is a {\it classical phase space}, that is, the space
of solutions of the {\it Euler--Lagrange equation} $\de S(a)=0$.

There exists a symplectic structure on the affine space $\sA_h$,
induced by the canonical 2-form given on the tangent space
$\Om^1(End E_h)$ at a connection $a$ by the formula
 \begin{equation}
 \om (\ga_1,\ga_2)=\int_{\Si} tr (\ga_1 \wedge\ga_2).
 \end{equation}
This form is $\sG_h$-invariant, and its restriction to
$\sA_F^{irr}$ is degenerate along $\sG_h$-orbits: at a connection
$a$, for a tangent vector $\ga \in\Om^1(End E_h)$, we have
 \[
 \om\in T\sG_h \text{if and only if}  \om=\na_a \phi \quad\text{for $
 \phi\in\Om^0(\End E_h)=Lie(\sG_h)^*$,}
 \]
and
 \[
 \int_{\Si} tr (\na_a \fie \wedge\om)=\int_{\Si}tr(\fie\wedge\na_a\om)=0.
 \]
Hence
 \begin{equation}
 \om\in T \sA_F \text{if and only if} \Om_0 (\na_a \fie,\om)=0.
 \end{equation}

Interpreting (3.70) as a moment map and using symplectic reduction
arguments, we get a non degenerate closed symplectic form $\om$ on
 the space
 \begin{equation}
 \sA_F /\sG_h= CLRep(\pi_1(\Si), \SU(2))
\end{equation}
of classes of $\SU(2)$-representations of the fundamental group
of the Riemann surface, and a stratification of this space. The
form $\om$ defines  symplectic structure on $CLRep(\pi_1(\Si),
\SU(2))^{irr})$ and  symplectic orbifold structure on
$CLRep(\pi_1(\Si), \SU(2))$.

\begin{rmk}
The symplectic structure $\om$ is defined in  pure topological
way (see \cite{G1}).
\end{rmk}

On the other hand, the form $\om_0$ on $\sA_h$ is the differential
 of the
1-form $D$ given by the formula
 \begin{equation}
D(\om)=\int_{\Si} tr((a) \wedge\om).
 \end{equation}
We consider this form as a unitary connection $A_0$ on the trivial
principal \hbox{$U(1)$-bundle} $L_0$ on $\sA_h$.

To descend this Hermitian bundle and its connection to the
 orbit space,
one defines  $\Theta$-cocycle (or $\Theta$-torsor) on the
 trivial line
bundle (see \cite{RSW}). This cocycle is  $U(1)$-valued function
$\Theta$ on $\sA_h\times\sG_h$ defined as follows: for any triple
$(\Si,a,g)$ where $(a,g)\in\sA_h\times\sG_h$, we can find a
triple $(Y,A,G)$ where $Y$ is a smooth compact 3-manifold, $A$ a
$\SU(2)$-connection on the trivial vector bundle $\sE$ on $Y$ and
$G$ is a gauge transformation of it, such that
 \begin{equation}
 \p Y=\Si, \quad a=A \vert_{\Si} \quad \text{and} \quad g=G\vert_{\Si}.
 \end{equation}
Then
 \begin{equation}
 \Theta (a, g)=e^{i(CS(A^G)-CS(A))}.
 \end{equation}

Recall that the Chern--Simons functional on the space $\sA
(\sE_h)$ of unitary connections on the trivial vector bundle is
given by the formula
 (3.68).
It can be checked that the function (3.79) does not depend
 on the
choice of the triple $(Y,A,G)$ (see \cite{RSW}, \S2).

The differential of $\Theta$ at $(a, g)$ is given by the formula
 \begin{multline}
 d\Theta(\om,\phi)= \\
 \frac{\pi i}{4}\,\Theta\int_{\Si}\Bigl(tr(g^{-1} d g\wedge g^{-1}\om
 g)-tr(a\wedge\na_{a^g}\phi)+2tr(F_{a^g} \wedge \phi) \Bigr),
 \end{multline}
where $\om\in\Om^1(End E_h)$ and $\phi\in\Om^0(\End
E_h)=Lie(\sG_h)$.

But the restriction of this differential to the subspace of flat
connections is much simpler:
 \begin{equation}
 d \Theta (\om, \phi)=\frac{\pi i}{4}\,\Theta\int_{\Si} tr(g^{-1}d g
 \wedge g^{-1}\om g),
 \end{equation}
and is independent of the second coordinate.

That this function is  a cocycle ndeed results from the functional
equation
 \begin{equation}
\Theta (a, g_1 g_2)=\Theta (a, g_1) \Theta (a^{g_1}, g_2).
 \end{equation}
Using this function as a torsor $\sA_h \times_{\Theta}U(1)$, we
get a principal $U(1)$-bundle $S^1(L)$ on the orbit space
$\sA_h/\sG_h$:
 \begin{equation}
S^1(L)=(\sA_h\times S^1)/\sG_h,
 \end{equation}
where the gauge group $\sG_h$ acts by
 \begin{equation}
g(a, z)=(a^g, \Theta(a,g) z),
 \end{equation}
on the line bundle $L$ with a Hermitian structure.

Following \cite{RSW}, let us restrict this bundle to the subspace
of flat connections $\sA_F$. Then one can check that the
restriction of the form $D$ (3.77) to $\sA_F$ defines a
$U(1)$-connection $A_{CS}$ on the line bundle $L$.

By definition, the curvature form of this connection is
 \begin{equation}
F_{A_{CS}}=i\cdot\Om.
 \end{equation}

Thus the quadruple
 \begin{equation}
(CLRep(\pi_1(\Si), \SU(2)),\om, L, A_{CS})
 \end{equation}
is a {\it prequantum system} and we are ready to perform the
Geometric Quantization Procedure described in section 1.

\subsection{Complex polarization of $CLRep(\pi_1(\Si), \SU(2))$}

The standard way of getting a complex polarization is to give to
Riemann surface $\Si$ of genus $g$ a conformal structure $I$.
 We get a complex
structure on the space of classes of representations
$CLRep(\pi_1(\Si), \SU(2))$ as follows: let $E$ be our complex vector
bundle and $\sA$ the space of all
 connections on
it. Every connection $a\in\sA$ defines the correspoding covariant
derivative $\na_a\colon\Ga(E)\to\Ga(E\tensor T^*X)$, a first
order differential operator with the ordinary derivative $d$ as
the principal symbol. A complex structure gives the decomposition
$\p+\op$, so any covariant derivative can be decomposed as
$\na_a= \nabla^{1,0} + \na^{0,1}$,
 where
$\na^{1,0}\colon\Ga(E)\to\Ga(E\tensor\Om^{1,0})$ and $\na^{0,1}
\colon\Ga(E)\to\Ga(E\tensor\Om^{0,1})$. Thus the space of
connections admits a decomposition
 \begin{equation}
 \sA=\sA'\times \sA'',
 \end{equation}
where $\sA'$ is an affine space over $\Om^{1,0}(ad E)$ and $\sA''$
is an affine space over $\Om^{0,1}(ad E)$.

The group $\sG$ of all automorphisms of $E$ acts as the group of
gauge transformations, and the projection $pr \colon \sA\to\sA''$
to
 the space
$\sA''$ of $\op$-operators on $E$ is equivariant with respect to
the $\sG$-action.

Giving $E$ a Hermitian structure $h$, we get the subspace
$\sA_h\subset \sA$ of Hermitian connections, and the restriction
of the
 projection $pr$
to $\sA_h$ is one-to-one. Under this Hermitian metric $h$, every
element $g\in\sG$ gives an element $\overline{g} = (g^*)^{-1}$
such that
 \[
 \overline{g} = g \text{if and only if}  g\in\sG_h.
 \]
Now for $g\in\sG$, the action of $\sG$ on the component $\sA''$ is
standard:
 \[
 \na^{0,1}_{g(a)}=g\cdot \na^{0,1}_a\cdot g^{-1}=\na^{0,1}_a-(\na^{o,1}_a g)
 \cdot g^{-1};
 \]
and the action on the first component $\sA'$ of $\p$-operators is
 \[
 \na^{1,0}_{g(a)}=\overline{g}\cdot\na^{1,0}_a\cdot\overline{g}^{-1}=
 \na^{1,0}_a-((\na^{0,1}_a g)\cdot g^{-1})^*.
 \]
It is easy to see directly that the action just described
preserves unitary connections:
 \begin{equation}
 \sG(\sA_h)=\sA_h,
 \end{equation}
and that the identification $\sA_h=\sA''$ is equivariant with
respect to this action.

It is easy to see that $(\na^{0,1}_a)^2 = 0\in\Om^{0,2}(ad E)$.
 Thus the orbit
space
 \begin{equation}
\sA'' /\sG=\bigcup \sM_i
 \end{equation}
is the union of all components of the moduli space of
topologically trivial $I$-holomorphic bundles on $\Si_I$. (Although this
union doesn't admit any good structure, as it contains all
unstable vector bundles, however we can use the notion {\it stack}
for such situation .)  Finally, the image of $\sA_F\in \sA_h$ is
the component $M^{ss}$ of maximal dimension ($3g-3$) of s-classes
of  semi-stable vector bundles. Thus by classical technique of GIT
of Kempf--Ness type we get:

 \begin{prop} {\it (Narasimhan--Seshadri)}
 \[
 CLRep(\pi_1(\Si), \SU(2))=M^{ss}.
 \]
 \end{prop}

 \begin{prop} The form $F_{A_{CS}}$ (3.85) is a $(1,1)$-form and
the line bundle $L$ admits  unique holomorphic structure
$\Oh(\Theta)$ compatible with the Hermitian connection $A_{CS}$.
 \end{prop}

On the other hand,  complex structure $I$ on $\Si$ defines a
K\"ahler metric on $M^{ss}$ (the so-called Weil--Petersson metric)
with K\"ahler form
 \begin{equation}\om_{\mathrm{WP}}=i F_{A_{CS}}=i\cdot\om.
 \end{equation}
This metric defines the Levi-Civita connection on the complex
tangent bundle $T M^{ss}$, and hence a Hermitian connection
$A_{LC}$ on the line bundle
 \begin{equation}
\det T M^{ss}= \Oh(4 \Theta),
 \end{equation}
and a Hermitian connection $\frac{1}{4}A_{LC}$ on $L$ compatible
with the holomorphic structure on $L = \Oh(\Theta)$. Thus we have

 \begin{prop}
 \[
 \frac{1}{4} A_{LC}=A_{CS}.
 \]
 \end{prop}

Finally, considering $M^{ss}$ as a family of $\op$-operators, we
get the Quillen determinant line bundle $L$ with a Hermitian
connection $A_Q$ with curvature form
 \begin{equation}
F_{A_Q}=i\cdot\om.
 \end{equation}
Hence we can extend the equality of Proposition 4:
 \begin{prop}
 \[
 \frac{1}{4} A_{LC}=A_{CS}=A_Q.
 \]
 \end{prop}

Now we would like to describe Higgs fields in  gauge terms.

A holomorphic Higgs field $\phi$ is given by a section from
$H^0(\Si_I, ad E \otimes K)$ and corresponding quadratic
differential
\begin{equation}
\pi (E, \phi) = tr(\phi^2).
\end{equation}
Every quadratic differential $w$ is a holomorphic function on the
dual space $\Om^{0,1}(\Si_I, K^{-1})$ of Beltrami differentials.
 This function is given
by the integral
\begin{equation}
H_w = \int_\Si tr(\phi^2) \cdot \dot I.
\end{equation}
 Obviously this function from the cotangent bundle
$T^* M^{ss}$ can be extended to the full moduli space of Higgs pair
$MH_\Si$.

\begin{prop} For every two holomorphic quadratic differentials
 $w, w'$
functions $H_w$ and $H_{w'}$ are Poisson commuting with respect
to the holomorphic symplectic structure.
\end{prop}
Indeed this question can be lifted to the affine space of
connections or Cauchy- Riemann operators on $E$ (3.87). The new
affine space is
\begin{equation}
T^* \sA'' = \sA'' \times \Om^{0,1}(ad E).
\end{equation}
But every function $H_w$ as a function on this flat space
doesn't depend on $\sA''$  - variables. But $\sA''$-variables and
$\Om^{0,1}(ad E)$ - variables are conjugate with respect to the
canonical symplectic form on the flat space. The zero-level of
the moment map for the complex gauge group action (see (3.88)) is
precisely the space of holomorphic Higgs fields.  Thus $H_w$ and $H_{w'}$
Poisson commute.

\subsection{Computations of ranks}

Computations of ranks in non-abelian setup is in a sense
absolutely new procedure using a new idea -  consideration
 of full collection of moduli spaces
 \begin{equation}
 M^{ss}_{2}, M^{ss}_{3}, \dots , M^{ss}_{g}, \dots
 \end{equation}
for all genus. In the complex polarization case ideas are coming
from the Donaldson theory for 4-folds (see, for example, \cite{T2}
and \cite{T3}) where the second Chern class $c_{2}(E)$ is playing
the role of the parameter $g$. In the real polarization case
 constructions are coming from  Conformal Field Theory.
 \begin{rmk}
 It will be very useful  to relate the Donaldson theory and CFT directly
  as it was done by Nakajima for Hilbert schemes $Hilb^{g} S$ of
  an algebraic surface $S$.
 \end{rmk}

So let $\Si$ be a smooth algebraic curve of genus $g$ and
 \begin{equation}
\sS (\Si)= Sym (H_{even}(\Si,\R)) \otimes \Lambda(H_{odd}(\Si, \R))
\end{equation}
 be the graded algebra where $\deg(\si)=4-i$ if $\si \in
H_i(\Si)$. Let $M^{ss}_g$ be the moduli space of rank 2
semi-stable vector bundles on $\Si$ with trivial determinant. On
the direct product $ \Si \times M^{ss}_g$ there isn't the
universal bundle $\sE$ but there exists the bundle $\ad \sE$ with
its Chern class
 \begin{equation}
- p_1(\sE)= c_1(\sE)^2 - 4 c_2(\sE)= c_2(\ad \sE) \in H^4(\Si \times
M^{ss}_g).
 \end{equation}
This class is defined by the slant product, that is, the
homomorphism $\mu\colon H_*(\Si)\to H^*(M^{ss}_g)$ defined by
 \begin{equation}
\mu (\si)=-\frac{1}{4} p_1(\sE)/ \si.
\end{equation}
So, $ \deg(\al)=i \implies \mu(\al) \in H^i(M^{ss}_g)$.

More precisely,
 \begin{equation}
-p_1(\sE)=-2[\Si] \otimes [\Theta]+\pt \otimes x - 4
\sum_{i=1}^g(a_i \otimes \mu(a_i)+b_i \otimes \mu(b_i)),
\end{equation}
where the cohomology class of the theta diviser
 \begin{equation}
[\Theta] \in H^2(M^{ss}_g, \Z);\quad  x \in H^4(M^{ss}_g, \Z);
 \end{equation}
$[\Si]$ is the fundamental class of $C$; $[\pt]$ is the class of a
point of $\Si$ and $ a_i, b_i, i=1,\dots,g$ is a symplectic basis
of $H^1(\Si, \Z)$.

\begin{dfn}
The Newstead invariant is a linear function
 $$
N_g \colon \sS (\Si) \to \Q
 $$
sending a typical monomial $\al= \si_1 \cdot \cdot \cdot \si_r$
where $\si_i \in H_{n_i}(\Si)$ of total degree
$d(\al)=\sum_{i=1}^{r} (4-n_i)$
 to 0 if $d \neq 6g-6= dim M^{ss}_g$ and to the index of
 intersection
 \begin{equation}
 \mu(\si_1) \cdot \cdot \cdot \mu(\si_r)
 \end{equation}
 if $d=6g-6= dim M^{ss}_g$.
Such homogeneous polynomial of degree $3(g-1)$ we will call {\it
Newstead polynomial}.
\end{dfn}
 \begin{rmk}
  Of course we can apply this constructions to any
 family of vector bundles, that is, to any vector bundle $\sE$ over
 $ \Si  \times B $ where $B$ is a base of this family. Of course in this
 case we have one
 polynomial  which we will denote by $ N(\sE)$.\end{rmk}

Every Newstead polynomial  is  $Diff \Si$-invariant. From this we
can see its shape:
\begin{enumerate}
\item $$
 \quad N_g( \mu([\Si])^m \cdot \mu([\pt])^n \cdot \mu(a_{i_{\min}
 })\cdot \mu(b_{j_{\min} }) \cdot \cdot \cdot) \neq 0 \implies i_{\min}
=j_{\min}.
 $$
 Indeed, we  use the commutation relations,
 the diffeomorphism invariance of
 the Newstead polynomials and the existence of an orientation
 preserving
 diffeomorphism $\phi \colon \Si \to \Si$ satisfying
 $$
 \phi^* ( a_i)= -a_i;
 \phi^*(b_i)=-b_i; \phi^*(a_j)=a_j, \phi^*(b_j)=b_j, i\neq j
 $$
 (the
 Dehn twist) and get the statement.
\item The value
 $ N_g( \mu([\Si])^m \cdot \mu([\pt])^n \cdot \mu(a_{i_1
})\cdot \mu(b_{i_1}) \cdots )$ is independent on the choice of
handles $i_j$ because of
 the existence of an orientation preserving diffeomorphism
interchanging any pair of handles;
 \item the expression $ \ga=2
\sum_{i=1}^{g} \mu(a_i) \cdot \mu(b_i)$ is independent on the
choice of symplectic basis $\{a_i, b_i \}$.
\end{enumerate}

From this it is easy to see that if
\begin{equation}
\ga =2 \sum_{i=1}^g \mu(a_i) \cdot \mu(b_i) \in H^6(M^{ss}_g )
\end{equation}
then
 \begin{equation}
N_g \in \C[[\pt],[\Si], \ga].
\end{equation}
Thus this polynomial is nonzero only on the monomials of the type
$ [C]^m \pt^n \ga^p;\text{ for } m+2n+3p=3g-3 $.

Since the nondegenerate pairing (3.102) only involves the classes
$[\Si], [\pt], \ga $  and does not depend on the conformal
structure on $\Si$, we can reduce our  of the graded algebra (3.97)
to the "universal" algebra
 \begin{equation}
\Bbb S_3=  Sym ( \al \cdot \Z \oplus \be \cdot \Z \oplus \ga \cdot
\Z),
\end{equation}
where $\deg \al=1; \deg \be=2, \deg \ga=3 $. Then the collection
of the polynomials $N_g$ becomes the linear map
 \begin{equation}
N \colon \Bbb S_3 \to \Q
\end{equation}
sending a monomial of degree $d$ not divisible by 3 to 0 and a
monomial of degree $3n$ to $N_{n+1}$ (3.102) of this monomial. So
we have non homogeneous polynomial $N$ with homogeneous components
of degree $3g-3$.

Moreover, there is the following recurrence relation between homogeneous
components of $N$ (3.106):
 \begin{equation}
N_g (\ga \cdot \al)=g \cdot N_{g-1}(\al)
\end{equation}
(see for example \cite{Th2}, Proposition 26). This is very
important point: we get the inductive statement relating the
Newstead polynomials for different genus.
\begin{rmk}
We will see later the fusion rule of the same type in the real
polarization case.
\end{rmk}

It is convenient to use the following natural normalization
\begin{equation}
 N^0_{3n}=\frac{1}{(n+1)!}
N_{n+1}.
\end{equation}
We call such function the normalized Newstead polynomial. For
this polynomial
 \begin{equation}
 N^0\colon  \Bbb S_3 \to \Q
\end{equation}
we have the equality
 \begin{equation}
 N^0(\ga \cdot z) = N^0(z)
 \end{equation}
 for the generator $\ga$ and arbitrary $z \in \Bbb S$.
So we can exclude the generator $\ga$ and consider $N^0$ as a
linear map
 \begin{equation}
N^0 : \Bbb S_2= Sym_*(\al \cdot \Z \oplus \be \cdot \Z) \to \Q.
\end{equation}
The  so-called Newstead Conjecture ( see [Th2]) predicts that
 \begin{equation}
N^0(\al^m \cdot \be^k) \neq 0 \implies k \le \frac{1}{3}(m+2k).
\end{equation}
 So it is
convenient to let $\al\cdot\be=\om$ and such that the Newstead
polynomial is
 \begin{equation}
N^0 : Sym_*(\al \cdot \Z \oplus \om \cdot \Z) \to \Q
\end{equation}
satisfying the relation
\begin{equation}
N^0(\al^{3k} \cdot \om^{n})=N_{3(n+k)+1}(\al^{n+3k } \cdot \be^n).
 \end{equation}
From this we have the main equality:
 \begin{equation}
N^0(\al^{3k} \cdot \om^{n})=\frac{(-4)^{n+k}\cdot
(n+3k)!}{(n+k+1)!} \cdot (2k)! \cdot(4^k-2)\cdot B_{2k},
\end{equation}
where $B_i$ is the $i$th Bernoulli number (see \cite{Th2}, formula
(29)).

For any linear form $L \in
 Sym_*(\al \cdot \Z \oplus \be \cdot \Z) \oplus \ga \cdot \Z$, the
 interior product by any homogeneous element $x$ gives the differential
  operator
 \begin{equation}
 \frac{\partial}{\partial x} L(z)=(\deg z+\deg x) L(x \cdot z).
 \end{equation}
  Thus we can consider the relations of the type (3.110)
 as a differential equation. From this we  get the Newstead
  polynomials as  solutions to such equations. Moreover we
  recognize  the
 generating functions of the pairing (3.102) as the standard elementary
 functions (see for example \cite{Th1}).

\begin{rmk} \begin{enumerate} \item
 For the simplest graded algebra $ Sym_*(\Z \oplus \cdots
\oplus \Z)$ the product $ z_1 \cdot z_2 \in  Sym_{d_1+d_2}$of
 two homogeneous elements $ z_i \in Sym_{d_i}$ is defined as
 $$
z_1 \cdot z_2(\si_1, \dots, \si_{d_1+d_2})=
 $$
 \begin{equation}
= \frac{1}{(d_1+d_2)!}\sum_{g \in \Sigma_{d_1+d_2}}
z_1(\si_{g(1)}, \dots, \si_{g(d_1)}) \cdot z_2(\si_{g(d_1+1)},
\dots, \si_{g(d_1+d_2)}),
\end{equation}
where $\Sigma_n$ is the symmetric group on $n$ letters
\item Witten's approach deals with the volumes of moduli spaces
\begin{equation}
vol (M_{g}^{ss}) = N^{0}(\al^{3g-3})
\end{equation}
and the other values of the Newstead polynomials can be read off
from these formulae.
\item Historically the
 procedure
 was the following: Thaddeus reconstructed values of the Newstead
 polynomial (3.115)
 from the beautiful form of the generating functions for
 ranks
 of   conformal block spaces given by the Verlinde formula
 (see below).
\item To use these computations Thaddeus proposed the following
conjecture: {\it The quantization of $M^{ss}_{g}$ is numerically
perfect} (see subsection 1.5).
\item To avoid Thaddeus's conjecture we  use the
Donaldson's construction from \cite{Do}.
\item Other regular way to get the pairing formuli and the computation of
the collection of Newstead polynomials is the Geometric
Approximation Procedure: there is a chain of birational
transformations or flips (see \cite{Re})
\begin{equation}
 \text{MP}_{\max}\leftrightarrow\dots\leftrightarrow\text{MP}_1
\leftrightarrow\text{MP}_0,
 \end{equation}
where every $\text{MP}_i$ is the moduli space of $t$-stable pairs
(see \cite{Th2}, \cite{BD} or \cite{Ber}), and can be realised as
a family of vector bundles over $\Si$. Thus using the same
construction  we can define  polynomials $NP_{g,i}$ analogeous of
the Newstead polynomials, and the linear forms $NP_i$,  the
analogs of (3.114) and (3.115). Recall that the moduli space on
the extreme left admits the $\PP^{1}$-bundle structure that is
regular map
 \begin{equation}
\text{MP}_{\max}\to M_g^{ss}
\end{equation}
to the moduli space of stable bundle with trivial determinant.
From this it is easy to see that the form
 \begin{equation}
NP^{0}_{\max}\cdot H=N^0
 \end{equation}
where $H$ is the Grothendieck class of the $\PP^{1}$-bundle,
$N^{0}$ is the Newstead form.
\item
Now the moduli space $\text{MP}_0$ on the extreme right-hand end
of (3.119) is the projective space of dimension $3g-2$.
 \begin{equation}
\text{MP}_0=\PP^{3g-2}.
\end{equation}
\item Using the description of every flip in the chain (3.119) we can
reconstruct all polynomials step by step begining with
$\PP^{3g-2}$.
\end{enumerate}
\end{rmk}

Using the Hirzebruch-Riemann-Roch theorem we have
\begin{equation}
\sum_{i} rk H^{i}(M_{g}^{ss}, \Oh(k \Theta)) = (ch(\Oh(k\Theta))
\cdot td (M_{g}^{ss}))_{3g-3}.
\end{equation}
Since the divisor $\Theta$ is positive the LHS of previous
equality is equal to
$$
rk \sH^{r}_{I} = rk H^{0}(M_{g}^{ss}, \Oh(k \Theta))
$$
by Kodaira vanishing theorem. On the other hand, Newstead showed
(see \cite{Ne}) that
\begin{equation}
c_{1}(M^{ss}) = 4 \cdot [\Theta] \quad \text{ and } p(M_{g}^{ss})
= (1 + \beta)^{2g-2}.
\end{equation}

Hence  the Todd class
\begin{equation}
td(M^{ss}) = exp(2 [\Theta]) (\frac{1/2 \sqrt{\beta}}{sinh(1/2
\sqrt{\beta})})^{2g-2}
\end{equation}
Using values of the Newstead polynomials we get the RHS of (3.123)
:
\begin{equation}
rk H^{0}(M_{g}^{ss}, \Oh(k \Theta)) = (\frac{k+2}{2})^{g-1}
\sum_{j}\frac{1}{sin (\frac{j \pi}{k+2})^{2g-2}}.
\end{equation}

\begin{rmk}
All of the computations are due to  Zagier and Macdonald (see
\cite{Th1} for the corresponding references).
\end{rmk}

\subsection{Hitchin's  connections}

Here we will see that
\begin{enumerate}
\item The sending $ \Si_I$ to $ M^{ss}(\Si_I)$ (3.1) is a
 faithful functor and
\item the sending $\Si_I \to H^0(M^{ss}(\Si_I), \Oh(k \Theta))$ is a
successful quantization.
\end{enumerate}
That is, for $g > 2$ we can reconstruct $\Si_I$ from $M^{ss}$ and
every vector bundle from the collection of vector bundles on the
moduli space $\sM_{g}$ of curves of genus $g$
\begin{equation}
 p_k \colon V_k \to \sM_g ; \quad p_k^{-1}(\Si_I) =  H^0(M^{ss}(\Si_I),
 \Oh(k \Theta))
  \end{equation}
admits a projective flat connection. So, locally the
projective space $\PP H^0(M^{ss}(\Si_I), \Oh(k \Theta))$ doesn't
depend on $I$.

The first statement is almost obvious since
\begin{equation}
Sing M^{ss} = K(\Si_I)
 \end{equation}
is the Kummer variety of our curve $\Si_I$. Blowing up the
 finite set $Sing K(\Si_I)$, considering the double cover
 with  ramification in the exceptional divisor and then blowing
  down exceptional
  divisor we get the Jacobian $J(\Si)$. It is easy to see that
   the restriction of non-abelian theta divisor $\Theta$ to
   the Kummer variety gives the divisor $2 \Theta$ on the
   Jacobian.
Now the classical Torelli theorem gives us our curve $\Si_I$.

Recall that in the abelian case  under deformations of
 a complex structure on $\Si$ we get vector bundles

$\pi \colon V_k \to \sM_g, \quad \pi^{-1}(\Si) = H^0(J(\Si), \Oh
(k \Theta))) $.

It is the classical result  ( we saw this already in
terms of the extended Kodaira-Spencer theory  \cite{W})  that $V_k$
(1.29) admits a projective flat connection.

In parallel to this fact  Hitchin proved that even for
 non-abelian
theta function case $V_k$ admits
 a projective flat connection.

To construct the holomorphic projective connection we can use
 the extended Kodaira-
Spencer theory again. First of all our goodness conditions (1.48)
 holds:
\begin{equation}
H^0(T M^{ss}) = H^1(M^{ss}, \Oh) = 0
 \end{equation}
by Narasimhan-Ramanan result \cite{NR}. Now we have to construct
a holomorphic section $W_0 \in H^0(S^2 (T M^{ss}))$. To do this,
consider a  semi-stable vector bundle $E \in M^{ss}$. Then the fiber
\begin{equation}
T M^{ss}_E = H^1(\Si, \ad E)
 \end{equation}
where $ad E \subset End E$ is the traceless part of the sheaf of
 endomorphisms. The
dual space (using the Serre duality)
\begin{equation}
T^* M^{ss}_E = H^0(\Si, \ad E \otimes K_\Si)
 \end{equation}
where $K_\Si$ is the canonical class of $\Si$.

Now a tangent vector to the moduli space $\sM_g$ of
 Riemann surfaces at a point
$\Si_I$ is given by a Beltrami differential
\begin{equation}
\dot I \in H^1(\Si_I, K_\Si^*).
 \end{equation}
We have the cup-product map
\begin{equation}
H^1(\Si_I, K_\Si^*) \otimes H^0(\Si, \ad E \otimes K_\Si) \to
H^1(\Si, \ad E)
 \end{equation}
which is nothing else but a quadratic form on $T^* M^{ss}_E$ and
being extended on full moduli space $M^{ss}$ it gives the required
section of the symmetric power of the tangent bundle. Using the
standard machinery of the end of section 5 we get
 a holomorphic projective connection.

The next question is about the curvature. Here it is quite
reasonable to return to the gauge theory terminology.

Recall that the affine space of all unitary connections $\sA_h$
is canonically a symplectic manifold. Since the tangent space to
$\sA$ at a connection $a \in \sA$ is given by the space of forms
$\Om^1(\fsu(2))$, the symplectic form is given by the formula
(3.74). This form doesn't depend on $a$, so it is closed. The
gauge group $\sG$ acts on the space $\sA_h$ in the usual way  and
preserves the form $\om$.

 For the moment map (3.70) the symplectic quotient is
\begin{equation}
m^{-1} (0) / \sG = CLRep(\pi_1(\Si), \SU(2)).
 \end{equation}
(see (3.70)). The symplectic reduction of the symplectic form
defines a symplectic structure on this space. Thus
\begin{equation}
(CLRep(\pi_1(\Si), \SU(2)), \om)
 \end{equation}
is a phase space of the classical mechanical system. Moreover we
have the prequantization line bundle (3.83) with the connection
(3.85).

For every connection $a$ its covariant derivative $\nabla_a$
defines the complex
\begin{equation}
\Om^0(\fsu(2)) \to \Om^1(\fsu(2)) \to \Om^2(\fsu(2))
 \end{equation}
and the tangent space to $CLRep(\pi_1(\Si),\SU(2))$ at a point $a$
is given as the first cohomology group of this complex
\begin{equation}
T CLRep(\pi_1(\Si), \SU(2))_a = H^1_a.
 \end{equation}

 Kahler polarization is just the choice of a conformal structure
 which can be done using
the Hodge star operator $*$ acting on the space of forms
\begin{equation}
* \colon \Om^1 \to \Om^1.
 \end{equation}
Our complex structure $I$ acts as a endomorphism
\begin{equation}
I = - *
 \end{equation}
and
\begin{equation}
\om( \ga, I \ga) = - \om( I \ga,  \ga) = - \int_\Si tr(\ga, *\ga).
 \end{equation}

Such complex structure defines the decomposition
\begin{equation}
\nabla_a = \nabla^{1,0} + \nabla^{0,1}
 \end{equation}
for every $a \in \sA$ and identifies $\sA$ with  the complex
affine space
 of Cauchy-Riemann operators on $E$ taking (0,1)-part of $\nabla_a$. Thus the
 the tangent space
 \begin{equation}
T^{1,0}_a CLRep(\pi_1(\Si), \SU(2))_I = \Om^{0,1} (\fsu(2)).
 \end{equation}
This complex structure is invariant with respect to the gauge
group $\sG$ action thus we get an integrable complex structure as
the result of the symplectic reduction (3.89).

This complex structure $*$ on $\Si$ induces the complex structure
$I$ on the space of classes of representations $CLRep(\pi_1(\Si),
\SU(2))$ and we have the main identification
\begin{equation}
CLRep(\pi_1(\Si))_I = M^{ss} (\Si_{-*}).
 \end{equation}

The last remark which we make before the curvature computations is
following:
\begin{prop} The map
\begin{equation}
f \colon H^1(\Si_I, T^{1,0}) \to H^0(M^{ss}, S^2(T^{1,0} M^{ss}))
 \end{equation}
given by the formula
\begin{equation}
f(\dot I)(\phi) = \int_\Si tr(\phi^2) \wedge \dot I
 \end{equation}
is an isomorphism.
\end{prop}

Indeed by the Hartog's principle we can extend every holomorphic
quadratic form from $T^* M^{ss}$ to $MH_\Si$. But every
holomorphic form is  constant on the fibers of the map $\pi$
(3.54) (which are compact abelian varieties). Hence the space of
holomorphic quadratic forms is the dual space to the space of
holomorphic quadratic differentials on $\Si_I$ and we are done.

Now let $0 \in B \subset \C^m$ be an open set and
\begin{equation}
p \colon \tilde{M^{ss}} \to B , \quad p^{-1} (0) = M^{ss}_I
\end{equation}
be a body of holomorphic deformations of complex structure.
Consider a section
\begin{equation}
t \in H^0(B, TB)
 \end{equation}
that is a holomorphic vector field on $B$. By the classical
 Kodaira-Spencer
theory it define a section
\begin{equation}
\xi \in H^0(B, R^1 T\tilde{M^{ss}}_{/p})
 \end{equation}
that is the family of Kodaira-Spencer classes for deformations
along $t$. Moreover considering our vector field $t$ as a
differential operator on elements of a small open cover of fibers
we get the section
\begin{equation}
D_t \in H^0(B, R^1(D^1(L))_{}/p)
 \end{equation}
that is the family of cohomology classes $H^1(M^{ss}_b, D^1(L))$
over points of $B$. This section describes deformations of a
polarization line bundles over $t$.

Having the holomorphic quadratic form (3.133)
\begin{equation}
W \in H^0(B, R^0(S^2(T_{/p})
 \end{equation}
we can lift it to the family of Laplace-Beltrami operators
$\Delta$ (using the Levi-Civita
 connection for the metric (3.140)) over elements of open covers
 of fibers of $p$,
and their differences  define
\begin{equation}
D \in H^0(B, R^1(D^2(L)_{/p}))
 \end{equation}
because the operators of second order have the same principal
symbols. Described in (1.86) general construction gives us the
projective connection
 on the vector bundles $V_k \vert_B$ (3.146) with the principal symbol
 \begin{equation}
- \frac{1}{2k - 4} \si \cdot \de (W) \in H^0(B, R^1
T\tilde{M^{ss}}_{/p}).
 \end{equation}
Moreover, the combination of these operators
\begin{equation}
\p_t + 1/(2k -4) i \Delta +i D
 \end{equation}
gives  globally defined holomorphic {\it heat operator} on $L_I$
\begin{equation}
\p_t + P_t
 \end{equation}
and the principal symbol of $P_t$ is a quadratic form lifted from
the target of the moment map $\pi$ (3.93). Thus by Proposition 14
these symbols  commute as  functions on the cotangent
bundle.

What we have to do now is just to compute the commutator of these
operators for two holomorphic vector fields $t$ and $t'$ from
(3.147):
\begin{equation}
\p_t P_{t'} - \p_{t'} P_t - [P_t, P_{t'}].
 \end{equation}
The last term  commutator is of order 3  a priori. But its
principal symbol  is the Poisson brackets of symbols of operators
(3.154)
 (see \cite{GS1}
). They are Poisson  commuting. Thus this commutator
\begin{equation}
\p_t P_{t'} - \p_{t'} P_t - [P_t, P_{t'}] \in H^0(M^{ss}, D^2(L))
 \end{equation}
for any complex structure from $B$.

\begin{prop} This operator has to be  constant.
\end{prop}

Indeed, in the cohomology sequence of the triple (1.81) there is
the part
\begin{equation}
\to H^0(M^{ss}, S^2 T) \to H^1(M^{ss}, D^1(L))
 \end{equation}
which is a monomorphism because the combination of this
homomorphism with the symbol projection $H^1(M^{ss}, D^1(L)) \to
H^1(T)$ is an isomorphism by Proposition 15 (3.144). Thus the
beginning part of this exact sequence
\begin{equation}
H^0(M^{ss}, D^1(L)) \to H^0(M^{ss}, D^2(L)) \to 0.
 \end{equation}
In particular,
\begin{equation}
H^0(M^{ss}, D^1(L)) = H^1(M^{ss}, D^2(L))
\end{equation}
and in the same vein
\begin{equation}
H^0(M^{ss}, D^1(L)) = \C
\end{equation}
Thus the commutator in (1.154) is a constant and our projective
connection is projectively flat.

Problem: is this connection unitary? Or is there another unitary
projective flat
 connection?

Summarizing that we have

\begin{thm} The Kahler quantization of moduli spaces of topological trivial
 rk2 vector bundles is succesful.
\end{thm}

\section{Symplectic geometry of moduli spaces of vector bundles
 }

\subsection{Goldman $U(1)$-action}

 The symplectic-geometrical part of the
non-abelian theory of theta functions is well developed by W.
Goldman, J. Weitsman,
 L. Jeffrey and  J. Hurtubise.

 The starting point is the Narasimhan-Seshadri identification
 $ M^{ss}(\Si) = CLRep(\pi_1, \SU(2))$ as a result of the
  Symplectic Reduction described in the previous subsection:
 $\sA_h $ is the space of $\SU(2)$-connections, $\sG_h$ is
 the $\SU(2)$-gauge group, and
 $ \om (\ga_1, \ga_2) = \int_\Si tr (\ga_1 \wedge \ga_2)$ is
 the $\sG_h$-invariant symplectic form on $\sA_h$ (see (3.74)).

The map $F \colon \sA \to \Om^2(ad)$ (3.66) sending
  connection to its curvature form is
 the moment map for $\sG_h$ - action. Hence
 $F^{-1}(0) = \sA_F$ is the space of flat connections and
 $( CLRep(\pi_1, \SU(2)), \om )$ is the result of
 the symplectic reduction with respect to this action.

 Any simple closed loop $C$ on $\Si$ defines the Goldman
  function on the space of classes of representations
  $CLRep(\pi_1, \SU(2))$:
\begin{equation}
c_C(\rho) = 1/\pi \cdot cos^{-1}(1/2 tr(\rho([C]))) \in [0, 1].
\end{equation}
In \cite{G2} Goldman proved that $c_C$ is a Hamiltonian of
$U(1)$-action on $CLRep(\pi_1, \SU(2))$. This action can be
described as follows:
\begin{enumerate}
\item for every $g \in \SU(2) , \quad g \neq \pm id$  there are exist
 unique element
\begin{equation}
log\quad g \in \fsu(2), \quad |log\quad g|^2 = 2,
\end{equation}
and a real number $t \in \R$  such that
\begin{equation}
 g = e^{t \cdot log\quad g};
\end{equation}
\item obviously  $[g, e^{t \cdot log\quad g}] = 0$ for every
$t \in \R$.
\item Suppose that the class $[C] = a_1 $ in a standard presentation of
the fundamental group of $\Si$:
\begin{equation}
 \pi_1 = <a_1, \dots ,a_g, b_1, \dots , b_g \vert
  \prod_{i=1}^g [a_i, b_i] = id>,
\end{equation}
\item for any $\rho \in CLRep(\pi_1, \SU(2))$
such that $\rho (a_1) \neq \pm id$ the formula
\begin{equation}
\rho (b_1) \to \rho (b_1) \cdot e^{t \cdot log \rho(a_1)}
\end{equation}
defines Hamiltonian $U(1)$-action on $CLRep(\pi_1, \SU(2))$.
\end{enumerate}
It is easy to see that
\begin{equation}
C_1 \cap C_2 = \emptyset \implies \{c_{C_1}, c_{C_2} \}_\om = 0
\end{equation}
where the Poisson brackets are taken with respect to the canonical
symplectic structure $\om$ (3.74). Thus we have to require a
maximal set of disjoint homotopy inequivalent loops which we call
circles . It easy to see (you will see this in a minute) that the
maximal number of circles in such set is equal to 3g-3. Moreover,
let a maximal set of disjoint inequivalent circles be
\begin{equation}
\{ C_1, ... , C_{3g-3} \}.
\end{equation}

\begin{figure}[tbn]
\centerline{\epsfxsize=3in\epsfbox{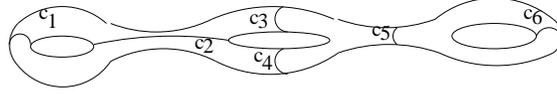}} \caption{\sl Maximal
collection of simple loops} \label{Fig 2}
\end{figure}

Then removing these circles we get
\begin{equation}
\Si - \{ C_1, ... , C_{3g-3} \} = \cup_{i=1}^{2g-2} \tilde{v_i}
\end{equation}
 is a finite set of trinions (or "pairs of pants").

\begin{figure}[tbn]
\centerline{\epsfxsize=3in\epsfbox{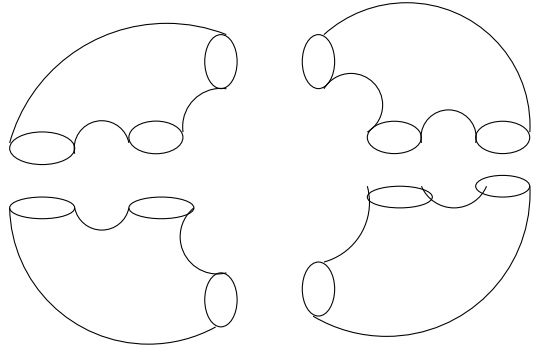}} \caption{\sl Pairs
of pants} \label{Fig 2}
\end{figure}

 A trinion (or "pair of pants")  plays the role of points
in algebraic geometry under investigations of geometry of
connections: just like we restrict functions to a point we can
restrict connections to a trinion. Roughly speaking, a trinion is
the minimal geometrical object admitting non trivial non-abelian
gauge theory.

 Riemann surface with a trinion decomposition (4.8) is called
{\it marked} Riemann surface.
 Any marked surface defines  the {\it dual
 trivalent graph} $\Ga$
\begin{enumerate}
\item  the set of vertices of which
\begin{equation}
V(\Ga) = \{ v_i \} = \{ \tilde{v_i} \}
\end{equation}
is the set of trinions;
\item  and with the set of edges
\begin{equation}
E(\Ga) = \{ e_i \} = \{ C_i \};
\end{equation}
\item and two vertices $v_i$ and $v_j$ are joined by the
 edge $e_l$ iff the corresponding circle
\begin{equation}
C_l = \p \tilde{v_i} \cap \p \tilde{v_j}.
\end{equation}
\end{enumerate}

\subsection{Real polarization}

Every marked Riemanian surface $\Si_\Ga$ defines a real
polarization of the prequantized phase space
\begin{equation}
 (CLRep(\pi_1, \SU(2)),\om, L, A_{CS}).
 \end{equation}

It is given in a very geometric way in the set-up of perturbation
theory of 3-dimensional  Chern--Simons theory. As we saw, a mark
of  Riemann surfaces is given by the choice of  maximal
collection of disjoint, noncontractable, pairwise nonisotopic
smooth circles on $\Si$. Any such system contains $3g-3$ simple
closed circles (4.7).

The isotopy class of a trinion decomposition $\{C_i\}$ gives the
map
 \begin{equation}
\pi_\Ga \colon CLRep(\pi_1, \SU(2)) \to \R^{3g-3}
 \end{equation}
with fixed coordinates $(c_1, \dots, c_{3g-3})$ such that
 \begin{equation}
 c_i (\pi_{\{C_i\}} (\rho))=\frac{1}{\pi} cos^{-1}(\frac{1}{2} tr
\rho([C_i]))\in [0, 1].
 \end{equation}
which are just Goldman functions (4.1).

Thus
\begin{enumerate}
 \item
this map is the moment map for the Hamiltonian
$U(1)^{3g-3}$-action on $CLRep(\pi_1, \SU(2))$, thus gives a real
polarization of the system (4.12).
 \item The coordinates $c_i$ are {\it action coordinates} for this
Hamiltonian system ( \cite{GS1}).
 \end{enumerate}

These functions $c_i$ are continuous on all $CLRep(\pi_1,
\SU(2))$ and smooth inside of the unit cube  $[0,1]^{3g-3}$.

Thus we have
 \begin{enumerate}
 \item The image of $CLRep(\pi_1, \SU(2))$ under $\pi_\Ga$ is
 a convex polyhedron
 \begin{equation}
 \Delta_\Ga\subset [0, 1]^{3g-3}.
 \end{equation}
\item The symplectic volume of $CLRep(\pi_1, \SU(2))$ is equal to
the Euclidean volume of $\Delta_\Ga $:
 \begin{equation}\int_{CLRep(\pi_1, \SU(2))}\om^{3g-3}=Vol
 \Delta_\Ga.
 \end{equation}
\item The  set of Bohr--Sommerfeld orbits of the real
polarization $\pi_\Ga$ of level $k$
 \begin{equation}
 BS_{k}(\pi_{\Ga},CLRep(\pi_1, \SU(2)),\om,L,A_{CS})
 \end{equation}
 is   the set of $1/2k$- integer points in
  the polyhedron $\Delta_\Ga$,
and
 \begin{equation}
 lim_{k\to\infty} k^{3-3g}\cdot N_{k-BS}=
 \int_{CLRep(\pi_1, \SU(2))}\om^{3g-3}
 =Vol \Delta_\Ga
 \end{equation}
 just as predicted by general quantization rules.
 \end{enumerate}

\begin{figure}[tbn]
\centerline{\epsfxsize=3in\epsfbox{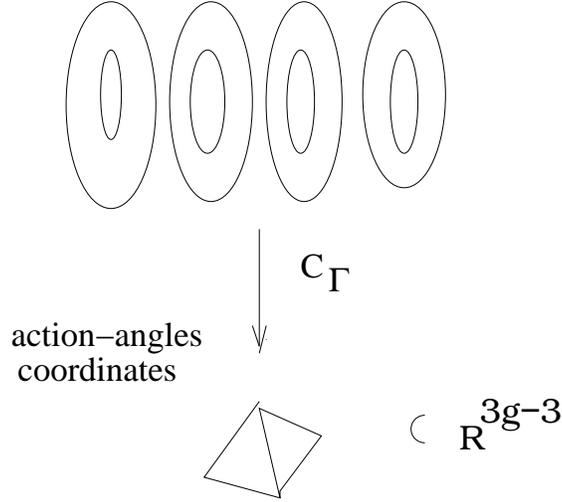}} \caption{\sl Moment
map} \label{Fig 2}
\end{figure}

To see the statements about Bohr-Sommerfeld fibers remark that
 every vector $w\in \R^{3g-3}$ can be interpreted as a {\it
differential $1$-form} on $\R^{3g-3}$, and by the usual
construction using the symplectic form $\Om$, this defines a
vector field $\xi_w$ tangent to the fibers of $\pi$. Integrating
such vector fields, one defines the collection of transformations
 \begin{equation}
 \{t_w\}= \{e^{\xi_w} \} \subset Diff^+ (CLRep(\pi_1, \SU(2)).
 \end{equation}

These transformations preserve the curvature form of the
connection $A_{CS}$. Thus (because $CLRep(\pi_1, \SU(2))$ is
simple connected), there exists a collection of gauge
transformations $\al_w\in\sG_L$ of $L$ such that
 \begin{equation}
 (t_w)^*(A_{CS})=A_{CS}^{\al_w}
 \end{equation}
 (see \cite{AB} for the description of the gauge-orbit with fixed
 curvature form)

We can view such gauge transformations as $U(1)$-{\it torsors},
just as in the description of the formulas for classical theta functions
for Abelian varieties (2.49). From this
  the subset $BS_{k} \subset \Delta_{\Ga}$
of Bohr-Sommerfeld fibers of level $k$ is the subset of
$1/2k$-integer points of $\Delta$.

To describe the polyhedron $\Delta_{\Ga}$ we have to describe the
space of flat $\SU(2)$-connections on a trinion $pp$ (pair of
pants): let
\begin{equation}
  \p pp = C_1 \cup C_2 \cup C_3.
\end{equation}
Then we have the map
\begin{equation}
  \pi \colon CLRep(\pi_1(pp), \SU(2)) \to \R^3
\end{equation}
with fixed coordinates $(c_1, c_2, c_3)$ just like in (3.139). It
is easy to check (see \cite{JW1}) the following fact

\begin{figure}[tbn]
\centerline{\epsfxsize=3in\epsfbox{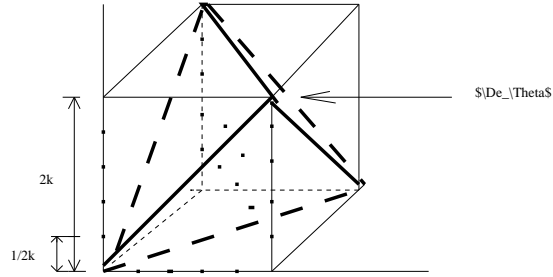}} \caption{\sl Moment
tetrahedron} \label{Fig 2}
\end{figure}

\begin{prop} \begin{enumerate}
\item the image
\begin{equation}
 \pi (CLRep(\pi_1(pp), \SU(2)) = \De_\Theta
\end{equation}
is the tetrahedron given by the following inequalities
\begin{equation}
0 \leq c_i \leq 1, \quad |c_1-c_2| \leq c_3 \leq min(c_1 + c_2, 2
- c_1 -c_2);
\end{equation}
\item  $\pi$ is 1-1 map, smooth inside of the tetrahedron and continious
everywhere.
\end{enumerate}
\end{prop}

Using the dual 3-valent graph $\Ga$ of the trinion decomposition
(4.8) we can describe the polyhedron $\De_\Ga$ (4.15):

\begin{prop} The polyhedron $\De_\Ga$ is given by inequalities (4.24) for every
 triple $(c_i, c_j, c_k)$ of the coordinates in $\R^{3g-3}$ corresponding to
 a boundary of any trinion in the decomposition (4.8).
 \end{prop}

\subsection{ Bohr-Sommerfeld fibers}

In these terms we can describe Bohr-Sommerfeld fibers  (using
action coordinates $c_i$ ) as  functions
 \begin{equation}
 w\colon E(\Ga) \to\frac{1}{2k}\{0,1,2,\dots,k\}
 \end{equation}
on the collection of edges of the 3-valent graph $\Ga $ to the
collection of $\frac{1}{2k}$ - integers, such that, for any three
edges $e_{1}, e_{2}, e_{3}$ meeting at a vertex $v$, the following
3 conditions hold:
 \begin{enumerate}
 \item $w(e_{1}) + w(e_{2}) + w(e_{3})\in \frac{1}{k}\cdot \Z$;

 \item $w(e_{1}) + w(e_{2}) + w(e_{3}) \leq 1$;

 \item and
 \begin{equation}
 |w(e_{1})-w(e_{2})|\leq w (e_{3}) \leq min( w(e_{1}) + w(e_{2}),
  2 - w(e_{1} - e_{2})).
 \end{equation}
 \end{enumerate}
Such a function $w$ is called an {\it admissible integer weight of
level} $k$ on the graph $\Ga$.

\begin{figure}[tbn]
\centerline{\epsfxsize=3in\epsfbox{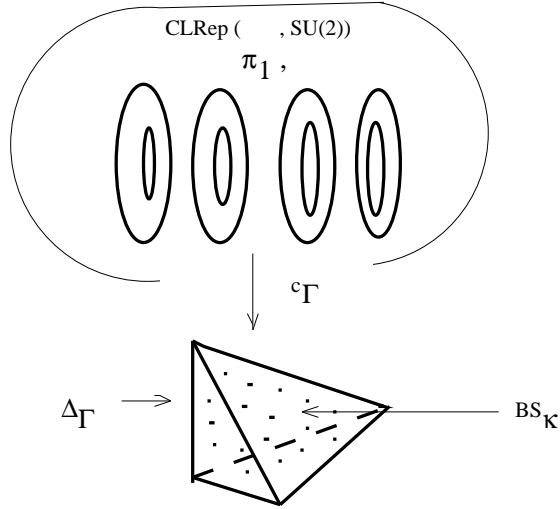}} \caption{\sl
Bohr-Sommerfeld fibers} \label{Fig 2}
\end{figure}

Let $W^{k}(\Ga)$ be the set of such weights.

\begin{rmk}
The geometrical meaning of the last inequality (4.26) is the
following: it is the triangle inequality for triangles on
2-sphere. Thus if we do our graph flat that is fix a cyclic
orientation around every vertex  and write down around every
vertex a triangle with sides orthogonal to edges and of lengths
$w(e_{1}), w(e_{2}), w(e_{3}))$, identifying equal sides we get
triangulated Riemann surface
 $\Si_{w}$  which is
a finite combinatorial cover of $S^{2}$
\begin{equation}
\phi_{w} \colon \Si_{w} \to S^{2}.
\end{equation}
For small genus, for example for the graph $\Theta$, this cover
has degree one and gives just a triangulation of the basic sphere
 $S^{2}$. But from metric point of view this is a $1/2k$-integer
  triangulation.
\end{rmk}
Here is the list of simple but very important properties of
weights:
\begin{prop}
\begin{enumerate}
 \item the number $|W^k(\Ga)|$ of admissible weights of level $k$
is independent of the graph $\Ga$ but depends on genus $g$ only;
 \item
 \begin{equation}
 |W^{k}(\Ga)| =|W_g^k| =\frac{(k+2)^{g-1}}{2^{g-1}}
\sum_{n=1}^{k+1}
 \frac{1}{(\sin(\frac{n\pi}{k+2}))^{2g-2}}
 \end{equation}
 is the Verlinde number (3.126).
 \end{enumerate}
\end{prop}

To describe the geometry of Bohr-Sommerfeld fibers given by an
 admissible integer weight $w$ (4.25) we have to consider the
 levels of this function:
for $\al\in \{0, \frac{1}{2k},\dots,1\}$ let
 \begin{equation}
w^{-1}(\al)=\Ga_1(\al)\cup\cdots\cup\Ga_n(\al)\subset\Ga
\label{eq5.1}
 \end{equation}
be the decomposition into connected components. Then every
component $\Ga_i(\al)$ is a 3-valent graph with $n_i$ univalent
vertices $a_1,\dots, a_{n_i}$ (see section 6).

Every gauge class of connections contains a connection $a_0$
adapted to a trinion decomposition (see \cite{JW1},
Definition~2.2). Fix the filtration
 \begin{equation}
Z(\SU(2))=\Z_2\subset U(1)\subset \SU(2),
 \end{equation}
and view it as the triple $ \{\Z_2, U(1), \SU(2)\}$

For a connection $[a]\in \pi^{-1}(w)$, we have the function
 \begin{equation}
e_w\colon E(\Ga)\to \{\Z_2, U(1), \SU(2)\}
 \end{equation}
 sending every loop $C_j$ to the element of $G$ conjugate to the
stabilizer of the monodromy of $[a]$ around this loop, and the
function
 \begin{equation}
v_w\colon V(\Ga)\to \{\Z_2, U(1), \SU(2)\}
 \end{equation}
sending  trinion $P_i$ to the stabilizer of the flat connection
$a\rest{P_n}$. Of course,
 \begin{align}
 & C_j\subset \partial P_n \implies v_w(P_n)\subset e_w(C_j); \notag \\
 & C_1 \cup C_2 \cup C_3=\partial P_n \quad \text{and} \quad
e_w(C_1)=e_w(C_2)=\SU(2) \implies
 \\
 & \qquad\qquad e_w(C_3)=\SU(2) \implies v_w(P_n)=\SU(2), \notag
 \end{align}
and so on.

Obviously
 \begin{equation}
e_w(C_j)= \text{$U(1)$ or $\SU(2)$.}
 \end{equation}

 \begin{prop} The functions $e_w$ and $v_w$ depend on  $w$ and not on
the choice of $[a]\in \pi^{-1}(w)$.
 \end{prop}

More precisely, they depend on the combinatorics of the
decomposition (4.29).

Thus $w$ defines the direct products
 \[
 \prod_{C\in E(\Ga)} e_w(C) \quad\text{and}\quad
 \prod_{P\in V(\Ga)} v_w(P),
 \]
and $\prod_{P\in V(\Ga)} v_w(P)$ acts on $\prod_{C\in E(\Ga)}
e_w(C)$ as follows: for
 \begin{gather*}
 g=(g_1,\dots, g_{2g-2})\in \prod_{P\in V(\Ga)} v_w(P)
 \quad\text{with $g_i\in v_w(P_i)$, and}\\
(t_1,\dots, t_{3g-3})\in \prod_{C\in E(\Ga)} e_w(C)
 \quad\text{with $t_n\in e_w(C_n)$,}
 \end{gather*}
 if $C_n\subset \partial P_i\cap\partial
P_j$ then
 \begin{equation}
 g(t_n)=g_i \circ t_n \circ g_j^{-1}.
 \end{equation}

\begin{rmk}
This  action described in detail in subsection 6.3 under the
consideration of the gauge theory on graphs.
\end{rmk}

 \begin{prop}[\cite{JW1}, Theorem~2.5] The fibre $\pi^{-1}(w)$ is given by
 \begin{equation}
 \pi^{-1}(w)=\prod_{C\in E(\Ga)} e_w(C) \Bigm/ \prod_{P\in V(\Ga)} v_w(P).
 \end{equation}
 \end{prop}

 \begin{cor}
 The fibre $\pi^{-1}(w)$ is isomorphic to
 \begin{equation}
 \pi^{-1}(w) = T^t \times [(S^3)^p \times (S^2)^s] / G_w,
 \end{equation}
where $t, p$ and $s$ are nonnegative integers and $G_w$ is the
finite Abelian group defined by $w$, or more precisely by the
combinatoric data (4.1); moreover,
 \begin{equation}
 H_1(\pi^{-1}(w))=\Z^t \oplus \Z_2^p.
 \end{equation}
 \end{cor}

Translations along the torus $T^t$ in (4.37) are induced by
Hamiltonians lifted from the target space $\R^{3g-3}$ of $\pi$.

Jeffrey and Weitsman \cite{JW1} use the normalization of the
action coordinates via branched covers to construct a covariant
constant section $s_w$ of the restrictions of $(L^k,A_{CS})$ to
$\pi^{-1}(w)$.

Following the real polarization quantization rules from section 1
(formula (1.21)) the wave function space is given by the direct
sum
\begin{equation}
\sH^{k}_{\Ga} = \oplus_{b \in BS_{k}} \C \cdot s_{b}
\end{equation}
where $s_{b}$ is a covariant constant section of the restriction
of our line bundle and connection to $\pi^{-1}(b)$.

 Geometrically  this combinatorial description of Bohr-Sommerfeld
 fibers can be done
as follows: consider the space $\R^{3g-3}$ with action coordinates
$c_i$ (4.14). This space contains the integer sublattice
$\Z^{3g-3}\subset \R^{3g-3}$, and we can consider the {\it action
torus}:
 \begin{equation}
 T^A=\R^{3g-3} / \Z^{3g-3}.
 \end{equation}

In particular, we get a map
 \begin{equation}
 \pi_A \colon CLRep(\pi_1, \SU(2)) \to T^A
 \end{equation}
which glues at most points of the boundary of $\Delta_{\Ga}$.

Now every integer weight $w$ (4.25) satisfying (1) and (2), but
{\it a priori} without the Clebsch--Gordan conditions  (3) (4.26),
defines a point of order $2k$ on the action torus
 \[
 w\in T^A_{2k}.
 \]
In particular, the collection $W^k(\Ga)$  of admissible integer
weights (subject to (3) (4.26)) can be considered as a subset of
points of order $2k$ on the action torus:
 \begin{equation}
 W^k(\Ga)\subset T^A_{2k}.
 \end{equation}
Thus for different $\Ga$ of genus $g$ we have {\it a priori}
different subsets $\{W^{k}(\Ga)\}$ of the set of points of order
$2k$ on the action torus $T^{A}$ (4.40).

The linear decomposition of the wave function space (4.39)
defines the space uniquely up to the coordinatewise
$|W^k(\Ga)|$-torus action because every covariant constant
section is defined up to $U(1)$-action. But using admissible
integer weights we can define this space as a space
 with a fixed basis
 \begin{equation}
\sH^{k}_{\Ga} = \oplus_{w \in W^{k}_\Ga} \C \cdot w
\end{equation}
and the integer sublattice
\begin{equation}
\Z^{|W^k_\Ga|} \subset    \sH^{k}_{\Ga}.
\end{equation}
All of these spaces are embedded in one large space
\begin{equation}
\sH^{k}_{\Ga} \subset \sH^k =  \oplus_{w \in T^A_{2k}} \C \cdot w
\end{equation}
(see (4.42)) with the integer lattice
\begin{equation}
\Z^{|T^A_{2k}|} \subset    \sH^{k}.
\end{equation}

\subsection{Delzant model}

Returning to the moment map (3.148) we come into the set-up of the
toric symplectic  geometry: let $(M, \om)$ be a compact
symplectic smooth manifold of dimesion $2n$ with a {\it smooth}
Hamiltonian action of $n$-dimensional torus $T^n$. Then
action-angle coordinates define the moment map
\begin{equation}
 \pi \colon M \to \De \subset \R^n
\end{equation}
whose image is a convex polyhedron $\De$ in Euclidean $n$-space.
This polyhedron satisfies so called Delzant conditions and
contains complete information on the symplectic geometry of $(M,
\om)$. That is, $\De$ determines the manifold, the symplectic form
and the $T^n$-action (see \cite{D}). Moreover, if a polyhedron
$\De \subset \R^n$ satisfies the Delzant conditions then we can
construct a smooth symplectic manifold with Hamiltomian
$T^n$-action such that the image of the moment map is precisely
$\De$.

Moreover, if a phase space $(M, \om)$ is prequantized
and $M$ has a Hodge structure  whose Kahler form is $\om$, then
this Hodge structure can also be reconstructed from $\De$ (see
\cite{G1}).

We saw that our space $CLRep(\pi_1(\Si), \SU(2))$ admits the
Hamiltonian $T^{3g-3}$-torus action, but
\begin{enumerate}
\item our phase space $CLRep(\pi_1(\Si), \SU(2))$ is singular if $g>2$;
\item the action isn't smooth over the boundary of $\De$.
\end{enumerate}

Nevertheless the polyhedron satisfies  Delzant conditions hence
determines the smooth compact Hodge manifold $DM_\Ga$ with the
smooth $T^{3g-3}$-action which is called  {\it Delzant model of
moduli spaces} (see \cite{T8}). For a complex Riemann surface $\Si_I$ the moduli
space $M^{ss}(\Si_I)$ is an algebraic variety. Every such variety
admits a birational morphism
\begin{equation}
 f \colon M^{ss}(\Si_I) \to DM_\Ga.
\end{equation}

The list of properties is the following
\begin{enumerate}
\item The Delzant model $DM_\Ga$ is canonically polarized by the line bundle
$H$ corresponding to the theta polarization of every $M^{ss}$.
\item  The Delzant model is rigid as a complex manifold and hence
the Kahler quantization of $(DM_\Ga, H)$ is succesful. The wave
function spaces are given as spaces of holomorphic sections
\begin{equation}
 H^0(DM_\Ga), H^k)
 \end{equation}
and dimensions of these spaces can be computed by the
Duistermaat-Heckman formula.
\item Bohr-Sommerfeld fibers for $CLRep(\pi_1(\Si), \SU(2))$ and $DM_\Ga$
as sets of points of $\De_\Ga$ coincide.
\item The rank of wave spaces (4.49) is equal to the number
 of Bohr-Sommerfeld fibers
thus the quantization of the dynamical system $(DM_\Ga, \om)$ is
numericaly perfect (see \cite{T8}).
\end{enumerate}

\section{Two versions of CQFT}

Result of any quantization procedure is  creation of a quantum
mechanics which is  QFT in dimension d = 1 when there is only one
variable - a time.
 But our main geometrical object - Riemann surfaces - are also the basic object
of  conformal field theories in dimension 2 (or, more precisely, (1,1)).

 Some years ago the constraint of conformal invariance was proposed as
the fundamental principle of quantum field theory. All two-dimensional
 conformal field theories are duing to the theory of
critical phenomena and to the string theory. So, the new  powerful
mathematical tool appeared since large data on the geometry of
Riemann surfaces is encoded in this rich mathematical structure.
Moreover, the results of a succesful Kahler quantization procedure
= the wave function spaces $H^0(M^{ss}(\Si), \Oh ( k \Theta))$
appeared first as spaces of conformal blocks of CQFT.

On the
other hand,
 the spaces
 $\sH_\pi^k = \oplus_{b \in BS_k \subset B} \C \cdot
s_b$ as  results of the Bohr-Sommerfeld quantization procedure
 appeared as fibers of projective flat vector bundles over
the moduli space $\sM_g$ of curves of genus $g$ in the
combinatorial version of CQFT in
 \cite{FS} and \cite{MS}. Moreover, a 2-dimensional conformal field theory
was formulated  as algebraic geometry of the moduli spaces of
Riemann surfaces.

A very important case of that a pinched Riemann sphere was
investigated by Belavin-Knijnik-Zamolodchikov in \cite{KZ} which
produce  the general ideology as well as thechnical fondations.

We can start with the simple "classical" version of CFT.

\subsection{WZW-version}

Recall that the space of classical theta functions has very
important  property to be a space of   unique irreducible
representation of some group or algebra (namely the Heisenberg
group). The space of holomorphic sections $H^0(M^{ss}(\Si), \Oh (
k \Theta))$ is also the space
  of the irreducible representation of the gauge algebra
of WZW RCFT. Such spaces have played a noted role in conformal
field theory.

It turn out that in  rational conformal field theories, one encounters the
"conformal blocks" of Belavin, Polyakov and Zamolodchikov coincidencing with
 our spaces  $H^0(M^{ss}(\Si), \Oh ( k \Theta))$ .

 First of all recall Wess-Zumino-Witten CQFT
  (WZW for short) set up: let $t$ be
a local parameter in a point $p \in \Si$. Then
  the {\it loop algebra}
\begin{equation}
 sl_2^t = sl(2, \C)\otimes \C[t, t^{-1}]
\end{equation}
 is the {\it horisontal algebra} of WZW.

The central extention by the cocycle $k \cdot c$
\begin{equation}
\ov{sl_2^t} = sl(2, \C) \otimes \C[t, t^{-1}] \oplus \C \cdot e
\end{equation}
is the k-{\it gauge algebra of} WZNW and the bracket of two
elements of the loop algebra is given by the formula
\begin{equation}
[X \otimes f, Y \otimes g] = [X, Y] \otimes fg + k \cdot tr(XY)
Res_0(gdf)
\end{equation}

The theory of Kac-Moody algebras tells us that this algebra admits
 unique irreducible representation $\sH^k$ such that there
exists a non-zero vector annihilated by
 $sl(2, \C)\otimes \C[t]$. Remark that the central element
  of $\ov{sl_2^t}$ acts by multiplication by $k$ (which is the
  level number).

On the other hand, the field $\C (\Si)$ of meromorphic functions
 on $\Si$  has the natural embedding:
\begin{equation}
i \colon sl(2, \C) \otimes \C (\Si) \to \ov{sl_2^{t}}.
\end{equation}
This field $\C (\Si)$ contains the subring $\C (\Si)_p$ of
regular functions in $\Si - p$. The structure of algebra on
\begin{equation}
sl(2, \C) \otimes \C (\Si)_p
\end{equation}
is given by the formula
\begin{equation}
[X \otimes f, Y \otimes g] = [X, Y] \otimes fg
\end{equation}
The restriction of the embedding (4.4) to this subring realizes it as the
subalgebra of $\ov{sl_2^t} $.

Thus the algebra $sl(2, \C) \otimes \C (\Si)_p$ acts on the
representation $\sH^k$
 described before.

 \begin{dfn} The space of conformal blocks $V_k$ of level $k$
 is the subspace of the dual
$(\sH^k)^*$ annihilated by the Lie algebra
 $sl(2, \C) \otimes \C (\Si)_p$.
\end{dfn}

Why such space is related to vector bundles on $\Si$? First of
all, we note, that every vector bundle $E$ on $\Si - p$ is
trivial and of course it is trivial on a
 small disc $D_p$ containing $p$. Fixing
these trivializations we send a vector bundle $E$ to the Laurent
series
\begin{equation}
\ga_E \in SL(2, \C)\otimes \C[t, t^{-1}].
\end{equation}
To get rid of the trivializations, we have to mod out by the
automorphisms groups of trivial bundles on $D_p$ and $\Si - p$.
So the set of vector bundles is
\begin{equation}
\{ E \} =   SL(2, \C) \otimes \C (\Si)_p \backslash  \fsl(2,
\C)\otimes \C[t, t^{-1}] / SL(2, \C)\otimes \C[t].
\end{equation}
Remark that we don't discuss the stability condition and our
 double coset space is so called {\it algebraic stack},
  not a moduli space.

In any case our, stack $\{ E \}$ is the quotient of the universal
homogeneous space
\begin{equation}
Q = SL(2, \C) \otimes \C[t, t^{-1}] / SL(2, \C)\otimes \C[t]
\end{equation}
and the line bundle $\Oh( \Theta)$ can be lifted to the universal
line bundle $L_\chi$ on this homogeneous space. The theorem of
Kumar and Mathueu provides the isomorphism
\begin{equation}
H^0(Q, L_\chi^k) = \sH^k
\end{equation}
and $H^0(M^{ss}, Oh(k \Theta))$ coincides with the subspace of
 $(\sH^k)^*$ invariant under Lie algebra
  $sl(2, \C) \otimes \C (\Si)_p$ action.

To see that this space is finite dimensional over $\C$, we need
only the standard experience with matrix divisors (for example,
like in the survey \cite{T2}). But to compute the dimension we
have to use more sophisticated versions of all objects.

First of all, on our Riemann surface $\Si$ we choose a finite set
of different points
\begin{equation}
 p_1, p_2, ....., p_n \subset \Si
\end{equation}
different from our fixed point $p \in \Si$ and labeled by the set of
irreducible representations of $sl(2, \C)$, that is, by positive
 integers $n_1, n_2, ...., n_n$ (twice spin-numbers).

Now our algebra $sl(2, \C) \otimes \C (\Si)_p$ acts on the space
$V_{n_i}$ of the irreducible representation $n_i$ by the formula
\begin{equation}
(X \otimes f) v  = f(p_i) X v
\end{equation}
where $p_i$ is the point corresponding to the representation
$n_i$ and $v \in V_{n_i}$.

New vector space is given now by the formula
\begin{equation}
V_k (\{ (p_i, n_i) \})_\Si = Hom_{sl(2, \C) \otimes \C (\Si)_p}(
\sH^k, \otimes_{i=1}^n V_{n_i}).
\end{equation}

These spaces have a number of good properties.
 For example, they don't depend on moduli of curves, the choice
  of points and so on. But the most important property is the
  behavior with respect to a simplest degeneration of
   curve $\Si$ to  curve $\Si_0$ with one double point.
   Geometrically this means that we have a smooth curve $\Si_0$
   of genus $g-1$ with a fixed pair of points $(q_+, q_-)$.
   Then
\begin{equation}
V_k (\{ (p_i, n_i) \})_\Si = \oplus_{n \in \Z^+} V_k (\{ (p_i,
n_i), (q_+, n), (q_-, n) \})_{\Si_0}.
\end{equation}
These properties provide the induction by genus (just like formula (3.107)
).

Now we are interested only in ranks of previous spaces.
 We may forget about curves, points and others geometrical
  features (how we did after (3.109) in subsection 3.7) and consider the monoid $M$ of formal sums
\begin{equation}
N = n_1 + n_2 + ... + n_n \in M
\end{equation}
providing the  present labeling of points. Every rank of space $V_k (\{ (p_i,
n_i) \})_\Si$ can be done as a symbol
\begin{equation}
rk_g(N) , \quad \text{where}\quad  N =  n_1 + n_2 + ... + n_n \in
M.
\end{equation}
The list of properties of these numbers is the following:
\begin{enumerate}
\item $rk_0(0) = 1;$
\item $rk_g(N) = \sum_{n \in \Z^+}rk_{g-1}(N +n +n)$;
\item let $g = g' + g''$, then
\begin{equation}
rk_g(N' + N'') = \sum_{n \in \Z^+} rk_{g'} (N' + n)\cdot
rk_{g''}(N'' + n).
\end{equation}
\end{enumerate}

Using these equalities, we can decent genus $g$ to 0 and  get so
called {\it fusion rules}:
\begin{enumerate}
\item $rk_0(0) = 1$;
\item $rk_0(N' + N'') = \sum_{ n \in \Z^+} rk_0(N' + n) rk_0(N'' + n)$.
\end{enumerate}

These fusion rules give the structure of commutative ring on our
monoid $M$:
\begin{equation}
 N \cdot N' = \sum_{ n \in \Z^+} rk_0(N + N' + n) \cdot n.
\end{equation}

This commutative ring is called the  {\it fusion ring}.

\begin{rmk} We don't prove here the factorisation rules (5.17) in the algebraic
geometrical set-up, that is, in the set-up of vector bundles. We
will get them as bypass results at the end of our long story. But
remark that there are exist approaches to prove these rules in
purely geometrical way.
\end{rmk}

We saw that our fusion ring is very close to be  the ring of
irreducible representations $R(sl(2, \C))$ where the multiplication
is  the tensor product of representations. But there is
 a small difference coming from level $k$ condition.

 According to the Clebsh-Gordan rules, a tensor product of
  irreducible representation has the decomposition:
\begin{equation}
n_1 \otimes n_2 = (n_1 + n_2) \oplus  (n_1 + n_2 - 2) \oplus ....
\oplus (n_1 - n_2)
\end{equation}
(here $n_1 \geq n_2$). But the level condition is $n \leq k$ and
 the multiplication is given by the formulas:
\begin{enumerate}
\item $ n_1 + n_2 \leq k \implies n_1 \circ n_2 = n_1 \otimes n_2$;
\item $n_1 + n_2 \geq k \implies n_1 \circ n_2 = (2k - n_1 - n_2) \oplus ... \oplus (n_1-n_2)$.
\end{enumerate}
From this it is easy to see that the fusion ring $M_k$ of level
 $k$ is the quotient of
 $R(sl(2, \C))$ by the ideal generated by $(k+1)$.

The fusion ring $M_k$ has the element
\begin{equation}
c = \sum_{n=0}^k n \circ n.
\end{equation}

Let $Spec M_k$ be the set of characters (homomorphisms to $\C$)
of $M_k$. Then after elementary computations we get the equality
\begin{equation}
rk_g(\emptyset) = \sum_{\chi \in Spec M_k} \chi(c)^{g-1}
\end{equation}
and
\begin{equation}
\chi(c) = \sum_{n=0}^k \vert \chi(n) \vert^2.
\end{equation}

Recall that every character of the representations ring $R(sl(2,
\C))$ is given by $z \in \C$ and
\begin{equation}
\chi_z(n) = \frac{sin((n+1)z)}{sin z}.
\end{equation}
This character is from $Spec M_k$, that is, it vanishes on $(k+1)$
iff
\begin{equation}
z = \frac{n \pi}{k+2}, \quad (1 \leq n \leq k+1).
\end{equation}
Now elementary computations and formula (8.35) give
\begin{equation}
rk_g(\emptyset) = rk V_k =\frac{(k+2)^{g-1}}{2^{g-1}}
\sum_{n=1}^{k+1}
 \frac{1}{(\sin(\frac{n\pi}{k+2}))^{2g-2}}
\end{equation}
 which is  the Verlinde number again.

 In the conformal blocks set up  very important observation is
  the absence of Hermitian structure on all vector spaces of this
  theory. Indeed we are working always with "sl" objects, not
   with "su". This feature is distiguishing from the quantization
   set up where Hermitian structure  appears manifestly.

\subsection{WZW-connection}

On the moduli space $\sM_{g,1}$ of  curves of genus $g$ with a fixed
(pinched) point  we
have vector bundles (1.29) interpreted as vector bundles of
conformal blocks. This identification induces another interpretation of the
Hitchin's  projective flat connection coming from CQFT (see
\cite{TK} or  beautiful survey  \cite{S}).

The construction and investigation of the connection in these
terms are almost parallel to the constructions described in
subsection 3.8.
 The difference is following: instead of the moduli space of
Higgs bundles (which is nothing else as the space of $SL(2,
\C)$-representations of the fundamental group of  Riemann
surface) we consider the infinite dimensional space $Q$ (5.9) and
lift
 all of our problems to it. Actually before that we have to develop a
 technique of working with
 infinite dimensional objects of such type rigorously. This was done
in \cite{BL}. Saveing the time-space here we describe this
approach only briefly.

Recall that the ind-homogeneous ("ind" is a term from \cite{BL})
 space $Q$ (5.9) describes
$\SU(2)$-vector bundles trivialized on $\Si - p$ where $(\Si, p)
\in \sM_{g, 1}$. Of course we consider the dense open subset
\begin{equation}
Q_0 \subset Q
\end{equation}
corresponding to  semi-stable bundles. For this space the
following fact was proved in \cite{DS}: the natural map
\begin{equation}
f \colon    Q_0 \to M^{ss}(\Si)
\end{equation}
is a locally trivial torsor in the etale topology. (The proof of
this fact is contained also in  unpublished preprint of Beilinson
and Kazhdan.)

In this set up Yves Laszlo constructed a family of second order
differential operators $\theta (\dot I)$ over the family of
deformations of complex structures. This construction commutes
with the construction given by the composite morphisms $\phi_s$
(1.86) with the cohomology sequence of the exact triple  (1.81)
which gives
\begin{equation}
H^1(\Si, T \Si) \to H^0(S^2(TM^{ss}(\Si))) \to H_h^1(d_s).
\end{equation}
This lifting is, in a  sense, a lifting of the Sugawara tensor
$T(\dot I)$ (modulo the ind-scheme technique). Namely for every
function $f$ on $\sM_{g,1}$ and a local section $\tilde{e}$ of
the family of dual to (5.10) spaces, the WZW-connection is given
by the formula
\begin{equation}
\nabla ( f \otimes \tilde{s}) = \dot I(f) \cdot \tilde{s} +
\frac{1}{2k-4} \cdot T(\dot I) \tilde{s} \otimes f \quad mod
(\tilde{s} \otimes f).
\end{equation}
From this point the equality of both  connections is obvious
because the only non trivial term in the previous formula is the
Sugawara term just as in the formula (3.152). Thus the Hitchin's
connection is the WZW-connection indeed.

\subsection{Monodromy representations}

All of our projective and vector connections are on vector bundles
over the Teichmuller space $\tau_g$. The modular
 group
\begin{equation}
Mod_g = Diff^+ \Si / Diff_0 \Si
\end{equation}
 acts on the Teihmuller space $\tau_g$ such that
\begin{equation}
 \sM_g = \tau_g / Mod_g.
\end{equation}
Thus all of them define  projective or vector representation of
the modular group $Mod_g$
\begin{equation}
\rho \colon Mod_g \to PGL(N_k, \C)
\end{equation}
where $N_k$ is rank of the vector bundle $V_k$.

For example, for the classical theory of theta functions we have
two vector representations
\begin{equation}
Mon_w^k \colon Mod_g \to PGL(k^g, \C)
\end{equation}
given by the  Welters projective flat connection, that is, by
applying the full procedure of extended Kodaira-Spencer theory
(see \cite{W}). On the other hand, the symplectic theory
 yields another combinatorial connection (see \cite{K1}
 and another monodromy representation
\begin{equation}
Mon_{BS}^k \colon Mod_g \to PGL(k^g, \C).
\end{equation}

Fortunately, the perfect quantization identifies both of these
connections and representations to give unique abelian
representation
\begin{equation}
Mon_w^k = Mon_{BS}^k = Mon^k_a.
\end{equation}
Moreover,  the monodromy group of this representation is finite
\begin{equation}
Mon_k^a = Sp(2g, \Z_{2k}).
\end{equation}

For non-abelian case we have the representation
\begin{equation}
Mon_H^k \colon Mod_g \to PGL(rk_g(\emptyset),\C)
\end{equation}
 given by the Hitchin connection.

On the other hand, using identifications from \cite{TUY} and \cite{TK}
 Kohno constructed this
representation in the pure combinatorial way
\begin{equation}
Mon_c^k \colon Mod_g \to PGL(rk_g(\emptyset),\C).
\end{equation}
 (The subscript
$"c"$ here stands for "combinatorial".) The lines decompositions (4.39)
of fibers
 is parallel to the decomposition of the highest
weight representation of the affine Lie algebra of $\fsl(2,\C)$
by eigenspaces of the operator $L_0$ from Sugawara construction
of the representation of the Virasoro Lie algebra (see
\cite{K2}). Manifestly this representation is projective
Hermitian. It  was investigated by Macsbaum and Funar. The main
result is the following: the monodromy group of this
representation
 is infinite if $k \neq 2,3,4,6$ and $k \neq  10$ for $g=2$.

Moreover, in paper \cite{An} J. Andersen proved {\it asymptotic 
faithfulness} of the collection of representations that is the intersection
 of kernels
$$
\bigcap_{k=1}^\infty ker Mon^k_H = 1 \text{ or } \Z_2 \text{ for } g = 2.
$$
(Of course this result is proved for any rank of vector bundles.)
More precisely he proved that that Berezin-Toeplitz endomorphisms associated to smooth functions on $M^{ss}$ are asymmptotically flat.

The appearance of  the Sugawara construction in the procedure of
 definition of this connection suggested the equivalence of
these representations. It was proved using the
 Borthwick-Paul-Uribe method in
 \cite{T3}. This method is quite different from CQFT and it will be
  productive to use generalized Sugawara construction \cite{Sc}
 to get this coincidence directly
 as the generalization of Kniznik-Zamolodchikov construction.

\section{3-valent graphs}

\subsection{Spin networks}

 Friedan and Shenker \cite{FS}, Moor and Seiberg
\cite{MS} proposed a new approach to CFT - the {\it axiomatic
description of a rational
 conformal field theory}. We will call it the MS-axiomatic. The main ingredient
 of this theory is 3-valent graphs which  appeared in constructions of
 section 4  already. We saw that every 3-valent graph is equivalent
 to a trinion-decomposed Riemann surface (see (4.9)-(4.11)).

In a similar way, we can extend these constructions to an
orientable surface $\Si$ of genus $g$ with $n$ holes.
Geometrically such surface with a trinion decomposition is
equivalent to its dual 3-1-valent graph $\Ga$ or a history
(3-valent graph of
 genus $g$ with $n$  parabolic edges which we describe in the special
 subsection).
In this case admissible integer weights are defined in the same
way as in (4.25)
 and (4.26) but the set of weights
 \begin{equation}
 W^k(\Ga, v_1, \dots, v_n)
 \end{equation}
is the set of functions $w$ (4.25)-(4.26) satisfying the
conditions
 \begin{equation}
 w(a_i) = v_i
 \end{equation}
for every parabolic edge $a_i, \quad 1 \leq i \leq n$, that is, the
values of these functions on parabolic edges are  fixed constants.
So in the same vein we have spaces
 \begin{equation}
\sH^{k}_{\Ga, v_1, \dots, v_n} = \oplus_{w \in W^{k}(\Ga, v_1,
\dots, v_n)} \C \cdot w
\end{equation}
 The {\it
Clebsch--Gordan conditions} (1), (2) and (3) from (4.26) provide a  very
 surprising
point of all these constructions: namely, let us multiply all our integer
 weights $w \in W^{k}(\Ga)$ by level $k$. We get the maps
\begin{equation}
 k \cdot w = j \colon E(\Ga) \to 1/2 \Z_{\geq 0}.
 \end{equation}
But we can consider the target set  as the set of irreducible
representations of $\SU(2)$:
\begin{equation}
\widehat{\SU(2)} = 1/2 \Z_{\geq 0}.
\end{equation}
So every integer weight  can be considered as a coloring of edges by
irreducible representations of $\SU(2)$. The conditions (1), (2)
and (3) (4.26) give that for every triple $e_{1}, e_{2}, e_{3}$
via one vertex $v$ the tensor product of representations
\begin{equation}
0 \in j(e_{1}) \otimes j(e_{2}) \otimes j(e_{3})
\end{equation}
that is this tensor product contains the trivial representation.
Or the tensor product $j(e_{1}) \otimes j(e_{2})$ admits an
intertwiner
\begin{equation}
i_{v} \colon j(e_{1}) \otimes j(e_{2}) \to j(e_{3})
\end{equation}
as  $\SU(2)$-equivariant homomorphism.

Such coloring of graph edges is called  $\SU(2)$-{\it spin
network} by Penrose \cite{P}( see also \cite{T1}). Thus the set
$W^{k}(\Ga)$ is the set of spin networks of level $k$ over a graph
$\Ga$. We will see in the following section that for every
integer weight $w \in
 W^{k}(\Ga)$ the function
\begin{equation}
 2k \cdot w = t \colon E(\Ga) \to  \Z_{\geq 0}.
 \end{equation}
is the type of a loop  on our graph $\Ga$. It means, that   a
colouring can be defined purely in terms of the geometry of graph
as a {\it topological space}.

The main fact here is the identification of the set $SNW^k(\Ga)$
of $\SU(2)$-spin networks  of level $k$
$$
SNW^k_\Ga = W^k(\Ga)
$$
and the set of integer weights $W^k(\Ga)$.

Our first task  is
to compare spaces $\sH^{k}_{\Ga}$
 for different graphs of genus $g$.
As a result of this comparing we  will construct, following Kohno
\cite{K1}, the holomorphic vector bundle with a holomorphic
projective flat connection over the moduli space $\sM_{g}$ of
Riemann surfaces of genus $g$. This construction is purely
topological. So first of all we have to construct combinatorial
models of the Teichmuller space and the trivalent graphs complex.
These ideas and constructions came from low dimensional topology
and classical and quantum Conformal Field Theories.

\begin{figure}[tbn]
\centerline{\epsfxsize=3in\epsfbox{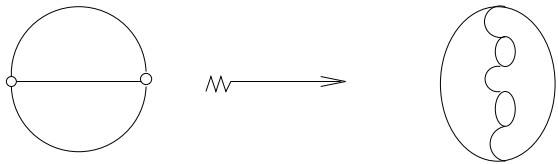}} \caption{\sl Pumping
up trick} \label{Fig 2}
\end{figure}

\subsection{3-dimensional topology}

 Recall again that any 3-valent graph is equivalent to  a trinion decomposition
 of a Riemann surface. It is easy to see that such surface is a boundary
  of uniquely determined  handlebody.
   We explain this fact more carefully.

Returning to  a maximal set of disjoint inequivalent curves $ \{
C_1, ... , C_{3g-3} \} $ (4.7) and removing these curves we get
a finite set of trinions $ \Si - \{ C_1, ... , C_{3g-3} \} =
\cup_{i=1}^{2g-2} \tilde{v_i} $ (4.8) (see Fig.2).

Every oriented trinion defines the handlebody inside  it and
gluing (sewing) these handlebodies together we get the handlebody
$H_\Ga$ such that
\begin{equation}
\p H_\Ga = \Si_\Ga
 \end{equation}
 Any diffeomorphism $h \in Mod_g$ (where $Mod_g$ is the modular group (5.30))
  gives us a compact  3-manifold
 \begin{equation}
X_{\Ga} = H_\Ga \cup_h  H_\Ga
 \end{equation}
where we are gluing two handlebodies with the boundaries by the
diffeomorphism
 \begin{equation}
h \colon \p H_\Ga = \Si_\Ga \to \Si_\Ga = \p H_\Ga.
 \end{equation}
 For the 3-manifold such presentation is called a {\it Heegaard splitting}.
 The subgroup
 \begin{equation}
M_\Ga \subset M_g
 \end{equation}
of diffeomorphisms which can be extended to diffeomorphisms of
the handlebody $H_\Ga$, acts on $Mod_g$ on the left and right and we
have the set of 3-manifolds
\begin{equation}
\{ X_\Ga \} = M_\Ga \backslash Mod_g \slash M_\Ga
 \end{equation}
as the set of double cosets.

 In  Heegard splitting (6.10) let us add a new handle to the surface and extend
 the diffeomorphism $g$ by the identity on this new handle. This procedure
  doesn't
 change our 3-manifold and is called  {\it stabilization}.
 Singer proved in \cite{Si} the following statement

 \begin{prop} The equivalence relation between elements $h \in Mod_g$ when $g \to
 \infty$ creating homeoomorphic 3-manifolds, is generated by isotopies, left
 and right
 actions of $M_\Ga$ and a stabilization.
 \end{prop}

 To construct an invariant of 3-manifolds we have to send a marked
 Riemann surface $\Si_\Ga$ to the wave function vector space $\sH^k_\Ga$
  (4.43) with fixed vector $w_{H_\Ga}$ given by the function (4.25) such that
  \begin{equation}
w_{H_\Ga}(e) = 0 \quad \text{for every} \quad e \in E(\Ga).
\end{equation}
Here
\begin{equation}
w_{H_\Ga} \in \sH^k_\Ga
\end{equation}
is called the {\it vacuum vector} for the handlebody $H_\Ga$.

 The space $\sH^k_\Ga$ admits the standard Hermitian structure which identifies
 this space with the dual.

 Now if we have a linear representation
 \begin{equation}
\rho \colon Mod_g \to GL (\sH^k_\Ga)
\end{equation}
 then
 \begin{equation}
I(X_{\Ga, g}) = \frac{< w_{H_\Ga}, \rho(g)(w_{H_\Ga})>}{\Vert
w_{H_\Ga} \Vert^2}
\end{equation}
has to be an invariant of $X_{\Ga}$ as 3-manifold since it
easy to see that
\begin{enumerate}
\item the vector $w_{H_\Ga}$ is invariant with respect to $M_\Ga$-action and
\item under normalization (6.17) this vector is invariant under stabilization.
\end{enumerate}

But the situation is a slightly more complicated then just
described. Namely, Kohno constructed nonlinear representation
$\rho$ but made it slightly projective,
 that is, determined up to the usual phase amplitude
 \begin{equation}
e^{\frac{\pi i \cdot k}{4 (k+2)} }.
\end{equation}
Thus the invariant (6.17) is determined up to the same phase
amplitude.

To construct such slightly projective representation (6.16) one has to
treat combinatorial constructions in the next special section.

\subsection{Geometry of graphs}

 A  graph is  one
dimensional complex containing edges -
 cells of one dimensional skeleton and vertices -cells of
 zero dimensional
 skeleton.

The combinatorial invariant of a pinched Riemann surfaces with
 a trinion decomposition
is a (3-1)-valent graphs or a {\it history}. Such graphs have
3-valent and 1-valent vertices. The last vertices are half
intervals with one vertex. Vertices of such edges are called
{\it parabolic}.

  Graphs and histories are Feynman diagrams - objects of the standard
  formalism
 in fields (for 3-valent graph case, $\phi^3$ - theories).

\begin{figure}[tbn]
\centerline{\epsfxsize=3in\epsfbox{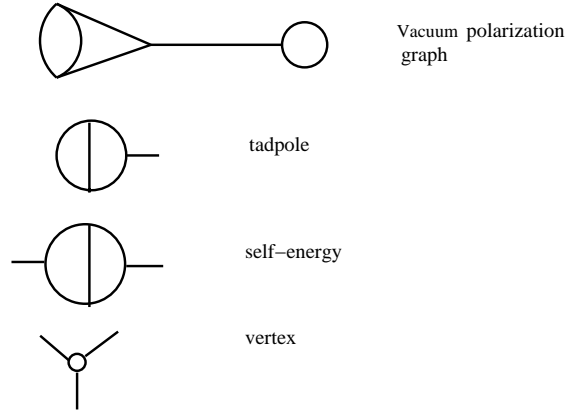}} \caption{\sl Standard
Feynman diagrams} \label{Fig 2}
\end{figure}

In classical and quantum field theory a 3-valent graph without
parabolic vertices is called a vacuum-polarization graph.  So we will
consider such
graph as a geometrical  object. Note that
 an closed 1-valent graph is a disjoint finite union of compact
  intervals and
 a closed 2-valent graph is a disjoint finite union of compact
 circles.

Every graph $\Ga$ has the following ingredients:
 \begin{enumerate}
\item $E(\Ga)$ is the set of edges;
 \item  $V( \Ga)$ is the set of vertices,
 \item  $S(v)$ is the star of a vertex $v$ that is the set of local edges incident  to  $v$,
 \item   $F(\Ga)=\{v\in e \}$ is the set of flags (edge + vertex).
 It is the same as
the set of oriented edges $\vec E(\Ga)$,
\item   $L( \Ga)$ is the set of circles, that is, of edges with only one vertex.
\end{enumerate}
We will consider homogeneous graphs that is the graphs of the same valence $d$ at every vertex.

For example for 3-valent graphs two projections   $p_e \colon F(\Ga)\to E(\Ga)$ and
 $v \colon F(\Ga)\to V(\Ga)$ are finite covers
 of degrees
  2 and 3 with the same ramification locus
  $W_e = W_v \subset F(\Ga)$, containing pairs
  $v\in e$, where   $e$ is a loop:
\begin{equation}
R_e = L( \Ga) \subset E( \Ga) ; \quad R_v = L( \Ga) \subset V(
\Ga)
\end{equation}
Let $\vert\ \vert$ be the cardinality of a finite set, then
 \begin{equation}
2 \cdot \vert E(\Ga) \vert-\vert L(\Ga) \vert=\vert F(\Ga) \vert=
3 \cdot \vert V(\Ga) \vert-\vert L(\Ga)\vert.
 \end{equation}

 Thus $ \vert V(\Ga) \vert=2g-2$ and $ \vert E(\Ga)
\vert=3g-3$, where $g>1$ is an integer number called  {\it genus}
of
 $\Ga$.

There exists only finite set $TG_{g}$ of 3-valent graphs of genus
$g$. But up to now nobody knows how many graphs it contains
\begin{equation}
\vert TG_{g} \vert = ?
\end{equation}
(besides of small $g = 2, 3,\dots, 11 $).

Consider any graph $\Ga$ with vertices of any valence $ > 2$. Let $e \in E(\Ga)$ be 
any edge with distinct two vertices $ \p e = v_1 \bigcup v_2$ and $v_1 \neq v_2$. Remove the edge $e$ and identify vertices $v_1$ and $v_2$ to get a single new vertex
 $v_{new}$. We get a new graph $\Ga'$ which is called  {\it contraction} of $\Ga$
 along the edge. Such contraction decreas $\vert E(\Ga) \vert$ and   $\vert V(\Ga) \vert$ by one. The minimal under this procedure graph $\Ga_{min}$ has only one vertex.  Obviously every graph can be obtained from a 3-valent graph by applying a chain of contractions. The inverse operation is {\it expansion} which is the blow up of a vertex $v$ with  $\vert S(v) \vert > 3$ to a new edge $e_{new}$ with
 two new vertices $v_{new}, v'_{new} = \p e_{new}$ and partitions
$$
(S(v_{new} - e_{new}) \bigcup  ( S(v'_{new}) - e_{new}) = S(v).
$$
 The maximal under this procedure graph is a trivalent graph $\Ga_{max}$.

Every 3-valent graph $\Ga$ has {\it closed relatives}: let $v, v'
$ be vertices of one edge $e$ and $e_{1}, e_{2}, e$ are edges via
$v$ and $e'_{1}, e'_{2}, e$ are edges via $v'$. Let us contract
the edge $e$ to obtain 4-valent vertex $v_{new}$. For this vertex the star $S(v_{new}) = e_{1} \bigcup e_{2} \bigcup e'_{1}
\bigcup e'_{2}$. Now divide this quadruple in two pairs, for
example, $(e_{1}, e'_{1})$ and $(e_{2}, e'_{2})$. Then  extend
the vertex $v_{new}$ to new edge $e_{new}$ such a way that we get
\begin{enumerate}
\item
two new vertices $v_{1}$ and
 $v_{2}$ such that $\p e_{new} = v_{1}, v_{2}$,
\item
 with stars  $S(v_{1}) = e_{1}, e'_{1}, e_{new}$ and
 $S(v_{2}) = e_{2}, e'_{2}, e_{new}$.
 \end{enumerate}
 We  get a new
   graph $\Ga'_{e}$. This is the composition of one contraction and one expansion
preserving valence of the graph.

   \begin{dfn} The transformation of graphs
   $$
   \Ga \to \Ga'_{e}
   $$
is called  {\it elementary transformation}.
   \end{dfn}
(We will see later that this transformation lies under {\it elementary move for  trinion decompositions.})

Sending the edge $e$ of $\Ga$ to the new edge $e_{new}$ of $\Ga'_e$ we get the identification 
\begin{equation}
  \na \colon E( \Ga) \to E(\Ga'_{e}).
   \end{equation}

\begin{figure}[tbn]
\centerline{\epsfxsize=3in\epsfbox{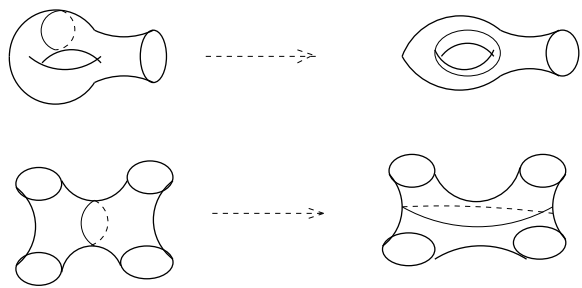}} \caption{\sl
Elementary transformations} \label{Fig 11}
\end{figure}

So, around any 3-valent graph there exists the collection of
3-valent graphs which differ from $\Ga$ by  elementary transformation. (This triple corresponds to the triple of vertices of the triangle of Fig. 13.)

By the induction on $g$ it is easy to prove that any 3-valent
 graph $\Ga$ of genus $g$ can be transformed to any other
 3-valent graph $\Ga'$ by a finite chain of elementary transformations.
 Thus the elementary transformations act on $TG_{g}$  transitively.

In the same vein we can consider the composition of two contractions and two
 expansions preserving valence of verticies and so on. 

Now for any graph $\Ga$ let 
$$
\sE(\Ga) = \{ \Ga_{max} \}
$$
be the set of all expansions of $\Ga$ to a trivalent graph. This set of 3-valent graphs we call  {\it nest} of 3-valent graphs.

Of course we can consider compositions of expansions for every vertex $v \in V(\Ga)$ with $\vert S(v) \vert > 3$ independentely. So for every such vertex
the set of expansions $\sE(v)$ is in 1-1 correspondence with the set  of 3-valent trees $T$ with parabolic edges identifyed with the set $S(v)$:
$$
\p T = S(v).
$$
Thus 
$$
\vert \sE(v) \vert = \sum_{T \vert \p T = S(v)} \frac{1}{\vert Aut T \vert}.
$$

 So for a graph $\Ga$ with one 
4-valent vertex 
$\vert \sE(\Ga) \vert = 3 $
 and so on.

The combinatoric behind of these constructions can be described perfectly  by the standard pairing scheme in the {\it Feynman diagram technique}.
This technique  was applyed  as well to so called
{\it ribbon graphs}. Recall that a ribbon graph structure on a graph is  cyclic
 orders choosen  on all stars of edges. Geometricaly it is equaivalent to an embedding of a graph into an oriented surface. Such cyclic orders are required by applications to  matrix models (see \cite{MP} and \cite{Mu}). 

For such graphs contractions and expansions preserve the ribbon structures. For a homogeneous d-valent graph $\Ga$ a ribbon graph structure is equivalent to the choice of a {\it connection} on this graph (see \cite{BGH} for the consideration of graphs as manifolds). 

Recall that a path of lenth 1 on  $\Ga$ is just an oriented edge
 $\vec e$. Let the set
$P_1(\Ga) = \vec E(\Ga)$ be the set  of  1-paths on  $\Ga$.
Every such path  $\vec e $ has  vertices of two types
\begin{equation}
v_s(\vec e), \quad v_t(\vec e) \in V(\Ga)
\end{equation}
 {\it source} and  {\it target} which are equal for  loops.

A path of length  $d$ in  $\Ga$ is an ordered sequence $(\vec e_1,
..., \vec e_d)$ of oriented edges such that for every $i$
\begin{equation}
v_t(\vec e_i) = v_s(\vec e_{i+1}).
\end{equation}
Let $e_{d+1} = e_1$, that is, our order is cyclic. A path $(\vec
e_1, ..., \vec e_d)$     is called irreducible if   $e_i \neq
e_{i+1}  $
 for every
$i$ (including $i = d+1 $).

 A connection $\na$ on a homogeneous graph $\Ga$ is a collection of
identifications of the stars 
\begin{equation}
\na_{\vec e} \colon S(v_{s}(\vec e)) \to  S(v_{t}(\vec e)) 
\end{equation}
such that 
$$
\na_{\vec e} (e) = e
$$
where $v_{s}, v_{t} = \p \vec e$.

A path   $(\vec
e_1, ..., \vec e_n)$     is called {\it geodesic} if   
\begin{equation}
e_{i+2} =  \na_{\vec e_{i+1}} (e_i) 
\end{equation}
 for every
$i$ (including $i = n+1 $) when we have closed geodesic.

Obviously for every pair of edges from the star of a vertex there exists
 unique  closed geodesic passing through the edges. Hence the set of closed geodesics defines our connection
 uniquely. 

In the same vein we can say that a subgraph $\Ga' \subset \Ga$ is totally geodesic if every geodesic starting in $\Ga'$ stay within $\Ga$.

Now for every closed geodesic  $l = (\vec
e_1, ..., \vec e_d = \vec e_1)$ we have the set of vertices $V_l$ of the edges in this consequence. For every $v \in V_l$ we have  a monodromy permutation
$$
m_{v} (l) \colon S(v) \to S(v)
$$
So every closed geodesic $l$ defines a conjugation class 
$$
m(l) \in Conj(S_d)
$$
where $S_d$ is the group of permutations with $d$ elements.

A closed geodesic is {\it flat} if $m(l) = id$.

Let $G_\Ga$ be the set of simple closed irreducible geodesics of $\Ga$. 
Then the map
$$
c_\na \colon G_\Ga \to  Conj(S_d)
$$
is an invariant of  connection $\na$ and the preimage of the unit class is the subset $G_\Ga^F$ of flat geodesics.
Moreover this map is an invariant for gauge classes of connections:
for every $v \in V(\Ga)$ consider the group 
$S_v = Aut S(v) = S_d $
where the last equality is non canonical, of course, and the product
$$
\sG = \prod_{v \in V(\Ga)} S_v.
$$
This group acts on the set $\sA_\Ga$ of all connections on $\Ga$ by the formula:
for every $\vec e$ with $\p \vec e = v_{s}, v_{t}$ 
\begin{equation}
g (\na)_{\vec e} = g_{s} \na_{\vec e} g_{t}^{-1}
\end{equation}
where $g_v$ are components of the decomposition of our discrete gauge group.

Then the set
$$
\sB = \sA / \sG
$$
is the set of orbits with respect to  the gauge group  and the function $c_\na$ 
is an invariant of the orbits.

Now gluing  2-cell to every simple closed flat geodesic we get a 2-complex
$C_\na(\Ga)$  which is a closed oriented surface corresponding to this ribbon graph.

A ribbon graph is called {\it plane} if the surface $C(\Ga) = S^2$ is 2-sphere.

For ribbon graphs the set of expansions is "reducible". In particular we can
choose from any pair of graphs correlated by  an elemetary transformation one
 of them. Recall that in the procedure of expansion we identify two expanded graphs if there is a ribbon graph isomorphism from one to the other preserving all local edges. In particular if $v \in V(\Ga)$ is d-valent vertex then
$\vert \sS(v) \vert = \frac{1}{d-1} C^{2d-4}_{d-2}$ which is the Catalan number.

Now forget the ribbon structure:  
let  $P_d(\Ga)$ be the set of irreducible paths of length  $d$ in
$\Ga$. Every path
 $(\vec e_1, ..., \vec e_d) \in P_d(\Ga)$ defines two vertices
\begin{equation}
v_s((\vec e_1, ..., \vec e_d)), \quad v_t((\vec e_1, ..., \vec
e_d)) \in V(\Ga)
\end{equation}
- source and  target and two maps
\begin{equation}
v_s \colon P_d(\Ga) \to V(\Ga)
\end{equation}
$$
v_t \colon P_d(\Ga) \to V(\Ga)
$$
$$
(v_s \times v_t) \colon P_d(\Ga) \to V(\Ga) \times V(\Ga).
$$
The preimage of the intersection
\begin{equation}
L_d(\Ga) = (v_s \times v_t)^{-1} ( V(\Ga)_\Delta )
\end{equation}
is the set of oriented irreducible loops of length   $d$. In
particular
$$
F(L_1)(\Ga)= L(\Ga),
$$
where  $F$ is the projection of a path to the graph.  For every
vertex á
 $v \in V(\Ga)$ we have the set of irreducible loops marked by $v$
\begin{equation}
L_d(\Ga)_v = (v_s \times v_t)^{-1} ( v).
\end{equation}
The union
\begin{equation}
L_\infty(\Ga)_v = \cup_{d=1}^\infty L_d(\Ga)_v
\end{equation}
admits a group structure
\begin{equation}
\pi_1^C(\Ga)_v = L_\infty(\Ga)_v.
\end{equation}
if we consider only irreducible fragments of compositions.

Obviously, this group depends on the marking point $v$.

Let $\pi_1(\Ga)$ be the standard fundamental group of $\Ga$ (as
1-complex). Then the natural epimorphism  $r \colon \pi_1^C(\Ga)
\to \pi_1(\Ga)$ is an isomorphism. Obviously, if $\Ga$ is a 3-valent graph 
of genus $g$ then $\pi_1(\Ga) = F_g $
is  free group with $g$ generators.

For any 3-valent graph  $  \Ga  $  a choice of  two  edges  $e, e'$  (may be not
   different) defines a new graph  $\Ga_{e, e'}$  by the
following procedure: let us choose a pair of inner points
 $v \in e$ and   $v' \in e'$ and joint them by the edge
$e_{new}$. We get a new  3-valent graph  with  the pair of new
vertices  and the triple of new edges
 $e_{new}, e{1/2} , e_{ 1/2}$, where the last two edges are just halves of the edges
  obtained by dividing  previous edgesá $e, e'$.

The inverse  operation is   much more simple  :  we just remove  an
edge and a pair of its vertices.

\begin{figure}[tbn]
\centerline{\epsfxsize=3in\epsfbox{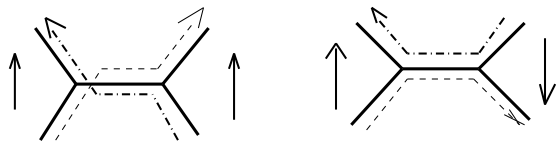}} \caption{\sl Two
paths} \label{Fig 3o}
\end{figure}
If  our  edges  were  in two  paths  $
  e \in p = (\vec e_{1}, ...,
\vec e_{n})$ and  $e' \in p' = (\vec e'_{1}, ..., \vec e'_{m})$,
then on the new
 graph  $\Ga_{e, e'}$ we would have (formally) two paths
$p_{new}, p'_{new}$.

It is easy to see  that
\begin{enumerate}
\item $ p \neq p' \implies p_{new} = p'_{new}$;
\item $ p = p' \implies p_{new} \neq p'_{new}$.
\end{enumerate}

By   induction one obtains:
\begin{enumerate}
\item Any   3-valent  graph  can be created by  a finite  chain  of such
elementary  operations starting with  the disconnected union of circles
 $0 \cup 0 ... \cup 0$;
\item  every path on $\Ga$  can be created by a finite chain  of
  elementary
 operations  starting from a path on a disjoint finite set of circles.
\end{enumerate}
These statements give the induction procedure which allows to prove  all simple
facts about trivalent graphs. There are two types of such problems: {\it Eulerian paths} and {\it coloring problems}.

 Recall that a path $p = (\vec e_{1}, ...,
\vec e_{n}) \in P_d(\Ga)$ is called {\it Eulerian} if it contains every oriented edge and exactly one time. A system of paths $P_E$ is called Eulerian if every oriented edge is contained and exactly one time in one and exactly one path from
this system. Such system is called minimal if $\vert P_E  \vert = min$. This is
an invariant of the graph $e(\Ga)$.

We know already that an elementary operation tranform an Eulerian system $P_E$
on $\Ga$ to the Eulerian system $P'_E$ on  $\Ga_{e, e'}$. Moreover
$$
\vert P'_E  \vert = \vert P_E  \vert \pm 1.
$$
From this it easy to see that for every 3-valent graph $\Ga$ of genus $g$
\begin{enumerate}
\item $1 \leq e(\Ga) \leq g+1$;
\item $e(\Ga) = g-1 mod  2$;
\item there exists a chain of elemetary operations with result graph $\Ga$ such that $e(\Ga') = 1$.
\end{enumerate}

We can consider every plane ribbon graph $\Ga$  as a graph on the plane $\R^2$.
Such embedded graph defines the geographic map with connected componets of the complement as the countries. This is the {\it 4-colors setup}. This problem is stated for embedded
 graphs.

Other "coloring" problem for graphs appears in the Feynman diagram technique
 (for the problem of dependence on the form of Lagrangian and of limits 
in matrix models). Let $C = (c_1, ... , c_N)$ 
is a set of colors. Any coloring is just a map
$$
c \colon E(\Ga) \to C
$$
such that for every vertex $v \in V(\Ga)$ the map $c \colon S(v) \to C$ is an
{\it embedding}. That is around every vertex all edges have different colors.
If the number $N$ of colors is minimal then  it is called the {\it chromatic number } $N_{min} = \chi (\Ga)$. Easely we have to consider graphs without 
1-loops. Then for 3-valent graphs
$$
3 \leq \chi(\Ga) \leq 4.
$$  
A little more efforts you need to prove the following fact

$\chi(\Ga) = 3$ if and only if there is a set $(l_1, ... , l_n)$ of disjoint
simple loops on $\Ga$ such that
\begin{enumerate} 
\item every loop $l_i$ contains an even number of edges;
\item for every vertex $v \in V(\Ga)$ there exists a loop which contains it:
$$
V(\Ga) \subset \bigcup_{i=1}^n l_i.
$$
\end{enumerate}

Both of these "coloring" problems are related by the following equivalence:
two  statements 
\begin{enumerate}
\item The 4-colors problem admits  positive solution;
\item For every plane 3-valent graph $\Ga$ without 1-loops $\chi(\Ga) = 3$
\end{enumerate}
are equivalent.
These problems of "classical theory" of graphs. Considering every homogeneous 
graph as a manifold (in style of \cite{BGH}) and following  Segal's ideas we 
have a reason
 to consider the "Loop space" of every graph.

Consider  standard interval $[0, 1]$ and divide it by a finite
set of points to a complex $[0, p_1, .... , p_{d-1}, 1]$. We can
consider this as a 2-valent history with two parabolic edges $[0,
p_1]$ and $[p_{d-1}, 1]$. The linking of the parabolic edges
 of such history gives  a cycle (or a loop) that is a closed 2-valent graph.

As we know every continuous combinatorial map
\begin{equation}
p \colon [0, p_1, .... , p_{d-1}, 1]\to \Ga
\end{equation}
 is called  path of length $d$ in $\Ga$.

If this path is irreducible, the preimage of every edge
\begin{equation}
      p^{-1}(e) = g_1 \cup ... \cup g_{col(e)}
\end{equation}
is a disjoint union of a finite set of compact 1-graphs and the
number of components is equal to $col(e)$ (without any orientation).

 Now we consider the set $\sL$ of irreducible loops that is closed
 paths (0=1). The type of  loop $l \in \sL $ is   the non negative
function
\begin{equation}
col \colon E(\Ga) \to \Z_{\geq 0}.
\end{equation}
(recall that we  consider irreducible unmarked loops only).

For any type of loops for every vertex $v \in V(\Ga)$ with the set
of edges
  $(e_{1,v}, e_{2,v}, e_{3,v})$  the following constraints hold:
 \begin{enumerate}
\item
\begin{equation}
 col (e_{1,v}) + col( e_{2,v}) + col( e_{3,v}) = 0 \mod 2;
\end{equation}
\item
\begin{equation}
|col( e_{1,v})-col(e_{2,v})| \leq col ( e_{3,v}) \leq col(
e_{1,v}) + col(e_{2,v}).
\end{equation}
\end{enumerate}
Indeed, consider a vertex $v$ and the star $S(v) = e_1, e_2, e_3$.
Let $n_{12}$ be the number equal to how many times  we come to $v$ along
 the edge $e_1$ and turn to the edge $e_2$. In the same vein we have
$n_{ij}$ for every pair of edges from the star. Then
\begin{equation}
 col(e_i) = \sum_{j \neq i} n_{ij}.
\end{equation}
From this we  get the statement immediately.

Now we say that a loop $l \in \sL$ has the type of level $k$ if
for
 every edge $e$
$$
 col (e) \leq k.
$$
In the same vein we can prove that the type function (6.36) gives
a type of
 loop of level $k$  iff
 for every vertex $v \in V(\Ga)$ with the set of local edges
  $(e_{1,v}, e_{2,v}, e_{3,v})$   almost the same constraint
   holds:
 \begin{enumerate}
\item
$$
 col (e_{1,v}) + col( e_{2,v}) + col( e_{3,v}) = 0 mod 2;
$$
\item  and a slightly  sophisticated  condition 2)
\begin{equation}
|col( e_{1,v})-col(e_{2,v})| \leq col ( e_{3,v}) \leq
\end{equation}
$$
 min (col( e_{1,v}) + col(e_{2,v}), 2k -
 col( e_{1,v}) - col(e_{2,v}).
$$
\end{enumerate}

Note that the inequality (6.38) is just the triangle inequality
and the inequality
 (6.40) is the {\it spherical} triangle inequality
 (with $\sqrt{k}$ proportional to the radius of the sphere).

Note also that it is quite reasonable to call every
loop  as {\it knot}.

So every knot defines the type function (6.36)  coloring edges of
our graph in $\Z_{\geq 0}$ colors. {\it Dividing this function by 2 we
get the map}
$$
1/2 col \colon E(\Ga) \to 1/2 \Z_{\geq 0}  = \widehat{\SU(2)}
$$
{\it where the last set is the set of irreducible representations of} $\SU(2)$. Conditions (6.38) and (6.40) mean that any pair $(\Ga, cal)$ is
a $\SU(2)$- {\it spin network}.

Recall that to get a marked Riemann surface from a trivalent graph $\Ga$ we
need to pump up our graph,
 that is, to replace every  vertex by a trinion and to glue
  trinions by tubes along the edges.

\subsection{ Circle complex and Graph complex}

  Let $\Si_{g,n}$ be an oriented surface of genus $g$ with $n$
  punctures and $C$ be a simple, closed, homotopy non-trivial,
   non-bounded pinchs loop on $\Si_{g,n}$. We call such loop a
     {\it circle} and  consider it up to the isotopy class.
     
{\it Circle
complex} or {\it Curve complex} $C(\Si_{g,n})$ has the set  of
circles as   0-skeleton $C_{0}(\Si_{g,n})$,
\begin{enumerate}
\item 1-skeleton $C_{1}(\Si_{g,n})$ is the set of
 pairs of classes of disjoint circles;
\item triangles correspond to triples of classes of disjoint cirles ;
\item every $i$-simplex corresponds to a set of $i+1$ classes
 mutually
disjoint circles,
\item so, simplexes of maximal dimension are collections of $3g-3$
 disjoint
circles just as (4.7).
\end{enumerate}
The main fact which we will use is the following: the group of
combinatorial automorphisms
\begin{equation}
Aut C(\Si_{g}) = Mod_{g}
\end{equation}
is the modular group (5.30).
\begin{rmk}
By the Roiden's theorem
\begin{equation}
Isom \quad \tau_{g} = Mod_{g}
\end{equation}
where the Teichmuller space $\tau_{g}$ (5.31) is equipped with
the Teichmuller metric. Moreover as a metric space $C(\Si_{g,n})$
is very closed to $\tau_{g}$. In particular, $C(\Si_{g,n})$ is
$\delta$-hyperbolic (\cite{MM}). The same facts are true for the
Weil-Peterson metrics.
\end{rmk}
Recall that in  the case of a pinched curve $\Si_{g, n}$ the circle
can't be isotopic to a boundary component.

\begin{figure}[tbn]
\centerline{\epsfxsize=3in\epsfbox{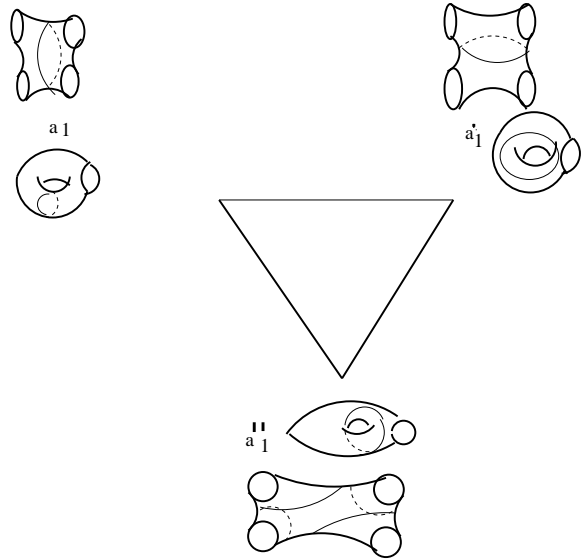}} \caption{\sl Triangles} \label{Fig 13}
\end{figure}

 For two special but very important  cases $\Si_{0, 4}$ and
  $\Si_{1,1}$ the
definition is slightly different
 (since there are no disjoint pairs of circles
): two vertices are connected by an edge if circles have 2
intersections in $\Si_{0,4}$-case and 1 intersection in
$\Si_{1,1}$-case. In both cases
\begin{equation}
C(\Si_{0,4}) = C(\Si_{1,1})
\end{equation}
is the ideal triangulation of  2-disc (so called Farey Graph,
see Fig. 14).

\begin{figure}[tbn]
\centerline{\epsfxsize=3in\epsfbox{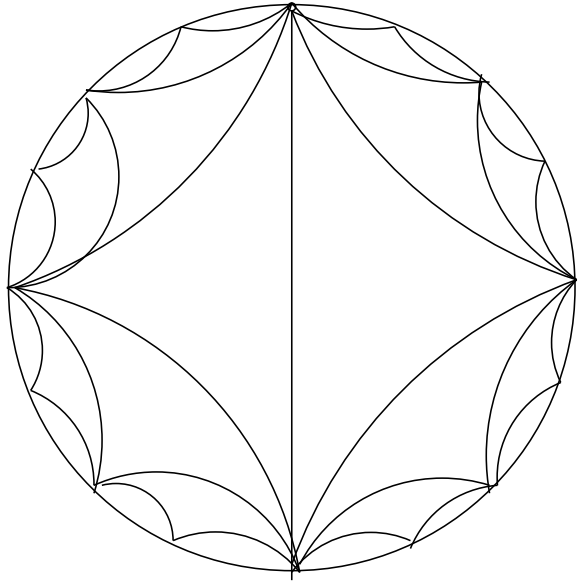}} \caption{\sl Farey
graph (the ideal triangulation of disk)} \label{Fig 13}
\end{figure}

  {\it Trinion decomposition
 graph } $C_1{TD(\Si)}$ has the set of trinion decompositions
  as 0-skeleton
\begin{equation}
V(C_{1}(TD(\Si))) = \{(C_{1}, ... , C_{3g-3})\}
\end{equation}
(see (4.7)) and a pair of trinion decompositions is connected by
an edge if the decompositions differ by an {\it elementary move}. That
is, in the system of circles $C_{1}, \dots, C_{3g-3}$ on a
Riemann surface
 $\Si$ we replace  only one circle, for example $C_{1}$, with
 another circle $C'_{1}$ intersecting twice if $C_{1}$ and
 $C'_{1}$ don't lie on the same trinion  and intersecting
 once otherwise (minimal intersections):
\begin{equation}
\{C_{1}, \dots, C_{3g-3}\} \to \{C'_{1}, \dots, C_{3g-3}\}.
\end{equation}
Recall that every trinion decomposition $C_{1}, ... , C_{3g-3}$
defines the 3-valent graph $\Ga$ (4.9)- (4.11). Thus we have the
{\it finite graph} $C_{1}(TG_{g})$ with the set of vertices
$$
V(C_{1}(TG_{g})) = TG_{g}
$$
(see (6.21)) and  a pair of graphs $\Ga, \Ga'$  is connected by an
edge if they differ by an elementary transformation 
(see Definition 9) corresponding to an elementary move of trinion
 decomposition (6.45). The fibers of the map
$$
\ga \colon V(C_{1}(TD(\Si)))  \to TG_{g}
$$
are described by the formula
$$
Mod_{g} / Aut H_{\Ga}
$$
where $H_{\Ga}$ is the handlebody (6.9) and $Aut$ is the group of
isotopy classes of diffeomorphisms of $\Si = \Si_{\Ga}$
extendable to diffeomorphisms of the handlebody.

The finite TG-graph $C_{1}(TG_g)$ is  1-skeleton  of the finite   $TG_g$-complex
$C(TG_g)$ which can be constructed by the following procedure: first of all let us glue by 2-cells all triangles coming from elementary transformations.
Over such triangles in TG-graph we have triangles from Fig. 13 of the infinite TD-graph. Up to such triangles we have the canonical connection $\na$ (6.22) on
 TG-graph and we can apply the "geodesic technique" (6.25)-(6.28). 
In particular,  we obtain  2-complex $C_\na(TG_g) = C(TG_g)$. The important
 fact (see Remark at the end of this subsection) proved in \cite{HT} is the following
\begin{enumerate}
\item $C(TG_g)$ is simply connected that is
$$
C(TG_g) = S^2,
$$
\item hence our TG-graph is plane;
\item hence we may ask about this graph all classical questions, for example,
on the chromatic number and so on.
\end{enumerate}

However  the standard way to check these properties is the following
\begin{enumerate}
\item to construct the infinite complex $C(TD(\Si))$ "by our hands";
\item to check described properties;
\item to desend all constructions to TG-complex using the map $\ga : C_1(TD) \to C_1(TG_g)$ and to get the complex $C(TG_g)$;
\item to check the required properties for this finite complex.
\end{enumerate} 

\begin{rmk}
It is quite natural to apply this geodesic technique to the infinite graph $C_1(TD(\Si))$ directly but we have to go passing the standard way.
\end{rmk}

 For our special cases  of Riemann surfaces it  can be proved that
\begin{equation}
C_{1}(TD(\Si_{0,4})) = C_{1}(TG(\Si_{1,2})) = C_{1}(\Si_{0,4})=
C_{1}(\Si_{1,2})
\end{equation}
is the Farey graph that is 1-skeleton of the Farey complex again.

\begin{figure}[tbn]
\centerline{\epsfxsize=3in\epsfbox{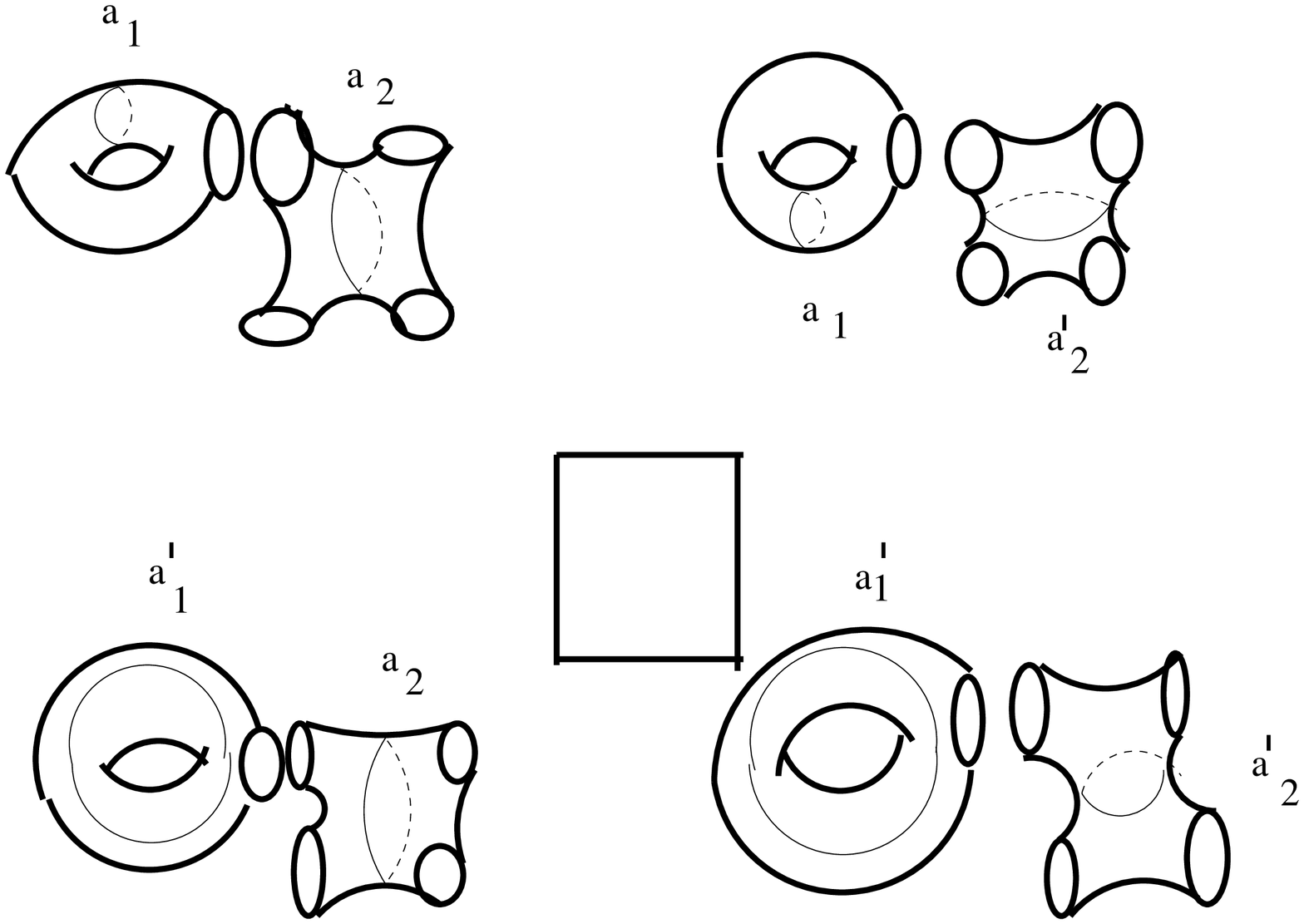}} \caption{\sl Squares} \label{Fig 13}
\end{figure}

The infinite complex $C(TD)$ is 2-dimensional. Its 2-cells
represent relations between special moves, that is, loops in
$TD$-graph
  and they can be triangles, squares,
 pentagones and hexagons. Every trinion decomposition of $\Si_{g,n}$
  contains trinion decompositions of surfaces with smaller $g$ and $n$ as
  fragments of it. We call such fragments subsurfaces.
 The relations,  coming from such subsurfaces are the following:

\begin{figure}[tbn]
\centerline{\epsfxsize=3in\epsfbox{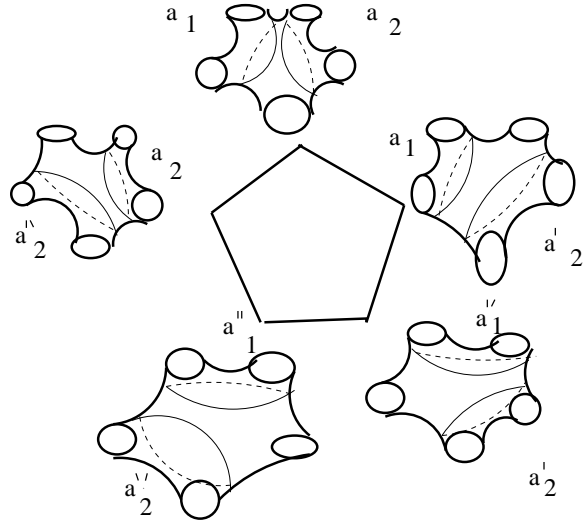}} \caption{\sl Pentagons}
 \label{Fig 13}
\end{figure}

\begin{enumerate}
\item $\Si_{1,1}$ and $\Si_{0,4}$ for triangles;
\item \begin{equation}
\Si_{1,1} \cup \Si_{0,4} = \Si_{1,3}
\end{equation}
for squares;
\item $\Si_{0,5}$ for pentagons and
\item $\Si_{1,2}$ for hexagons
\end{enumerate}
(see Figures 13, 15-17).

\begin{figure}[tbn]
\centerline{\epsfxsize=3in\epsfbox{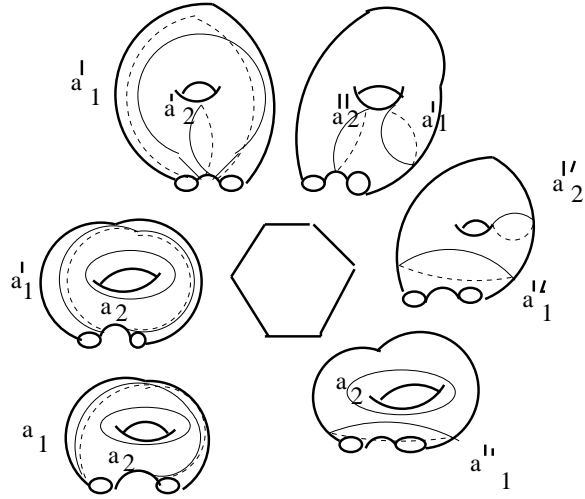}} \caption{\sl Hexagons} \label{Fig 13}
\end{figure}

Amazing fact is that we can define 2-cells of this complex in
terms of  $TD$-graph. It means that $C(TD)$-complex carries no
more information than its 1-skeleton:
\begin{enumerate}
\item triangles are contained in subcomplexes
 $C(TD(\Si_{1,1}))$ or  $C(TD(\Si_{0,4}))$ which are
 the Farey complexes.
Every edge of this complex  is a side of two triangles
corresponding
 to two subsurfaces
 $C(TD(\Si_{1,1}))$ or  $C(TD(\Si_{0,4}))$;
\item moreover, if we fix one vertex in any Farey subgraph of $TGD$-graph,
 we get
a circle which is a vertex of $C(\Si_{g,n})$ (indeed, fixing any
Farey subgraph of $TD$-graph we get the collection $\{ C_1, \dots
, C_{3g-3} \}$ of circles  of
 trinion decompositions with fixed set of circles
  $( C_2, \dots , C_{3g-3} )$,
 since our Farey subgraph corresponds to a choice of $C_1$ );
\item if $C$ is a vertex of  $C(\Si)$ (as $C_1$  before ) then there
exists a Farey subgraph with fixed vertex but not unique (indeed
we can extend
 our circle $C_1$ to a complete set  $\{ C_1, \dots , C_{3g-3} \}$
 and consider all
trinion decompositions with fixed  $( C_2, \dots , C_{3g-3} )$;
we get a Farey subgraph $F_{C}$  (just a possible choice) with
fixed point $p \in F_C$);
\item now we can define a map
\begin{equation}
\phi \colon Aut C_1(TD) \to Aut C_{0}(\Si_{g,n})
\end{equation}
sending a circle $C \in  C_{0}(\Si_{g,n})$ to  unique vertex in
$C_{0}(\Si_{g,n})$
 corresponding to the marked Farey graph  $g(p) \in g (F_C)$;
\item   square, pentagon or hexagon ($n$-gones for small $n$)
 can be  union of  triangle and
$n-1$-gon  gluied along an edge. To distinguish this trivial
situation from "non trivial" or {\it
 alternating}
 we consider a loop with $n$ consecutive points  such
that  every consecutive triple $p_{i},p_{i+1},p_{i+2} $ and the
pair of edges
 between  arn't contained in any Farey subcomplex;
\item at last, a hexagon can be {\it almost} alternating when its three
consecutive
vertices can be extended to the set of vertices of a square loop.
\end{enumerate}

Using these simple rules Margalit proved in \cite{M} that
\begin{enumerate}
\item the  map $\phi$ (6.47) can be extended to a
 well defined map
 \begin{equation}
 Aut C_1(TD) \to Aut C(TD),
\end{equation}
\item and this map can be extended to an isomorphism
\begin{equation}
\phi \colon Aut C(TD) \to Aut C(\Si),
\end{equation}
\item such that
\begin{equation}
 Aut C_{1}(TD)= Aut C(TD)=  Aut C(\Si)= Mod_g.
\end{equation}
\end{enumerate}
 TD-complex  $C(TD)$ was introduced by Hatcher and Thurston
 in \cite{HT} where they proved that this complex is connected and
  simple connected. But  the equality (6.50) was proved quite
  recently in \cite{M}.

 \begin{rmk} Returning to the Rational Conformal Field Theory
  WZW we can consider the vector space spanned by  different
  conformal blocks (see Definition 8) interpreted as the
  holomorphic sections of some power of theta-line bundle on the
  moduli space of semi-stable vector bundles on a Riemann
   surface $\Si$ (see section 3).
    This surface can be
  formed by gluing a number of trinions. The conformal
  blocks are obtained by summing over the intermediate states
  passing through the circles along which we glue our trinions.
  Different gluing procedures lead to different bases of the same
  vector space. Hence vectors of
  this space obtained by one gluing procedure can be expressed as
  linear combinations of vectors found by another way of gluing
   the same surface. These linear transformations are given by
   {\it dual matrices}, which  have to satisfy some
   {\it consistency conditions}. To describe these conditions it is
   quite convinient to use the finite complex $C(TG)$.
    First of all, we
   have some very simple duality matrices defined by the elementary moves. The complex $C(TG)$ is connected, hence every duality matrix can be represented as a
    product of the simple ones (maybe not uniquely). To get
   an  unambiguous description of all duality matrices, we have to
     use consistency conditions on simple duality matrices which
     are given by loops of simple matrices. So we have to find all
     independent conditions. Every loop of simple matrices
     corresponds to the statement that the product of the corresponding
     duality matrices is equal to one. Fixing a collection of
      {\it fundamental loops}, like triangles, squares , pentagons
      and hexagons before, and filling faces of such loops
     gives us the graphs complex $C(TG)$. The relations for the
     loops are complete iff the complex is simply connected.
   Thus we get a system of (polynomial) equations for duality
    matrices. Solutions to this system of equations give the {\it
    moduli space of rational conformal field theories}. Sending
 the   level $k$ to infinity we get a fusion algebra which for our
     $\SU(2)$-case is the representation algebra $R(sl(2,\C))$
      from subsection 5.1.
\end{rmk}

Returning to wave function spaces $\sH^k_\Ga$ (4.43)
 we will see that the system
of these spaces defines a vector bundle
\begin{equation}
 \sH^k_{TG} \to C(TG_g) = S^2
\end{equation}
  on the $TG$- complex  with a special structure-some local
linear decomposition on line subbundles.

Remark that the bundles is a subbundle of the trivial bundle
(see (4.45))
\begin{equation}
 \sH^k_{TG} \subset \C^{(2k)^{3g-3}} \times C_{TG} = (\oplus_{\al \in T^A_{2k}} \C_\al)
\times C_{TG}.
\end{equation}

 All these constructions  suggest  existence of
a flat connection $a_K$. Since the complex  $C_{TG}$ is simply
connected the bundle $\sH^k_{TG}$ has to be trivial with the
trivialization defined up to scaling. From this we get the
projective  bundle
\begin{equation}
\PP \sH^k \to \sM_g
\end{equation}
on the moduli space of curves of genus $g$ with the projective
flat connection. But this bundle is a slight projective bundle
with the phase amplitude (6.18) only.

Returning to the graph $C^{1}_{TG_g}$ we will see that
\begin{enumerate}
\item the system of duality matrices defines a flat connection
on $C_{1}(TG_g)$;
\item this connection has  trivial projective monodromy, thus
\item this connection defines a representation
\begin{equation}
\rho \colon Aut C_{1}(TG_g) \to PGL(rk_{g}(\emptyset)), \C).
\end{equation}
\item this representation is equal to the desired Kohno representation.
\end{enumerate}

\subsection{Gauge theory on graphs.}

A connection $a$ on the trivial $\SU(2)$-bundle on $\Ga$ is a map
\begin{equation}
a\colon P_1(\Ga)\to \SU(2)
\end{equation}
such that for the  orientation reversing involution $i_e$ we have
\begin{equation}
a(i_e(\vec e))=a(\vec e)^{-1}.
\end{equation}
Then the ``path integral'' is
\begin{equation}
a( (\vec e_1,\dots, \vec e_d))=a(\vec e_1) \cdot \cdots \cdot
a(\vec e_d) \in \SU(2)
\end{equation}
that is, it is the map
\begin{equation}
a\colon P_d(\Ga)\to \SU(2)
\end{equation}
such that for the orientation reversing involution $i_\De$ we have
\begin{equation}
a(i((\vec e_1,\dots, \vec e_d)))=a((\vec e_1,\dots, \vec
e_d))^{-1}.
\end{equation}
In the same vein we have the monodromy map for loops
\begin{equation}
a\colon L_d(\Ga)_v \to \SU(2).
\end{equation}
Obviously every such connection is {\it flat }.

Let $\sA(\Ga)$ be the space of connections that is the space of
functions (6.56) subjecting to the constraint (6.57):
\begin{equation}
\sA(\Ga)=\{ a\in \SU(2)^{P_1(\Ga)} \vert a(i_e(\vec e))=a(\vec
e)^{-1} \}.
\end{equation}

Every element $\w{g}$ of the gauge transformations group
$\sG(\Ga)$ is a function
\begin{equation}
\w{g}\colon V(\Ga)\to \SU(2)
\end{equation}
that is
\begin{equation}
\sG(\Ga)=\SU(2)^{V(\Ga)}
\end{equation}
with the componentwise multiplication.

This group acts on the space of connections $\sA(\Ga)$ by the
following rule
\begin{equation}
\w{g}(a(\vec e))=\w{g}(v_s(\vec e)) \cdot a(\vec e) \cdot
\w{g}(v_t(i(\vec e))).
\end{equation}
Recall that $\w{g}(v_t(i(\vec e)))= \w{g}(v_t(\vec e))^{-1}.$

The space of gauge orbits
\begin{equation}
\sB(\Ga)=\sA(\Ga) /\sG(\Ga)
\end{equation}
is the space of classes of representations of the fundamental
group of $\Ga$ as it is expected because all our connections are
flat:
\begin{equation}
 \sA(\Ga) /\sG(\Ga)=CLRep(\pi_1(\Ga), \SU(2))
\end{equation}
 is the quotient of the space of representations  by the adjoint action
\begin{equation}
 CLRep(\pi_1(\Ga), \SU(2))=\Hom(\pi_1(\Ga) , \SU(2)) /Ad \SU(2)
\end{equation}
- the space of classes of representations of this group.

Obviously the fundamental group of 3-valent graph $\Ga$ of genus
$g$ is the free group with $g$ generators. Thus  the quotient
(6.66)
\begin{equation}
CLRep(\pi_1(\Ga), \SU(2)) = uS_g
\end{equation}
is the {\it unitary Schottki space}. (see the definition (4) in 3.2
after the formula (3.44)).

The gauge transformation group $\sG(\Ga)$ contains the diagonal
subgroup
\begin{equation}
\SU(2)_\De \subset \sG(\Ga)
\end{equation}
of constant functions (6.63). The action of this subgroup defines
the space of constant orbits
\begin{equation}
 CL \sA(\Ga)=\sA(\Ga) / \SU(2)_\De.
\end{equation}

For every vertex $v\in V(\Ga)$ we have the subgroup
\begin{equation}
\sG(\Ga)_v=\{\w{g}\in \sG(\Ga) \vert \w{g}(v)=\id \}
\end{equation}
of gauge transformations preserving the framing at $v$. This is a
normal subgroup, and
\begin{equation}
\sG(\Ga) /\sG(\Ga)_v=\SU(2)_\De
\end{equation}
Thus the full group of gauge transformations is a semidirect
product of $\sG(\Ga)_v$ and $\SU(2)_\De$.

The quotient
\begin{equation}
\sA(\Ga) /\sG(\Ga)_v=Rep(\pi_1(\Ga, v), \SU(2))
\end{equation}
is the orbit space of framed (at $v$) connections. This space
depends on
 the choice of a point $v$.

Thus the quotient map (6.67)
\begin{equation}
P\colon \sA(\Ga)\to CLRep(\pi_1(\Ga), \SU(2))
\end{equation}
can be decomposed as follows
\begin{equation}
\sA(\Ga) \xrightarrow{\,P\,} Rep(\pi_1(\Ga), \SU(2))
\xrightarrow{/\SU(2)\De} CLRep(\pi_1(\Ga), \SU(2))
\end{equation}
or
\begin{equation}
\sA(\Ga)\xrightarrow{/\SU(2)_\De} CL\sA(\Ga) \xrightarrow{P_{cl}}
CLRep(\pi_1(\Ga), \SU(2)).
\end{equation}
The involution $i_\De$ acts on the space $\sA(\Ga)$ of connections
\begin{equation}
i^*_\De\colon \sA(\Ga)\to \sA(\Ga)
\end{equation}
and there exists an element
\begin{equation}
\w{g}_i=\begin{pmatrix}0&1\\-1&0\end{pmatrix} \in\SU(2)_\De
\subset \sG(\Ga)
\end{equation}
such that
\begin{equation}
i^*_\De=\w{g}_i.
\end{equation}
 Thus the involution $i^*_\De$ (3.23) acts
trivially on $CL\sA(\Ga)$. Recall that a gauge fixing is a
section of the projection $P$ (6.75).

Recall that the set $\Conj (\SU(2))$ of conjugacy classes of
elements of $\SU(2)$ can be described as the interval $[0,1]$
with respect to the map
\begin{equation}
conj\colon \SU(2)\to \Conj (\SU(2))=[0,1]
\end{equation}
sending a matrix $g\in \SU(2)$ to
\begin{equation}
 conj g=\frac{1}{\pi} \cdot \cos^{-1}\Bigl(\frac{1}{2}\Tr g \Bigr)\in
 [0,1].
 \end{equation}

Using this map coordinatewise one obtains the map
\begin{equation}
conj\colon \CL\sA(\Ga)\to [0,1]^{3g-3}=\prod_{e\in E(\Ga)}
[0,1]_e,
\end{equation}
which is obviously surjective.

The map $conj$ is the composite
\begin{equation}
 \sA(\Ga) \xrightarrow{\SU(2)_\De} \CL\sA(\Ga)\to
[0,1]^{3g-3}=\prod_{e\in E(\Ga)} [0,1]_e
\end{equation}
and the involution $i_\De$ preserves its fibers.

 \subsection{Abelian gauge theory and $U(1)$-spin networks}

 Abelian gauge theory on $\Ga$ is a good model for non-Abelian gauge
theory. Abelian spin networks can be considered as
theta-characteristics under numeration of theta functions.
 A $U(1)$-connection $a$ is a map
\begin{equation}
a\colon P_1(\Ga)\to U(1)
\end{equation}
such that
\begin{equation}
a(i_{e}(\vec e))=a(\vec e)^{-1}.
\end{equation}
The ``path integral'' is given by
\begin{equation}
a( (\vec e_1,\dots, \vec e_d))=\prod_{\vec e\in P_1(\Ga)} a(\vec
e)\in U(1),
\end{equation}
and so on. Again let
\begin{equation}
\sA_{U(1)}(\Ga)=\bigl\{ a\in U(1)^{P_1(\Ga)} \bigm| a(i_e(\vec
e))=a(\vec e)^{-1} \bigr\}.
\end{equation}
be the space of $U(1)$-connections. Then we have the same
involution
\begin{equation}
 \sA_{U(1)}(\Ga) \xrightarrow{i^*_e} \sA_{U(1)}(\Ga).
\end{equation}
Every element $\w{u}$ of the gauge transformation group
$\sG_{U(1)}(\Ga)$ is a function
\begin{equation}
\w{u}\colon V(\Ga)\to U(1).
\end{equation}
This group acts on the space of connections $\sA_{U(1)}(\Ga)$ by
the same rule:
\begin{equation}
\w{u}(a(\vec e))=\w{u}(v_s(\vec e)) \cdot a(\vec e) \cdot
\w{u}(v_t(i(\vec e))).
\end{equation}

However, there is one important difference between Abelian and
non-Abelian theories. Namely, the diagonal group
\begin{equation}
U(1)_\De \subset \sG_{U(1)}(\Ga)
\end{equation}
has  trivial action. Thus
\begin{equation}
\sB_{U(1)}(\Ga)=\sA_{U(1)}(\Ga) /\sG_{U(1)}(\Ga)=\Hom(\pi_1(\Ga) ,
U(1))=U(1)^g.
\end{equation}
Remark that this is abelian unitary Schottki space, that is, the
zero fiber of the projection $\pi$ (2.22) on which we are
constructing our delta functions (2.46) and (2.48) for the
application of the coherent state transform to get an analytical
presentation of theta functions on the complex Schottki space
$(\C^*)^g$. We would like to do this for non-abelian case too.

Let
\begin{equation}
P_a\colon \sA_{U(1)}(\Ga)\to U(1)^g
\end{equation}
be the projection map.

The abelian and non-abelian theories are related by the following
map:
\begin{equation}
 d\colon \sA_{U(1)}(\Ga)\to \sA(\Ga)
 \quad\hbox{given by}\quad
 d(a(\vec e))=\begin{pmatrix} e^{i \phi}, 0\\ 0, e^{- i \phi}
\end{pmatrix}.
\end{equation}
This map $d$ is equivariant with respect to every involution
$i^*_e$. Obviously the image $d (\sA_{U(1)}(\Ga) ) \subset
\sA(\Ga)$ is a 2-section of the projection $conj$, that is, the
composite
\begin{equation}
 conj \circ d\colon \sA_{U(1)}(\Ga)\to \prod_{e\in E(\Ga)}
[0,1]_e
\end{equation}
 is the factorization by the involution $i^*_\De$.

It is easy to check

\begin{prop} In the chain of maps
\begin{equation}
\sA_{U(1)}(\Ga) \xrightarrow{\,d\,} \sA(\Ga) \xrightarrow{\,P\,}
U(1)^g\in \CLRep (\pi_1(\Ga))
\end{equation}
the composite $d \circ P$ is the projection map $P_a$ (5.75).
\end{prop}

The set of irreducible representations of $U(1)$ is
\begin{equation}
\widehat{U(1)} = \Z
\end{equation}
and a spin $U(1)$-network on a graph $\Ga$ is given by a function
\begin{equation}
w \colon E(\Ga) \to \widehat{U(1)} = \Z
\end{equation}
(like (4.46)) such that for every triple of edges meeting at a
vertex $v \in V(\Ga)$ the following condition holds
\begin{equation}
 w(e_1) + w(e_2) + w(e_3) = 0.
\end{equation}
Thus the triangle inequality  (Clebsh-Gordan conditions)
becomes the equality.

Again a spin network is of level $k$  when $w(e) \leq 2k$.

Considering the basis $(a_1 , \dots , a_g, b_1, \dots b_g)$
(2.32) we have the decomposition of the jacobian of the curve
$\Si_\Ga$
\begin{equation}
J_{\Si_\Ga} = Hom(\pi_1(\Si_\Ga), U(1)) = \prod_{i=1}^g
U(1)_{a_i}\times \prod_{i=1}^g U(1)_{b_i}
\end{equation}
which coincides with  the decomposition (2.21),(2.32).

Now we can consider the coordinates (4.14) of the map (4.41) as
elements of
 the target group $U(1)$.  Hence the map $\pi_\Ga$ (4.41) has the
  form
 \begin{equation}
\pi_\Ga \colon J_{\Si_\Ga} \to U(1)^{E(\Ga)} = T^A
\end{equation}
where the target is the action torus (4.40). The image of this
map is a g-dimensional torus
\begin{equation}
\De_\Ga = T^g_- \subset T^A.
\end{equation}
Now we can send every spin network that is a function (6.98)
satisfying (6.99) to a point of order $k$ on $T^A$  dividing by
$k$
\begin{equation}
w(e) \to 1/k w(e).
\end{equation}
 We get the map of the set $SNW_a^k$
of abelian spin works of level $k$
\begin{equation}
1/k \colon SNW_a^k \to (T^A)_k.
\end{equation}
Amusing fact is that
\begin{equation}
1/k (SNW_a^k) \subset \De_\Ga = T^g_-.
\end{equation}
Moreover
\begin{equation}
1/k(SNW_a^k) = W^k(\Ga) = BS_k(\Ga).
\end{equation}
Thus we get the canonical identification of the set $\Z^g/ k\Z^g$ of
theta characteristics of level $k$ (2.48) with the set of spin
networks
\begin{equation}
(\Z^g / k\Z^g) = SNW_a^k.
\end{equation}
So in the classical case we can use $U(1)$-spine networks for the
enumeration of theta characteristics. In the same vein we will use
$\SU(2)$-spin networks
 for the enumeration of non-abelian theta characteristics.

Also we saw  that the harmonic analysis on $T^g_+ = U(1)^g$ is
 the standard Fourier analysis on a torus. We use the Fourier
decomposition
 to define our delta-functions (2.46) and (2.48) which CST $C_t^{-i \Om}$
 (2.42) sends to classical theta functions. Now we have to extend these
 constructions to the non-abelian case.

 \subsection{Harmonic analysis of $\SU(2)$-spin networks}

Recall that a spin network defines a map $ j \colon E(\Ga)\to
\widehat{SU(2)}$ with constraints. For a triple of representations
$j_{e_1}, j_{e_2}, j_{e_3}$  an
 intertwiner is  trivial component of the tensor product
$j_{1} \tensor j_{2} \tensor j_{3}$. Such a component exists iff
under the traditional identification $\widehat{SU(2)}=\half
\Z^{+}$ the Clebsch--Gordan conditions
 \begin{equation}
 j_{1}+j_{2}+j_{3}\in \Z \quad\text{and}\quad \vert j_{1}-j_{2} \vert \le
 j_{3} \le j_{1}+j_{2}
 \end{equation}
 (the triangle inequality) hold. Function $ j \colon E(\Ga)\to
\widehat{SU(2)}$ defines a  spin network $\Ga_j$ iff these
conditions hold for every triple $j_{v,1}, j_{v,2}, j_{v,3}$ of
representations around every vertex $v\in V(\Ga)$.

Any intertwiner $ i_v\in j_{v,1} \tensor j_{v,2} \tensor j_{v,3}$
is defined uniquely so one  omits any labeling of the vertices
and denotes  $\SU(2))$-spin network by the symbol $\Ga_{j}$.

 Spin network $\Ga_{j}$
is of  level $k$ if $j_{i} \le k$ for every edge $e_i$ and in the
last inequality we take $min(j_{1}+j_{2}, k - j_{1}-j_{2})$
instead of $j_{1}+j_{2}$. Of course if $k' > k$ then every spin
network of  level $k$ is of level $k'$ automatically. We would
like to fix a level and  denote a spin network of level $k$ by the
symbol $\Ga_{j^k}$ and the set of all spin net works by the symbol
 $SNW^k(\Ga)$. Recall that we have the identification $SNW^k(\Ga) = W^k(\Ga)
 $ (see the end of subsection 6.1).

Spin network $\Ga_{j, k}$ is a purely combinatorial object but we
can identify it with
 a function on the special space - the unitary Schottky space $uS_g$
given by the formula
 \begin{equation}
 uS_{g}=\SU(2)^{g} / Ad_{\mathrm{diag}} \SU(2)
 \end{equation}
where $Ad_{\mathrm{diag}}\SU(2)$ is the diagonal adjoint action
on the direct product. This space doesn't depend on a three
valent graph modulo the choice of appropriate topological data.
Only one invariant of a graph is used, namely its
 genus $g$.However, it is quite natural to use a graph $\Ga$
as a starting point of the construction.

Consider the product
 \begin{equation}
 \SU(2)^{E(\Ga)}=\prod_{e\in E(\Ga)} \SU(2)_e
 \end{equation}
with $\SU(2)$ components enumerated by edges of $\Ga$ which we
can identify
 with the space of flat connections $\sA$ (6.88) (under a choice of any
  orientation of $\Ga$).

Let $dx$ be the Haar measure on $\SU(2)$ normalized by the
condition $\int_{\SU(2)}dx=1$ and $\vec dx$ the product measure
on $\SU(2)^{E(\Ga)}$ normalized by $\int \vec dx=1$. Then by the
Peter--Weyl formula, any function $f\in L^{2}(\SU(2)^{E(\Ga)},
\vec dx)$ has the decomposition
 \begin{equation}
 f(x)=\sum_{\vec{\rho}\in \widehat{\SU(2)^{E(\Ga)}}}
 \Tr[B_{\vec{\rho},f} \vec{\rho}(x)],
 \end{equation}
where $\widehat{\SU(2)^{E(\Ga)}}$ is the space of irreducible
representations of $\SU(2)^{E(\Ga)}$, and $B_{\vec{\rho},f}$ are
endomorphisms of the space $V_{\vec{\rho}}$ of the representation
$\vec{\rho}$, given by
 \begin{equation}
 B_{\vec{\rho},f}=\frac{1}{\dim V_{\vec{\rho}}}\int_{\SU(2)^{E(\Ga)}}
 f(x) \vec{\rho}\,^{-1}(x) \vec dx.
 \end{equation}
Every irreducible representation of $\SU(2)^{E(\Ga)}$ is given by
tensor product of irreducible representations of $\SU(2)$ $
\vec\rho=\rho_{1}\tensor \cdots \tensor\rho_{3g-3}$. So every spin
network $\Ga_{j}$ of genus $g$ defines a representation of
$\SU(2)^{E(\Ga)}$ by the tensor product of all labels
 $\vec j=\bigotimes_{e\in E(\Ga)} j_e$.

Of course this is an analog of the standard Fourier decomposition
(2.41) or even (2.46). Here the representation $\vec \rho$ is a
label of {\it frequency} and an endomorphism $B_{\vec \rho, f}$
is the {\it Fourier coefficient } that is the analog of a {\it
number}. The formula (6.113) is nothing else than the integral
formula for a Fourier coefficient.

Every endomorphism of the space $V_{\vec j}$ is a vector in the
tensor product
$$
 \Bigl(\bigotimes_{e\in E(\Ga)} j_e\Bigr)
 \tensor \Bigl(\bigotimes_{e\in E(\Ga)} j_e\Bigr)^*=
 \Bigl(\bigotimes_{e\in E(\Ga)} j_e\Bigr)
 \tensor \Bigl(\bigotimes_{e\in E(\Ga)} j_e\Bigr),
 $$
since we are dealing with $\SU(2)$-representations.

But after symmetrizing components of the final product can be
labeled by elements of the set $F(\Ga)$, and we can decompose it
as
 \begin{equation}
 \Bigl(\bigotimes_{e\in E(\Ga)} j_e\Bigr)
 \tensor \Bigl(\bigotimes_{e\in E(\Ga)} j_e\Bigr)^*=
 \bigotimes_{v\in V(\Ga)}(j_{v,1} \tensor j_{v,2} \tensor j_{v,3}).
 \end{equation}

For every triple of representations around a vertex $v\in V(\Ga)$ we
have a vector $ i_v\in j_{v,1} \tensor j_{v,2} \tensor j_{v,3}$
 and their tensor product gives us the vector
 \begin{equation}
 B(\Ga_j)=\bigotimes_{v\in V(\Ga)} i_v \in  End V_{\vec j}\,.
 \end{equation}
So one has the function
 \begin{equation}
 f_{\Ga_{j}}(x)= Tr[B(\Ga_j) \vec j(x)]\in L^{2}(\SU(2)^{E(\Ga)},
 d \vec x).
 \end{equation}

This function is invariant with respect to the action on
$\SU(2)^{E(\Ga)} = \sA $ of the gauge group on $\Ga$  $
\SU(2)^{V(\Ga)}=\prod_{v\in V(\Ga)} \SU(2)_v$ with components
enumerated by the vertices and an orientation of $\Ga$.

It is easy to see that the endomorphism $B(\Ga_j)$ (5.96) is an
intertwiner with respect to such action. Then by the standard
result of harmonic analysis on groups our function $f_{\Ga_{j}}$
is invariant with respect to
 this action. That is, $f_{\Ga_{j}}$ is a function on the homogeneous
space
 \begin{equation}
 uS_g = \sA  / \sG.
 \end{equation}

But the set $\{ \Ga_j\} = SNW^k(\Ga) = W^k(\Ga)$ enumerates the
collection
 of functions $f_{\Ga_j}$ (6.115) each of which is an analog of the
  function
 (2.48) enumerated by the set $(\Z^g / k\Z^g)$ of theta characteristics.

 \begin{prop} Graph $\Ga$ and  function $f_{\Ga_j}$ determine
  spin network
 $\Ga_j$ uniquely.
 \end{prop}
 Indeed,
 \begin{enumerate}
 \item knowing  graph we can reconstruct the map $P \colon \sA \to uS_g$;
 \item lifting  function $f_{}\Ga_j$ on $\sA = \SU(2)^{E(\Ga)}$ and using
 the harmonic decomposition of this function on $\sA$
  (6.61)
 we determine  the representation $j_1 \otimes \dots \otimes j_{3g-3}$ and
 the space $V$ of this representation;
 \item  by  the Peter-Weyl formula we determine the endomorphism $B(\Ga_j)$
  (6.114);
  \item decomposing this endomorphism on blocks (6.113)
  correspoding to
   vertices we reconstruct the function $j$ (6.4).
   \end{enumerate}
 and we are done.

 Now we are ready to construct  CST for the pair $uS_g \subset S_g$
 (see the end of 3.2) since the complex Schottki space $S_g$ is
 the complexification of $uS_g$.

\subsection{RCFT  and Kohno representation}

The usual data for RCFT include an algebra $A$ and
a discrete set $R$ of representations
 (a finite set for level $k$) $\{ \rho_i \}$
(just like in subsection 5.1) such that
\begin{equation}
\rho \in R \implies \rho^* \in R,
\end{equation}
and a distinguished element $1$ such that $1^* = 1$. Moreover for every
ordered triple $\rho_i, \rho_j, \rho_k$ we have a finite (rationality of RQFT)
dimensional space $V^{\rho_i}_{\rho_j, \rho_k}$.
\begin{equation}
dim V^1_{\rho_i, \rho_j} = 1 \quad \text{ if } \quad \rho_i =
\rho_j^*
\end{equation}
and 0 otherwise.

\begin{rmk} So  we can extract a {\it modular tensor category}
from RCFT.
\end{rmk}

Having such $RCFT$ and 3-1-valent graph $\Ga$ and a spin network
$\Ga_j$ we get from every vertex $v \in V(\Ga)$ the spaces with
edges $e_1, e_2, e_3$ through it thus we have the space of
intertwiners
 \begin{equation}
V^{j_1}_{j_2, j_3} \in j_1 \otimes j_2 \otimes j_3
\end{equation}
(like in (6.6)). Taking tensor products of them together over all
vertices
 we get the space
 \begin{equation}
\sH^k_{\Ga_j} = \otimes_{v \in V(\Ga)} V^{j_1}_{j_2, j_3}
\end{equation}
 and taking a direct product over all possible networks (labelings) we get
 the space
 \begin{equation}
\sH^k_{\Ga} = \oplus_{j \in SNW^k(\Ga)} \sH^k_{\Ga_j}.
\end{equation}
For $\SU(2)$-case we have
\begin{enumerate}
\item $ dim V^{j_1}_{j_2, j_3} = 1$;
\item $ dim \sH^k_{\Ga_j} = 1$;
\item the space (6.121) is precisely our space (4.39), (4.43)
 of wave functions under
Bohr-Sommerfeld quantization.
\end{enumerate}
It is useful to fix bases in all our spaces. (In $\SU(2)$-case
where all spaces are 1-dimensional we can do it using the
standard Wigner 3j-symbols (for recalling and quantum groups
extention see \cite{KR}). We equipe
 the wave function spaces with  bases.

 This construction hasn't
to be invariant with respect to changes in trinion decompositions
 (or changes of bases). But having the complex $TG_g$ we know
 the obstackles.

 To get such invariance we have to produce matrices
 of transformations to new bases and check that these matrices satisfy the
 same algebraic identities as the changes of decompositions themselves.
 We have to check that these matrices satisfy the same algbraic
 equations as
 our transformations itselves (see Fig.11, Definition 9 and 1)-4) around (6.47)).

 Fortunately there are 5 elementary moves only described in \cite{MS} with
  traditional symbols
 \begin{equation}
F, \quad S, \quad, T, \quad, \Om, \quad , \Theta
\end{equation}
satisfying the required algebraic relations.

Recall that elementary transformations of 3-valent graphs admit
relations 1)-4) around (6.47). So our intertwiners vector spaces must
satisfy the  system of identities described below:
\begin{enumerate}
\item  3 transformations are diagonal matrices of phases (6.18)
 (in an appropriate basis)
\begin{equation}
\Om^{\pm}, \quad \Theta, \quad T \colon  V^{\rho_i}_{\rho_j, \rho_k}
\to V^{\rho_i}_{\rho_j, \rho_k};
\end{equation}
\item "genus 1 equality "
\begin{equation}
S \colon V^{\rho_i}_{\rho_j, \rho_j^*} \to
\sum_l V^{\rho_i}_{\rho_l, \rho_l^*};
\end{equation}
\item the fusing matrix
\begin{equation}
F \colon V^{\rho_i}_{\rho_j, \rho_k} \otimes V^{\rho_k}_{\rho_t, \rho_m} \to
\sum_{l} V^{\rho_j}_{\rho_t, \rho_l} \otimes V^{\rho_l}_{\rho_j, \rho_m};
\end{equation}
\item the braiding matrix
\begin{equation}
B^\pm = F \Om^\pm F^{-1}
\end{equation}
subjecting the ordinary Reidemeister relations
\begin{equation}
B_{12} B_{23} B_{12} = B_{23} B_{12} B_{23}.
\end{equation}
\end{enumerate}
We can see that the most general topological (geometrical) setting
for applying of discussed TQFT is a 3-manifold containing a
3-valent graph (instead a loop = knot). So the most striking
difference between such objects and spin networks is following:
in spin networks the holonomy of one edge around another is
trivial that is any "braiding" relations are absent. In spin
networks it isn't necessary to distinguish over and under crossing
when a 3-valent graph is projected on a plane. Otherwise all
rules for evaluating spin networks closely parallel the skein
rules for the Jones polynomials. Here we are translating RCFT
setting to Penrose's interpretation of the quantum gravity.

Recall that all braid properties of all matrices in the game are
arising from the holonomy of the Knizhnik-Zamolodchikov equation
(see \cite{KZ}) and the polynimial equations for these matrices
expected from the Conformal Field Theory were proposed by Moore
and Seiberg in \cite{MS}.

Consider a 1-3-valent graph $\Ga$ which corresponds to a  trinion
decomposition of a pinched Riemann sphere $ \Si = \PP^{1} - p_{1}
- \dots - p_{n} $.

Let
\begin{equation}
\{I_{\mu}\}, \quad \mu = 1,2,3
\end{equation}
be any orthonormal basis of $sl(2,\C)$ with respect to the
Cartan-Killing form and
\begin{equation}
\Om = \sum_{\mu} I_{\mu} \otimes I_{\mu} \in End (1) =
End(sl(2,\C)).
\end{equation}
Let $\pi_{i}(I_{\mu})$ be the action on i-th component of the
tensor product $j_{1} \otimes \dots \otimes j_{n}$ (of some
labeling of the graph) and
\begin{equation}
\Om_{ij} = \sum_{\mu} \pi_{i}(I_{\mu}) \pi_{j}(I_{\mu}) \in
End(j_{1} \otimes \dots \otimes j_{n}).
\end{equation}
Consider the configuration space
\begin{equation}
\C^{n} - \De = \{(z_{1}, \dots, z_{n}) \vert i \neq j \implies
z_{i} \neq z_{j}\}
\end{equation}
and the trivial vector bundle on it with the fiber $ j_{1} \otimes
\dots \otimes j_{n}$. It is easy to check that the KZ-equation
\begin{equation}
\p \Phi / \p z_{i} = \frac{1}{k+2} \sum_{i \neq j}
\frac{\Om_{ij}}{z_{i}-z_{j}} \Phi, \quad 1 \leq i \leq n
\end{equation}
defines a slightly projective flat connection
\begin{equation}
\om = \frac{1}{k+2} \sum_{i < j} \Om_{ij} d log(z_{i}-z_{j})
\end{equation}
on this trivial bundle. Any solution of the KZ-equation gives a
{\it covariant constant} (horizontal) section of this bundle.

Using the detailed information on the monodromy of this equation
Tsuchiya and Kanie proposed the operator formalism  two
dimensional conformal field theory in \cite{TK}. We can use this
information too but we would like to present the matrices of
MS-axiomatic description in more background independent way.

First of all consider the sphere with 4 holes $\Si_{0,4}$ and
 two graphs $\Ga$ and $\Ga'_{e}$ corresponding to two types of
  trinion decompositions of it (see Fig.11).

Symbolically we can note every labeling of $\Ga$ by the matrix
 \begin{equation}
\begin{pmatrix}
j_{2}&&&&j_{3}\\
&&j&&\\
j_{1}&&&&j_{4}
\end{pmatrix}
 \end{equation}
Then for the graph $\Ga'_{e}$ the matrix has the form
 \begin{equation}
\begin{pmatrix}
j_{3}&&&&j_{4}\\
&&j'&&\\
j_{2}&&&&j_{1}
\end{pmatrix}
 \end{equation}
Remark that the sets of parabolic edges are the same
\begin{equation}
P(\Ga) = P(\Ga'_{e}) = e_{1}, e_{2}, e_{3}, e_{4}
\end{equation}
but the subsets $E_{v}$ are different. Thus we have two spaces
\begin{equation}
\sH^{k}_{\Ga, j_{1}, j_{2}, j_{3}, j_{4}} = Hom_{\SU(2)}(j_{1}
\otimes j_{2}, j_{3} \otimes j_{4})
\end{equation}
(6.3) with the basis given by the set of weights $ W^{k} (\Ga,
j_{1}, j_{2}, j_{3}, j_{4}) $ (6.1) and
\begin{equation}
\sH^{k}_{\Ga'_{e}, j_{1}, j_{2}, j_{3}, j_{4}} =
Hom_{\SU(2)}(j_{2} \otimes j_{3}, j_{4} \otimes j_{1})
\end{equation}
with the basis $ W^{k} (\Ga'_{e}, j_{1}, j_{2}, j_{3}, j_{4}) $.

Obviously these spaces of (conformal) blocks coincide but the
bases are different.

\begin{dfn}
The matrix
\begin{equation}
F_{e} \colon \sH^{k}_{\Ga, j_{1}, j_{2}, j_{3}, j_{4}} =
Hom_{\SU(2)}(j_{1} \otimes j_{2}, j_{3} \otimes j_{4}) \to
\sH^{k}_{\Ga'_{e}, j_{1}, j_{2}, j_{3}, j_{4}} =
Hom_{\SU(2)}(j_{2} \otimes j_{3}, j_{4} \otimes j_{1})
\end{equation}
sending the first basis to the second is called the fusion matrix.
\end{dfn}

Let the symbol of the first basis be
\begin{equation}
 \begin{pmatrix}
j_{2}&j_{3}\\
 j_{1}&j_{4}
\end{pmatrix}
\end{equation}
and of the second be
\begin{equation}
 \begin{pmatrix}
j_{3}&j_{4}\\
 j_{2}&j_{1}
\end{pmatrix}
\end{equation}
Then the fusion matrix
\begin{equation}
F \begin{pmatrix}
j_{2}&j_{3}\\
 j_{1}&j_{4}
\end{pmatrix} = \begin{pmatrix}
j_{3}&j_{4}\\
 j_{2}&j_{1}
\end{pmatrix}
\end{equation}
can be considered as a transformation of solutions of KZ-equation
with coefficients
\begin{equation}
F_{ij} \begin{pmatrix}
j_{2}&j_{3}\\
 j_{1}&j_{4}
\end{pmatrix}
\end{equation}
In these notations it is easy to check the following equalities
\begin{enumerate}
\item \begin{equation}
\sum_{j} F_{ij}\begin{pmatrix}
j_{2}&j_{3}\\
 j_{1}&j_{4}
\end{pmatrix} F_{jk}\begin{pmatrix}
j_{3}&j_{4}\\
 j_{2}&j_{1}
\end{pmatrix} = \delta_{ik};
\end{equation}
\item \begin{equation}
F_{ij}\begin{pmatrix}
j_{2}&j_{3}\\
 j_{1}&j_{4}
\end{pmatrix} = F_{ij}\begin{pmatrix}
j_{4}&j_{1}\\
 j_{3}&j_{2}
\end{pmatrix}
\end{equation}
\item for $\Si_{0,5}$ five fusion matrices are subjecting the
pentagone equation (see 3) below (6.47)).
\end{enumerate}

Now using fusion matrices and  simply connectedness   of $C(TG)$
from 5.2 we can identify uniquely all spaces (4.39), (4.43) and
(6.121).

Using the same symbols we can define the collection of {\it
braiding matrices}. Again we can describe it as monodromies of
covariant constant solutions to the KZ-equation. Symbolically we
have the chain of transformations
\begin{equation}
\begin{pmatrix}
j_{2}&j_{3}\\
 j_{1}&j_{4}
\end{pmatrix} \to \begin{pmatrix}
j_{1}&j_{3}\\
 j_{2}&j_{4} \end{pmatrix} \xrightarrow{ F } \begin{pmatrix}
j_{3}&j_{4}\\
 j_{1}&j_{2}
\end{pmatrix} \to \begin{pmatrix}
j_{3}&j_{2}\\
 j_{1}&j_{4}
\end{pmatrix}
\end{equation}
\begin{dfn} The composition of these transformations
 is called a braid matrix.
\end{dfn}
Again we denote coefficients of this matrix by the symbol
\begin{equation}
B_{ij}\begin{pmatrix}
j_{2}&j_{3}\\
 j_{1}&j_{4}
\end{pmatrix}
\end{equation}
This braiding matrix can be diagonalized by the fusion matrix:
\begin{equation}
B \begin{pmatrix}
j_{2}&j_{3}\\
 j_{1}&j_{4}
\end{pmatrix} = F \begin{pmatrix}
j_{2}&j_{3}\\
 j_{1}&j_{4}
\end{pmatrix}^{-1} D F \begin{pmatrix}
j_{2}&j_{3}\\
 j_{1}&j_{4}
\end{pmatrix}
\end{equation}
where $D = \{d_{i}\}$ is the diagonal matrix with elements
\begin{equation}
d_{i} = (-1)^{j_{2} + j_{3} - i} e^{\pi i (\frac{i(i+1)}{k+2} -
\frac{j_{2}(j_{2}+1)}{k+2} - \frac{j_{3}(J_{3}+1)}{k+2})}.
\end{equation}
Moreover
\begin{equation}
B_{ij} \begin{pmatrix}
j_{2}&j_{3}\\
 j_{1}&j_{4}
\end{pmatrix} = (-1)^{j_{1} + j_{4} - i - j} e^{\pi i (-1)
(\frac{i(i+1)}{k+2} + - \frac{j(j+1)}{k+2} -
\frac{j_{1}(j_{1}+1)}{k+2} - \frac{j_{4}(J_{4}+1)}{k+2})} F_{ij}
\begin{pmatrix}
j_{1}&j_{3}\\
 j_{2}&j_{4}
\end{pmatrix}
\end{equation}
All of these equalities can be checked using monodromies of
KZ-equation for $n = 4$ (or for $n = 3$ plus $\infty$ as in
\cite{K1}).

Using braid matrix $B$ we can defined matrices $\Om^{\pm}$ (6.123)
by the formula (6.127).

For every edge $e \in E(\Ga)$ we have the diagonal transformation
of the space $\sH^{k}_{\Ga}$ (4.41) given in the basis $\{w\}$
(4.39) by the formula
\begin{equation}
T_{e}(w) = e^{2\pi i (\frac{w(e)(w(e) + 1)}{k+2} -
\frac{k}{8(k+2)})} w.
\end{equation}

Now consider the first "elliptic" elementary transformation in
Fig.11. This transformation preserves the shape of the graph
$\Ga_{1,1}$ but changes "a" and "b" standard generators of the
fundamental group
 of 1-torus. This graph has one loop edges $a$, one 3-valent vertex
$v$ and one parabolic edge $p$ and corresponds to 1-holed torus.
For 2-holed torus  the graph $\Ga_{1,2}$ has two edges $a$ and
$a'$, two vertices $v_{1}$ and $v_{2}$ and two corresponding
parabolic edges $p_{1}$ and $p_{2}$.

So we have the spaces (6.120)
\begin{equation}
\sH^{k}_{\Ga_{1,1}, j} \quad \text{ and } \sH^{k}_{\Ga_{1,2},
j_{1}, j_{2}}.
\end{equation}
(of course the number $j$ has to be integer in this case).

Using fusion matrices around edges "a" we construct the
isomorphism
\begin{equation}
f \colon \sH^{k}_{\Ga_{1,2}, j_{1}, j_{2}} \to \oplus_{j}
\sH^{k}_{\Ga_{1,1}, j}
\end{equation}
where $j$  satisfies the Clebsh-Gordan conditions with $j_{1}+
j_{2} + j \leq k$. The source space admits the basis $w_{kl}$
which corresponds  to the admissible weights $w(a) = k$, $w(a') =
l$. In this basis using braid matrix we define the linear operator
\begin{equation}
\sB (w_{kl}) = \sum_{n} B_{nk} \begin{pmatrix}
j_{1}&j_{2}\\
 l & l
\end{pmatrix} w_{nl}
\end{equation}
on the space $\sH^{k}_{\Ga_{1,2}, j_{1}, j_{2}}$.

Now we are ready to define {\it switching operators} representing
the action of the modular group $Mod_{1,1}$ on the space
$\sH^{k}_{\Ga_{1,1}, j}$. More precisely, it is the diffeomerphism
of the elliptic elementary transformation in Fig. 11.

\begin{dfn}
An endomorphism
\begin{equation}
S_{k}(j) \in End \sH^{k}_{\Ga_{1,1}, j}
\end{equation}
is called a  switching operator if the following 3 conditions are
satisfied
\begin{enumerate}
\item \begin{equation}
S_{k}(j)^{2}= (-1)^{j} e^{\pi i \frac{j(j+1)}{k+2}} \circ id;
\end{equation}
\item
\begin{equation}
(S_{k}(j) T_{e=a})^{3} = S_{k}(j)^{2};
\end{equation}
\item let \begin{equation}
S_{k}(j_{1}, j_{2}) = f^{-1} \circ (\oplus_{j}S_{k}(j)) \circ f
\end{equation}
be the endomorphism of $\sH^{k}_{\Ga_{1,2}, j_{1}, j_{2}}$. Then
\begin{equation}
S_{k}(j_{1, j_{2}}) \circ T_{e=a}^{-1} \circ T_{e=a'}\circ
S_{k}(j_{1}, j_{2})^{-1} = \sB.
\end{equation}
\end{enumerate}
\end{dfn}
(Geometrically the operator $S_{k}(j_{1}, j_{2})$ has to be an
intertwiner of diffeomerphsims corresponding to sliding the hole
$p_{2}$ along "a" and "b" paths on 2-torus).

The existence of solutions for small $k = 1, 2, 3$ can be varified
 by elementary computations and for general $k$ can be obtained
 from $\SU(2)$-WZW model over 2-holed torus (see for example
 \cite{TUY}).

 \begin{rmk}
 We can see that for a closed graph (without parabolic edges) we
 are using braiding matrices in the equality (6.160) only. Moreover
 in this equality we are using the combination $\sB$ of braiding
 matrices. The geometrical meaning of this observation is the
 following: the construction of switching operators almost
 independent on over and under crossings when 3-valent graph is
 projected on a plane. So we can hope that in high genus closed
 case $\sM_{g}$ evaluation for spin networks and ones for RCFT
 coincide.
 \end{rmk}

Now following Kohno \cite{K1} and using the fusion matrices $F$
(6.143), switching operators $S$ (6.156) and diagonal operators
$T_{e}$ (6.152) we construct the slightly projective
representation (6.54).

\begin{figure}[tbn]
\centerline{\epsfxsize=3in\epsfbox{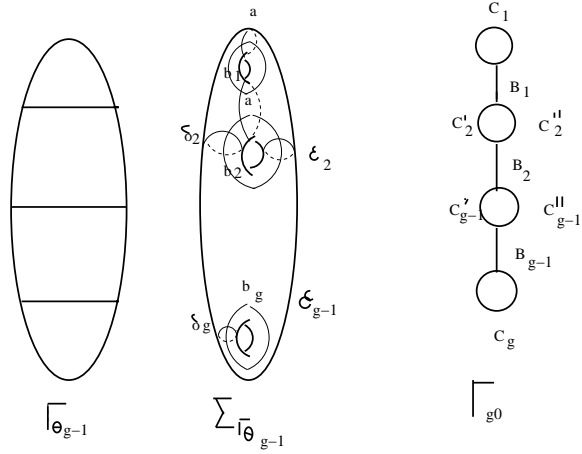}} \caption{\sl Graphs $\Theta_g$ and
 $\Ga_{g0}$} \label{Fig 14}
\end{figure}

Consider  multi theta graph $\Theta_{g-1}$ of genus $g$ as in
\cite{T3} (see Figure 18). On the Riemann surface
$\Si_{\Theta_{g-1}}$ there exists the collection of simple circles
\begin{equation}
a_{1}, \dots, a_{g}, b_{1}, \dots, b_{g}, \de_{2}, \dots,
\de_{g}, \ep_{2}, \dots, \ep_{g-1}
\end{equation}
as in Fig.12b in \cite{K1}. The collection of circles
\begin{equation}
a_{1}, \dots, a_{g}, \de_{2}, \dots, \de_{g}, \ep_{2}, \dots,
\ep_{g-1}
\end{equation}
is the collection (4.7) of disjoint inequivalent circles. These
curves correspond to edges of the graph $\Theta_{g-1}$. The
identification of these circles with corresponding edges give us
the collection of diagonal operators
\begin{equation}
T_{a_{1}}, \dots, T_{a_{g}}, T_{\de_{2}}, \dots, T_{\de_{g}},
T_{\ep_{2}}, \dots, T_{\ep_{g-1}}
\end{equation}
(6.152) on the space $\sH^{k}_{\Theta_{g-1}}$ (4.39), (4.43).

These operators represent the Dehn twists around the circles
(6.162).
 But it isn't enough to generate full $Mod_{g}$. We have to add
the Dehn twists around $\{ b_{i} \}$. To represent these twists
 we have to change the graph.

Consider the graph $\Ga_{g0}$ (see Fig. 18) which is a chain of circles $C_{1},
\dots, C_{g}$ joint by chain of consequent edges $ B_{1}, \dots,
B_{g-1}$ such that $B_{i}$ joints $C_{i}$ and $C_{i+1}$. Then
circles $C_{2}, \dots, C_{g-1}$ contain two edges
\begin{equation}
C_{i} = C'_{i} \cup C''_{i}.
\end{equation}
such that the set of edges
\begin{equation}
E(\Ga_{g0}) = C_{1}, C_{g}, C'_{2}, \dots, C'_{g-1}, C''_{2},
\dots, C''_{g-1}, B_{1}, \dots, B_{g-1}.
\end{equation}

 On
the Riemann surface $\Si_{\Ga_{g0}}$ there are circles which
correspond to edges $\{ B_{i} \}$. We cut this surface along
these circles to obtain two 1-holed tori $\Si_{1,1}^{1}$ and
$\Si_{1,1}^{g}$ and $g-2$ 2-holed tori $\{\Si_{1,2}^{i}\}, \quad
i= 2,\dots, g-1$. The 1-holed torus $\Si_{1,1}^{1}$ conains the
curve $b_{1}$ and $\Si_{1,1}^{g}$ contains $b_{g}$ from the
collection (6.44). On the
 2-holed torus $\Si_{1,2}^{i}$ we have the curves $ b_{i}, \quad
 i = 2,\dots, g-1$. The operator corresponding to the Dehn twist
  about $b_{i}$ acts on the correspoding block of the 1-holed and
   2-holed tori (6.153). So we have the collection of representatives
    of
   these Dehn twists
\begin{equation}
T_{A_{1}} \circ S_{k}(j(B_{1})) \circ T_{A_{1}}, T_{A''_{2}} \circ
S_{k}(j(B_{1}), j(B_{2})) \circ T_{A''_{2}}, \dots,
\end{equation}
$$
T_{A''_{g-1}} \circ S_{k}(j(B_{g-2}), J(B_{g-1})) \circ
T_{A''_{g-1}}, T_{A_{g}} \circ S_{k}(j(B_{g-1})) \circ T_{A_{g}}.
$$
(We identify the space $\sH^{k}_{\Theta_{g-1}}$ to
$\sH^{k}_{\Ga_{g0}}$ with the basis (4.39) enumerated by spin
networks $\{\Ga_{g0}\}$ of level $k$).

Using the description of generators of $Mod_{g}$  and
relations in \cite{Wa} we can check that all these operators give
the slightly projective representation (5.38) (see \cite{K1}).

So now we have two slightly projective representations (5.37) and
(5.38) corresponding to two versions of CFT WZW and MS-axiomatic
vertion. Fibers of the corresponding flat vector bundles are wave
function spaces for two different kind of polarizations: complex
 and real. In paper \cite{K1} Kohno considered his version of
 CFT as the combinatorial description of the space of conformal
 blocks that is the space of holomorphic sections (3.127).
But the direct identification can be constructed at the end of
this text only.

\section{Analytical aspect of the theory of
 non abelian theta functions}

\subsection{The unitary Schottky space}

This section is absolutely parallel to the subsection 2.4. The
foundation for the following constructions is  preprint
\cite{FMNT}.

The basis of the space $\sH^{k}_{\Ga}$ (6.3) is enumerated by
 the set of spin networks corresponding to the set of functions
 $\{ f_{\Ga_{j}}\}$ on the unitary Schottki space $uS_{g}$
 (6.110). On the other hand the unitary Schottky space is the fiber
 of the real polarization (4.13):
 \begin{equation}
 \pi_{\Ga}^{-1}(0, \dots , 0) = uS_{g}.
 \end{equation}

The unitary Schottky space $uS_g = \SU(2)^{g} /
 Ad_{\mathrm{diag}} \SU(2)$ is a singular manifold. To avoid this
problem we will consider geometrical objects on it as the
corresponding objects on the direct product $\SU(2)^{g}$
invariant with respect to
 $ Ad_{\mathrm{diag}} \SU(2)$-action.
 Our group $\SU(2)$ is embedded as a real part into its
  complexification
$\SL(2, \C)$. The principal fact is that the theory of finite
dimensional (non unitary) representations of $\SL(2, \C)$
completely the same answers to all questions. Obviously
\begin{equation}
 \widehat{SU(2)}= \widehat{\SL(2, \C)} =  \half \Z^{+}
 \end{equation}
and we can note both type of representations by the same symbols
$j \in \half \Z^{+}$. Remark that all spaces $V_j$ of
representations of both types of groups are precisely the same,
$rk V_j = (2j+1)$. Moreover all intertwiners in tensor
 products of representations are the same. But of course $\SL(2,
\C)$-action doesn't preseves the inner product.

So here for convenience of the reader and also to establish
notation we start by briefly recalling   from \cite{Ha} the
coherent state transform for Lie groups.

\subsection{$g$-extended Hall's construction of  CST for $\SU(2)$}

Using the heat kernel $\rho_{t}$ for a Laplacian $\Delta$ on
$\SU(2)$ associated with a given bi-invariant metric  B. C. Hall
constructed in \cite{Ha} the map
$$
C_{t}  :  L^{2}(\SU(2),dx)\mapsto {\cal H}(\SL(2, \C))
$$
\begin{equation}
C_{t}(f)(z)  =  \int_{\SU(2)} f(x) \rho_t (x^{-1} z)dx, \qquad
 f\in L^{2}(\SU(2),dx), z\in \SL(2, \C)
\end{equation}
where $dx$ is the normalized Haar measure on $\SU(2)$ and  ${\cal
H}(\SL(2, \C))$ is the space of holomorphic functions on $\SL(2,
\C)$. Recall $\rho_t$ has  unique analytic continuation to
$\SL(2, \C)$ also denoted by $\rho_t$. This map is called
$\SU(2)$-averaged coherent state transform and for each $f\in
L^{2}(\SU(2),dx)$, the function $C_{t}f(x)$ is just the solution
of the heat equation,
\begin{equation}
\frac{\partial  u}{\partial t} = \frac{1}{2}\Delta u,
\end{equation}
with initial condition given by $u(0,x)=f(x)$. Therefore,
$C_{t}(f)$ is given by
\begin{equation}
C_{t}f(z) = ({\cal C}\circ \rho_t \star f)(z) = \left({\cal
C}\circ e^{t\frac{\Delta}{2}}f \right)(z),
\end{equation}
where $\star$ denotes the convolution in $\SU(2)$ and ${\cal C}$
denotes
analytic continuation from $\SU(2)$ to $\SL(2, \C)$.

Again for every $f\in L^{2}(\SU(2),dx)$ there exists  the
expansion given by the Peter-Weyl's theorem:
\begin{equation}
f(x) = \sum_{j \in \widehat{\SU(2)}} \mbox{tr} (j (x)B_{j ,f} )
\end{equation}
where the sum is taken over the set of equivalence classes of
irreducible representations of $\SU(2)$ and $B_{j ,f}\in End \
V_j$ is given by
\begin{equation}
B_{j ,f} = \left( \dim V_j \right) \int_{\SU(2)} f(x) j
(x^{-1})dx,
\end{equation}
$V_j$ being the representation space for $j$. Then we obtain:
\begin{equation}
C_{t}f(z) = \sum_{j \in \widehat{\SU(2)}}
e^{-t\frac{\lambda_{j}}{2}}\mbox{tr} (j (z)B_{j ,f } )
\end{equation}
where $\lambda_{j}$ is the eigenvalue of $-\Delta$ on functions
of the type
$$\mbox{tr} (B_j(x) ), \quad B\in End(V_j).
$$
Remark that in this formula $j \in \widehat{\SL(2, \C)}$ but $B_j
$ is the same  endomorphism of the same space as for $j \in
\widehat{\SU(2)}$.

It turns out that there exists a natural extension of $C_{t}$ to
appropriately chosen distributions on $\SU(2)$ and we get
holomorphic functions on $\SL(2, \C)$ corresponding to
holomorphic sections of line bundles on it.

This program can be extended to $\SU(2)^g \subset \SL(2, \C)^g$
wordwisely. Moreover all constructions can be decented to
quotients
\begin{equation}
uS_g = \SU(2)^{g} /  Ad_{\mathrm{diag}} \SU(2), \quad S_g =\SL(2,
\C)^g/  Ad_{\mathrm{diag}} \SL(2, \C).
\end{equation}

A natural $g$-extension of Hall's construction of  CST can be done
immediately but an extention of the previous construction to the
complexification
\begin{equation}
uS_g = \SU(2)^{g} /  Ad_{\mathrm{diag}} \SU(2) \subset S_g
=\SL(2, \C)^g/  Ad_{\mathrm{diag}} \SL(2, \C)
\end{equation}
admits new parameters. We are starting  by precise description of
them.

 Let  $\Omega = (\Omega_{\al \be})_{\al , \be = 1}^g$ be a complex
symmetric matrix with  positive imaginary part that is a point of
Siegel upper half space $H$ and
\begin{equation}
\Omega = Re \Omega +  i\cdot Im \Omega
\end{equation}
be the decomposition on real and imagine  parts.

Consider on $uS_g$ the invariant Laplacian given by
\begin{equation}
\Delta^{(im \Omega)} = \sum_{\al, \be = 1}^g \frac{im \Omega_{\al
\be}}{2\pi i}  \frac{\partial^2 }{\partial x_\al \partial x_\be}
\end{equation}
where $x_\al,x_\be $ are invariant vector fields on $\SU(2)$.

Let $d\nu_t^{(im \Omega)}$ be the averaged heat kernel mesure on
$uS_g$ corresponding to the heat equation
\begin{equation}
\frac{\partial u}{\partial t} = \frac{1}{ 4}\Delta_\C^{(im
\Omega)} u  = \sum_{\al , \be = 1}^n \frac{im \Omega_{\al \be}}{8
\pi }\left( \frac{\partial^2}{\partial x_\al \partial x_\be } +
\frac{\partial^2 }{
\partial y_\al \partial y_\be }\right) u
\end{equation}
that is if $\mu_t^{(im \Omega)}$ denotes the fundamental solution
of this heat equation then $d\nu_t^{(im \Omega)}$ is obtained by
averaging $\mu_t^{(im \Omega)} dw$ with respect to the left
action of $uS_g$ where $dw\equiv dxdy$ is the Haar measure on
$S_g$. We get the transform
$$
C_t^{(im \Omega)} : L^2(uS_g,dx) \rightarrow  \sH (S_g)\cap
L^2((S_g, d\nu_t^{(im \Omega)})
$$
\begin{equation}
C_t^{(im \Omega)}(f)(w) = \left( {\cal C} \circ e^{\frac{t}{
2}\Delta^{(im \Omega)}}f \right) (w) \quad , \ t > 0.
\end{equation}
Now consider  non-self-adjoint Laplace operator on $uS_g$ of the
form
\begin{equation}
\Delta^{(-i\Omega )} = - \sum_{\al, \be = 1}^g \frac{i }{ 2 \pi
}\Omega_{{\al \be}} \frac{\partial^2}{ \partial x_\al \partial
x_\be}
\end{equation}
It is easy to see that the imaginary part of this Laplacian does
not affect the unitarity properties of the extended Hall's CST
defined in this way. Therefore we have the following
\begin{prop} For any $\Omega\in H_g$ and $t>0$ the transform
$$
C_t^{(-i\Omega)}  ={\cal C} \circ e^{\frac{t }{ 2}
\Delta^{-i(\Omega)}} \colon L^2(uS_g,dx) \rightarrow \sH (S_g)
\cap L^2((S_g,(d\nu_t^{(im \Omega)})
$$
is unitary.
\end{prop}
Indeed we can decompose this transform as
\begin{equation}
C_t^{(-i\Omega)}  = \left( {\cal C} \circ e^{\frac{t }{
2}\Delta^{(im \Omega)}} \right) \circ  e^{\frac{t }{ 2}
\Delta^{(-i re \Omega)}}.
\end{equation}
Then the unitarity of ${\cal C} \circ e^{\frac{t }{ 2}\Delta^{(im
\Omega)}}$ is well known and that of the operator $e^{\frac{t }{
2} \Delta^{(-i re \Omega)}} = e^{-i \frac{t }{ 2}\Delta^{(re
\Omega)}}$ (for any $t\in \R$) follows from the fact that
$\Delta^{(re \Omega)}$ is a self-adjoint operator.

Again  every $f\in L^{2}(uS_g,d\vec x)$ we can consider as $
Ad_{\mathrm{diag}} \SU(2)$-invariant  function on $\SU(2)^{g}$
and  use  the expansion given by the Peter-Weyl theorem:
\begin{equation}
f(x) = \sum_{\vec j \in \widehat{\SU(2)^g}} \mbox{tr} (\vec j
(\vec x)B_{\vec j ,f} )
\end{equation}
where the sum is taken over the set of equivalence classes of
irreducible representations of $\SU(2)^g$ and $B_{\vec j ,f}\in
End \ V_{\vec j}$ is given by the usual formula. It has to be an
intertwiner with respect to $ Ad_{\mathrm{diag}} \SU(2)$. Recall
that every class $\vec j \in \widehat{\SU(2)^g}$ is a tensor
product $\vec j = j_1 \otimes ... \otimes j_g$ where $j_i \in
\widehat{\SU(2)}$ and the same is true for $\SL(2, \C)^g$. Its
image under $C_t^{(-i\Omega)}$ is given by
\begin{equation}
C_{t}^{(-i\Omega)}f(\vec z) = \sum_{\vec j \in \widehat{\SU(2)^g}}
e^{-t\frac{\lambda_{\vec j}^{(-i\Omega)}}{2}}\mbox{tr} (\vec j
(z)B_{\vec j ,f } )
\end{equation}
where $\lambda_{\vec j}^{(-i\Omega)}$ is the eigenvalue of
$-\Delta^{(-i\Omega)}$ on functions of the type
$$\mbox{tr} (B_{\vec j} (\vec z)), \quad B\in End(V_{\vec j}).
$$
Again in this formula $j \in \widehat{\SL(2, \C)^g}$ and $B_{\vec
j} $ is the same  endomorphism of the same space as for $j \in
\widehat{\SU(2)^g}$. These endomorphisms are intertwiners with
respect to $ Ad_{\mathrm{diag}} \SL(2, \C)$-action thus
$C_{t}^{(-i\Omega)}f(\vec z)$ is a function on $S_g$.

This function is the analytic continuation to $\SL(2, \C)^g$ of
the solution of the complex heat equation on $\SU(2)^g$
\begin{equation}
\frac{\partial u }{ \partial t} = \frac{1 }{ 2}
\Delta^{(-i\Omega)} u
\end{equation}
with initial condition given by $f$.

Now this coherent states transform can be extended from
$L^{2}(uS_g,d\vec x)$ to the space of distributions
$L^{2}(uS_g,d\vec x)'$ (recall that the unitary Schottky space
$uS_g$ is compact). More precisely we have
\begin{prop} For any $\Omega\in H_g$ and $t>0$ the transform
\begin{equation}
 C_t^{(-i\Omega)} \colon (C^\infty(uS_g)' \to \sH ((\SL(2,
\C))^g)
\end{equation}
$$
f = \sum_{\vec j \in \widehat{\SU(2)^g}} \mbox{tr} (\vec j (\vec
x)B_{\vec j ,f} ) \mapsto \sum_{\vec j \in \widehat{\SU(2)^g}}
e^{-t\frac{\lambda_{\vec j}^{(-i\Omega)}}{2}}\mbox{tr} (\vec j
(z)B_{\vec j ,f} )
$$
is  well defined linear map.
\end{prop}

Indeed, the action of the Laplace operator and of its powers on
the space of distributions is defined by duality from the
corresponding action on $C^\infty(uS_g)$. For $f \in
(C^\infty(uS_g)'$ with $f = \sum_{\vec j \in \widehat{\SU(2)^g}}
\mbox{tr} (\vec j (\vec x)B_{\vec j ,f} )$ there exists positive
a constant $C_1$  such that for all $\vec j \in
\widehat{\SU(2)^g}$
\begin{equation}
\Vert  B_{\vec j ,f}   \Vert < C_1 \cdot e^{\Vert \vec j   \Vert}
\end{equation}
where
$$
\Vert \vec (j_1, ... , j_g)   \Vert = \sqrt{\sum_{l=1}^g (2j_l +
1)^2}.
$$
To show that the image of any distribution $f$ is an holomorphic
function it is enouph to observe that there exists a positive
constant $C_2$  such that for all $\vec j \in \widehat{\SU(2)^g}$
\begin{equation}
 |e^{-t\frac{\lambda_{\vec j}^{(-i\Omega)}}{2}}|\cdot \Vert  B_{\vec j ,f}
\Vert < e^{C_2 \cdot (\Vert \vec j   \Vert)^2}.
\end{equation}
for $\Vert \vec j   \Vert$ sufficiently large. The existence
follows from positivity of  imaginary part of $\Om$.

Our non-abelian  theta functions on $S_g$ are images under the CST
transform defined in previous Proposition  of linear combinations
of delta functions type distributions.

Again  let us consider the distribution
\begin{equation}
\theta^{uS_{g}}_{k} (x) = \sum_{\vec j \in \widehat{\SU(2)^g}}
\mbox{tr} (\vec j (\vec x) )
\end{equation}
that is $B_{\vec j ,f} = \id$

It is nothing else as  delta-function at identity of level $k$.
Let us apply  CST (7.20) to this distribution:
\begin{equation}
C_{1/k}^{(-i \Om)} (\theta_k^{uS_{g}}) = \sum_{\vec j \in
\widehat{\SU(2)^g}} e^{-t\frac{k \lambda_{\vec
j}^{(-i\Omega)}}{2}}\mbox{tr} (\vec j (z) )
\end{equation}

For every spin network $\Ga_{j} \in SNNW^{k}_{\Ga}$ we can start
with the function
\begin{equation}
\theta^{\Ga_{j}} (x) = f_{\Ga_{j}} \sum_{\vec j \in
\widehat{\SU(2)^g}} \mbox{tr} (\vec j (\vec x) )
\end{equation}
to get the analytic function
\begin{equation}
C_{1/k}^{(-i \Om)} (\theta^{\Ga_{j}})
\end{equation}

\begin{rmk} Again we are specializing our continious
positive parameter $t$ coming from
 heat kernels by its discret contrepart $1/k$.
 \end{rmk}

This collection of holomorphic functions on the Schottki space
\begin{equation}
\{ C_{1/k}^{(-i \Om)} (\theta^{\Ga_{j}})   \}
\end{equation}
depending on a period matrix of {\it any}  abelian variety is the
combinatorial version of non abelian theta functions - theta
functions with characteristics $\Ga_{j} \in SNW^{k}_{\Ga}$.

\section{BPU-map}

\subsection{Geometry of Lagrangian cycles}

We can see that the ranks of wave spaces  of the Kahler
quantization and the Bohr-Sommerfeld quantizations are precisely
the same that is the quantization is numerically perfect.
 We can identify projectivizations of both
 wave function spaces. But, is this identification natural?
  It will be the case if we can send every
covariant constant section to a holomorphic section of the
 line bundle by some canonical way.

The natural way was used firstly that is the  coherent state
transform or the Segal-Bargmann isomorphism introduced in the
context of the quantum theory as a transform from the Hilbert
space of of square integrable functions on the configuration
space to the space of holomorphic functions on the phase space.
In the abelian case (see subsection 2.4) considering the zero
fiber $T^g_+$ of the polarization (2.22) we can interpret our
complex torus $A$ as a complexification of this real torus and
input our geometrical situation to the situation of the classical
coherent state transform.

 In the very classical  context a configurations space is
 just $\R^g$ as the real part of the phase space $\C^g$ but the
 construction was extended to the cases when
\begin{enumerate}
\item $\R^n$ is replaced by $T^n$ and $\C^n$ by any complex
torus $A$ such that $T^n$ is a special Lagrangian subtorus (the
case of subsection 2.3);
\item $\R^n$ is replaced by $SU(2)$ and $\C^n$ by $SL(2, \C)$
(Hall, \cite{Ha}) in the section 7;
\item $\R^n$ is replaced by $SU(2)^{[g]} / diag Ad SU(2)$
 (= unitar
Schottky space) and $\C^n$ by $SL(2, \C)^{[g]}/diag Ad SL(2, \C)$
(= complex Schottky space) in the subsection 7.1;
\item $\R^n$ is replaced by any Bohr-Sommerfeld Lagrangian cycle
(more precisely a Legendrian cycle in the boundary of pseudo
convex domain) and $\C^n$ by any compact Kahlerian symplectic
manifold (Borthwick-Paul-Uribe map \cite{BPU}).
\end{enumerate}

The last case is the point of applications of the beautiful part
of calculus (Fourier integral operators, Legendrian distributions
and so on) and is a very strong tool to compare results of the
Bohr-Sommerfeld and Kaehler quantizations.

Actually the BPU-map is the extention of the classic WKB-method
"acting in the opposite direction": WKB sends a wave function to
a Lagrangian cycle with semi density in a phase space and BPU
sends any Bohr-Sommerfeld cycle with a half form to a wave
function (state) up to $U(1)$-scaling.

To start comparing the complex and symplectic geometries we
consider a polarized Kahler manifold $(M_{I}, \om, L, \nabla)$
 such that the curvature $F_{\na} = 2 \pi i \om$.
We consider every Lagrangian cycle as a class of a
 differential map of a compact smooth $1/2 dim M_{I}$-manifold
\[
\phi \colon S \to M_{I} \quad \text{ such that } \quad \phi^* \om
= 0
\]
non degenrated at some point (that is $ker d \phi_{p} = 0$ ) up to
the natural action of the group $\Diff^+$ of orientation
preserving diffeomorphisms of the manifold. Fondations and
constructions of the theory of such cycles are contained in
\cite{GT}.
 The main property of any
 such cycle is that  {\em it can't be contained by any proper
 algebraic subvariety}  (in particular, in a divisor). Thus any
 holomorphic object is uniquely determined by its restriction
 ( that is by
the operation $\phi^*$ ) to a Lagrangian cycle $\phi \colon S \to
X$. So restrictions to $S$ can serve as boundary conditions for
holomorphic sections of line bundles with curvature proportional
to $\om$.

 Geometrically, if we consider Lagrangian cycles as supports of boundary
conditions for holomorphic objects, they have the minimal possible
dimension. Usual boundary conditions deal with boundaries of
complex domains of real codimension 1. Thus it is only for
Riemann surfaces that Lagrangian boundary conditions coincide
with the usual boundary conditions. In this case, in the modern
theory of integrable systems the restriction of holomorphic
objects to a small circle around a point reduces many analytical
problems to algebraic geometry of curves (see for example the
survey \cite{DKN}). It seems reasonable to expect that
restrictions to Lagrangian submanifolds give a higher dimensional
generalization of the modern version of the theory of integrable
systems.

There are no invariants of an embedding of a Lagrangian
submanifold $S$ in a symplectic manifold. There are two ways of
getting invariants:
 \begin{enumerate}
 \item consider families of Lagrangian manifolds admitting
 invariants
(in particular limit singular subcycles); or
 \item to endow the submanifolds an additional structure (such as a
 section of
some bundle or a Hermitian connection on the trivial line bundle).
 \end{enumerate}
 The restriction of the pair $(L, \nabla)$ to a Lagrangian cycle $S$
gives this type of additional structure. It defines the space of
covariant constant sections: $ H^0_\na ((L, \nabla)\rest{S}) $
which can be nontrivial. Indeed, the restriction  to any
Lagrangian submanifold $S$ gives  topologically trivial line
bundle on $S$ with flat connection. A connection of this type is
defined by its monodromy character $
 \chi\colon \pi_1(S)\to U(1),
 $
and it admits a covariant constant section iff this character is
trivial.

 We know that a Lagrangian cycle $S$ is a  level $k$
Bohr--Sommerfeld (BS$_k$) cycle if the character is trivial.

Moreover, such a section defines a  trivialization of the
restriction , which identifies $C^\infty$ sections with complex
valued functions on $S$:
 \[
 \Ga (L^{k}\rest{S})=C^{\infty}_\C(S).
 \]
Thus the restriction to $S$ defines an embedding
 \begin{equation}
 \res\colon H^0(M_{I},L^{k}) \hookrightarrow C^{\infty}_\C(S).
 \end{equation}
up to constant.
 \begin{dfn} The image
 \begin{equation}
 \res(H^0(M_{I},L^k))=\sH_S\subset C^\infty(S)
 \end{equation}
is called  {\em analog of the Hardy space} of level $k$.
 \end{dfn}

Fixing a half-form $\hF$ on $S$; we call a pair $(S,\hF)$ a {\em
half-weighted Lagrangian} cycle or a Lagrangian cycle marked with
a half-form. Now we can identify the space of functions with the
space of half-forms
 \begin{equation}
 \Ga(L^{k})\rest{S} \cdot\hF =\Ga(\De^{1/2}),
 \end{equation}
where $\De$ is the complex volume bundle on $S$. This space is
self adjoint with respect to the natural Hermitian form given by
the integral of product of half-forms.

\subsection{Legendrian distributions}

Following Borthwick, Paul and Uribe \cite{BPU}, we can construct a
distribution in some completion of $C^\infty(S)\cdot\hF$. Its
restriction to the image of $\sH_I^k$ gives a covector or a
state.  BPU method uses usual codimension 1 boundary conditions
rather than Lagrangian boundary conditions, and the original
Hardy spaces for strictly pseudoconvex domains rather than analog
of Hardy space. We refer the reader to the beautiful paper
\cite{BPU} for the details because  in the paper \cite{GT} using
Rawnsley's results we related this construction closely to the
standard constructions of Algebraic geometry. The original
construction is the following:
  \begin{enumerate}
   \item Our
Hermitian connection on $L^*$ defines a contact structure on the
unit circle bundle $P$ of $L^*$.
\item The disc bundle in $L^*$ is a strictly
 pseudoconvex domain, and there is the Szeg\"o orthogonal projector
 $\Pi\colon L^2(P)\to \sH$ to the Hardy space of boundary values of
 holomorphic functions on the disc bundle.
  \item The contact manifold $P$
 is a principal $U(1)$-bundle, and the natural $U(1)$-action on $P$
commutes with $\Pi$ and gives a decomposition $\sH=\bigoplus_k
 H^0(M_{I},L^{k})$ of the Hardy space.
   \item If we fix a metaplectic
structure on $P$, we can lift every $BS_{k}$ submanifold to a
 Legendrian submanifold $\La\subset P$ over it, marked with the lifted
half-form $\hF$.
 \item $\La$ has an associated space of Legendrian
distributions of order $m$, which is the Szeg\"o projection of
space of
 conormal distributions to $\La$ of order $m + \half\dim M_{I}$ (see
\cite{BPU}, 2.1).
 \item A half-form on $\La$ is identified with
the symbol of a Legendrian distribution of order $m$ (see
\cite{BPU}, 2.2); thus at the level of symbols, all Legendrian
distributions look like delta
 functions or their derivatives.
  \item For a Legendrian submanifold $\La$
with a half-form we fix the Legendrian distribution of order
$\frac{1}{2}$ with symbol $\hF$ which is the Szeg\"o projection
of the delta function
 $\de_\La$.
 \end{enumerate}

So we have:
   \begin{enumerate}
    \item
For every lifting $\La\subset P$ of a $BS_{k}$ submanifold $S$
marked with a half-form $\hF$ we have a vector
\begin{equation}
BPU_{k}(\La,\hF)=\Pi^k_{\hF} (\de_\La)\in H^0(M_{I},L^{k}),
\end{equation}
 where $\Pi^k_{\hF}$ is the Szeg\"o projection to the
$k$th component of the Hardy space of the distribution with
symbol $\hF$.
 \item Every such lifting is defined up to $U(1)$-action on $P$;
 thus a pair $(S,\hF)$ defines a point of the projectivization
\begin{equation} BPU_{k} (S,\hF)=\PP (\Pi^k_{\hF} (\de_\La))\in \PP
 H^0(M_{I}, L^{k}).
 \end{equation}
 \end{enumerate}
To apply this construction to our case $(M^{ss}_{\Si_{\Ga}}, \om,
\Oh(k \Theta), A_{CS})$ (4.12) we have to observe that
 \begin{enumerate}
 \item This construction holds literally in the case
 when $\phi \colon S
 \to M_{I} $ is an immersion or $S$ has the structure of
 a smooth orbifold.
 \item By (3.92) and Propositions 12 and 13 there exists
 the canonical {\em geodesic} lifting of  Bohr-Sommerfeld cycles
 to Planckian cycles.
  \end{enumerate}

The following analog of Serre's Theorems~A and~B proved in
\cite{BPU}, Section~3:

 \begin{thm} If\/ $n$ is large enough then for any halfweighted
 $BS_k$-cycle $(S,hF)$\/ $BPU_{nk}(\La,\hF)\ne0$.  \end{thm}

There are two or three canonical ways to endow any Lagrangian
submanifold $S$ with a half-form:
 \begin{enumerate}
 \item If $X$ is a K\"ahler manifold with a metaplectic structure. Then
this metaplectic structure defines a metalinear structure on $S$
(see for example Guillemin \cite{Gu2}), and the K\"ahler metric
$g$ defines a half-form $\hF_g$ on $S$. (This method is of course
the most important for our applications.)
\item The graph of a metasymplectomorphism with symplectic volume as
square of the half-form.
 \end{enumerate}

Every half weighted cycle $(S, hF)$ is {\it weighted } by the
volume form $hF^2$ and defines the number
\begin{equation}
 v(S, hF) = \int_S hF^2
\end{equation}
-the {\it volume } of this weighted cycle.

 We write
 \begin{equation} BS_{k}^v
 \end{equation}
 for the family of half-weight $BS$ cycles
of a volme $v$. Then this space is an infinite dimensional {\it
complex Kahler manifold } (see \cite{T5} and \cite{GT}).

  The BPU construction gives a
``rational'' map
 \begin{equation} \PP \fie_k\colon BS^v_{k} \to \PP
H^0(M^{ss}_{I}, L^{k})^*
 \end{equation}
  the differential of which
is dual to the restriction monomorphism (1.27)
 (see \cite{T7} and \cite{GT}).

Recall that to send a half-weighted $BS_{k}$-cycle $(S, hF)$ to a
section $BPU_{k}(S, hF)$ we have to lift it to a Planckian cycle.
If the canonical class $K_{M_{I}}$ is proportional to
$c_{1}(L^{k}) = k \cdot[\om]$, we can do this up to finite set of
possibilities which we call a {\it fixing of theta-structure of
level} $k$ (see Definition 16 below). More precisely for every level $k$ we have the tower
 of covers:
 \begin{equation}
 M_{k-theta} \to M_{SNW(k)} \to M_{\Ga_{g}} \to \sM_{g}
 \end{equation}
where $\sM_{g}$ is the moduli space of curves of genus $g$ (5.31),
$M_{\Ga_{g}}$ is the moduli space of marked Riemann surfaces (4.8),
 $M_{SNW(k)}$ is moduli space of marked curves $\Si_{\Ga}$ with
fixed a spin network (4.42) and $M_{k-theta}$ is a curve with fixed
theta-structure for the fixed half-weighted  $k$-Bohr-Sommerfeld
fiber of the projection $\pi_{\Ga}$ (4.13). Let
\begin{equation}
\phi_{k-theta} \colon M_{k-theta} \to \sM_{g}
\end{equation}
be the composition of all these covers. Then the preimage of our
vector bundle $\sH^{k} \to \sM_{g}$ (1.29) on this cover
\begin{equation}
\PP \phi^{*}_{k-theta}  \sH^{k} = \PP(\oplus_{w \in SNW_{k}(\Ga)} \Oh
\cdot \theta_{w})
\end{equation}
where
\begin{equation}
\theta_{w} = BPU_{k-theta} (\widetilde{\pi_{\Ga}^{-1}(w)},
hF_{can}).
\end{equation}
Here $\widetilde{\pi_{\Ga}^{-1}(w)}$ is the corresponding to the
theta-structure  Planckian cycle (8.35).

The coincidence of the Hitchin, WZW and Kohno connections implies

\begin{thm} The lifting of the Hitchin connection to the cover
 $\phi_{k-theta}$ (8.10) is
Hermitian.
\end{thm}
Obviously this Hermitian structure isn't invariant with respect to this cover.

To use  the BPU construction described in (8.4), (8.5) of this
 section we must repeat some of the details.  We stay in
the situation of a complex polarization of a prequantized phase
space $(M_{I}, \om, L, \nabla)$ where $(L, \nabla)$ is
$\Oh(\Theta)$.

Consider the principal $U(1)$-bundle $P$ of the dual line bundle
$L^*$, the unit circle bundle in $L^*$. Let $D\subset L^*, \quad
\partial D=P $
 be the unit disc subfibration with boundary.

 Our complex polarization $I$, that is a compatible complex
 structure $I$ on $M$, defines a complex structure $D_I$ on $D$ as
 a strictly pseudoconvex domain.(Recall that $L$ is a positive line
 bundle).

   The Hermitian
 connection on $L^*$ iduced by $\nabla$ is given by 1-form $\al$
  on $P$ which defines a
 contact structure on $P$ with volume form \[ \frac{1}{2\pi} \al \wedge
 d\al^n, \quad \text{where} \quad n=\dim_{\C} M_{I}.  \] The null space of
 $\al$ at a point $p\in P$ is the maximal complex subspace of the tangent
 space for any complex structures.

 For every our $I$ the Hardy subspace
 \begin{equation}
 \sH_I\subset L^2(P)
 \end{equation}
consists of boundary values of holomorphic functions on $D_I$. We
have the Szeg\"o orthogonal projector
 \begin{equation}
 \Pi_I\colon L^2(P)\to \sH_I.
 \end{equation}
The natural action of $U(1)= Aut P$ on $P$ as a principal bundle
commutes with $\Pi_I$ and decomposes the space $\sH_I$ as a
Hilbert direct sum of isotypes:
 \begin{equation}
\sH_I=\bigoplus_{k=0}^{\infty} \sH_I^k ; \quad \sH_I^k
=H^0(M_{I}, L^{k})),
 \end{equation}
with only positive characters. Thus
 \begin{equation}
 \Pi_I=\bigoplus_{k=0}^{\infty} \Pi_k.
 \end{equation}
 Remark that in our situation the vector bundle $\sH^k \to \sM_{g}$
  (1.29) is embedded
 to the trivial
 bundle $L^2(P) \times \sM_{g}$ with the space (8.1) as a fibre
 and there
 exists the Szego projection
 \[
 \Pi \colon  L^2(P) \times \sM \to \sH = \oplus_{k=0}^{\infty} \sH^k,
 \]
 which is (8.2) fiberwise.

 Recall that the vector bundles  $\sH^k \to \sM_g$   admit
 projective flat connections. So locally the space $\PP\sH^k_I$ doesn't
 depend on  complex structure $I$.  Then there exists a
central extension $\widetilde{Mod_{g}}$:
  \begin{equation} 1\to U(1)\to\widetilde{Mod_{g}}\to Mod_{g} \to 1,
 \end{equation}
where the centre $U(1)= Aut P$ acts on $P$ as the proper group of
a principal bundle. This action induces a natural representation
 \begin{equation}
 \rho\colon \widetilde{Mod_{g}} \to Op(L^2(P))
 \end{equation}
to the operator algebra of the space of functions.

For every $g \in \widetilde{Mod_{g}} $ the projectors $\Pi_I$
(8.4) define the linear transformation
 \begin{equation}
p_{I, g(I)} \circ  \Pi_{g(I)} \circ \rho(g)  \circ \Pi_I \in
Op(\sH_I)
\end{equation}
 which commutes
with the action of $U(1)$. So this transformation is decomposed
  as a Hilbert direct sum of isotypes
   \begin{equation}
  Op(\sH_I) = \oplus_k       End H^0(L^k).
  \end{equation}
  These transformations  can be projectivized to the representation
 \[
 \PP\rho_I^k\colon Mod_{g} \to Aut\PP H^0(L^k).
 \]
Recall that our $M_{I}$ and the principal bundle $P$ have given
metaplectic structures.

Now if $(S, hF)$ is a half-weight $BS$-orbifold, the complex
conjugate of a covariant constant section gives a lift of it to a
half weighted Legendrian cycle  $(\La,\hF)$, and the BPU --
construction (8.4) defines a section $ BPU(\La,\hF)\in
H^0(M_{I}, L) $.
 In the same vein, we get lifting $\La_{k}$ of $S$ from $BS_{k}$
 to the Legendrian cycle on $P$ which is a cyclic $(k)$-cover of
 it and
the system of sections
 \begin{equation}
 BPU_{k}(\La_{k},\hF)\in
H^0(X, L^{k}).
 \end{equation}

 The Lagrangian cycle $S$ can be reconstructed from the set of
 its BPU
images as the quasi-classic limit as $1/k=\text{Planck's
constant}\to0$: the wave fronts of distributions concentrate on
$S \subset M_{I}$ (see \cite{BPU} and the references given there).

Now for two Lagrangian cycles $(S_1,\hF_1)$ and $(S_2,\hF_2)$, the
asymptotic behaviour of the scalar product
$\Span{BPU_k(S_1,\hF_1), BPU_k(S_2,\hF_2)}$ (see (1.9)) can be
computed in terms of the intersection $S_1\cap S_2$ (see
\cite{BPU}). In particular,
 \begin{equation}
S_1\cap S_2=\emptyset \implies BPU_k(S_1,\hF_1) \perp
BPU_k(S_2,\hF_2)
 \end{equation}
asymptotically as $k\to\infty$. (For the orbifold case these
asymptotics are somewhat weaker, but are still quite expressive
for geometric corollaries). This asymptotic technique comes from
the physical interpretation of this set-up as the ``classical''
geometric quantization. More precisely the asymptotic analysis of
quantum states gives
 \begin{multline}
 \Span{BPU_k(S,\hF), BPU_k(S,\hF)} \ \sim\
 \left(\frac{k}{\pi}\right)^{\half dim M_{I}}\int_{S} \vert\hF \vert^2 \\
 +O(k^{\half(dim M_{I}-1)}),
 \end{multline}
 and if $S_1\cap S_2=\emptyset$ then
 \begin{multline}
 \Span{BPU_k(S_1,\hF_1), BPU_k(S_2,\hF_2)} \ \sim\
 \left(\frac{k}{\pi}\right)^{\half(dim M_{I}-1)}+ \\
 O(k^{\half(dim M_{I}-2)}).
 \end{multline}
 From this  and the description of the differential of BPU-map we get
  \begin{prop} Let $\pi\colon M_{I}\to B$ be any real
 polarization  with  transversal set of
BS$_{k}$ - fibers.  Then we can fix such half forms on these
fibers then for $k\gg0$, the BPU vectors in $H^0(M_{I},L^{k})$ are
linear independent:
  \[ rk
\Span{BPU_{k}(BS_{k} \cap B)}=\#(BS_k \cap B).  \]
\end{prop}

 \begin{rmk} A more sophisticated analysis of the asymptotics of quantum
states extends this observation as follows: let
 \[
S_1,\dots, S_{N_{\max}}\subset BS_{k}
 \]
be a maximal collection of {\em disjoint\/} $BS_{k}$ cycles. Then
 \[
N_{\max} \le rk H^0(L^{k}),
 \]
and the right-hand side is given by the Riemann--Roch theorem.
\end{rmk}

\subsection{Geodesic lifting}

The lifting of a BS-Lagrangian cycle in $M_{I}$ to a Legendrian
cycle on $P$ we have described is defined up to the natural
$U(1)$-action on $P$ and the states $BPU_k(S,\hF)$ are defined
up to a phase. But in our situation we can do this almost
canonically (up to a finite ambiguity called a choice of a
theta-structure) and get an actual basis of $H^0(M_{I}, L^{k})$.

To describe this almost canonical lifting we must consider the
Lagrangian Grassmannization of the tangent bundle of $M_{I}$ as
described in \cite{T3}. Pointwise, the tangent space $(TM_{I})_x$
at a point $x\to M_{I}$ is $\C^n$ with  constant symplectic form
$\Span{\ \,,\ }=\om_x$ and  constant Euclidean metric $g_x$,
giving the Hermitian triple $(\om_x,I_x,g_x)$. Define the
Lagrangian Grassmannian $(LaGr)_x =LaGr(TM_{I})_x$ to be the
Grassmannian of oriented Lagrangian subspaces in $(TM_{I})_x$.
Taking this space over every point of $M_{I}$ gives the oriented
Lagrangian Grassmannization of $TM_{I}$
 \begin{equation}
 \pi\colon LaGr(TM_{I})\to M_{I}
\quad\text{with} \quad \pi^{-1} (x) =(LaGr)_x.
 \end{equation}
A complex structure $I_x$ on $(TM_{I})_x$ gives the standard
identification
 \begin{equation}
(LaGr)_x=U(n) / SO(n).
 \end{equation}
This space admits a canonical map
 \begin{equation}
 det\colon(LaGr)_x\to U(1)=S^1_x
 \quad\text{sending $u \in U(n)$ to $det u \in   U(1)=S^1$.}
 \end{equation}
 Taking this map over every point
of $M_{I}$ gives the map
 \begin{equation}
 det\colon LaGr(TM_{I})\to S^1(L_{-K}),
 \end{equation}
where $S^1(L_{-K})$ is the unit circle bundle of the line bundle
$\bigwedge^{n} TM_{I} =det TM_{I}$, with the first Chern class
 \[
 c_1(det TM_{I})=-K_{M_{I}}
 \]
 the canonical class of $M_{I}$ (see for example \cite{T2} and
\cite{T3}).

We have already noted that our Lagrangian cycles does not usually
have an orientation defined a priori. Thus we must consider the
Lagrangian Grassmannian $\La(TM_{I})$ forgetting orientations.
Then we get a map
 \[
det\colon \La(TM_{I})\to S^1(L_{-K/2})
 \]
in place of (8.28).

Now for every oriented Lagrangian cycle $S\subset X $, we have the
Gaussian lifting of the embedding $i\colon S\to X$ to a section
 \begin{equation}
 G(i)\colon S\to La(TM_{I})\vert_{S},
 \end{equation}
 sending $x\in S$ to the subspace $TS_x\subset(TM_{I})_x$. The composite
of this Gauss map with the projection $det$ gives the map
 \begin{equation}
 {det}\circ G(i)\colon S\to S^1(L_{-K/2})\vert_{S}.
 \end{equation}
Thus every Lagrangian cycle $S$ defines a Legendrian subcycle
 \begin{equation}
 \La={det}\circ G(i) (S)\subset S^1(L_{-K/2}).
 \end{equation}
The Levi-Civita connection of the K\"ahler metric defines a
Hermitian connection $A_{LC}$ on $L_{-K/2}$.

 \begin{dfn} A Lagrangian cycle $S$ is  almost geodesic if the
Legendrian cycle $\La={det}\circ G(i)(S)$ is horizontal with
respect to the Levi-Civita connection $A_{LC}$ on $L_{-K/2}$.
\end{dfn}

We  use Proposition 13 now. The line bundle $L_{-K/2}$ is $\Oh(2
\Theta )$ and the Levi-Civita connection is induced by the
connection $\nabla$ on $L = \Oh(\Theta)$. Then we have
 \begin{prop} A Lagrangian cycle $S$ is $BS$ if and only if it is
almost geodesic.
 \end{prop}

 The Hermitian structures of our line bundles define a map
 \begin{equation}
 \mu_2\colon S^1(L^*)\to S^1(L_{-K/2})
 \end{equation}
of the principal $U(1)$-bundles of these line bundles, which
fibrewise is minus  isogeny of degree $2$. Thus {\em every}
Lagrangian cycle $S$ defines an oriented Legendrian subcycle (see
(8.31))
 \begin{equation}
 \La=\mu_2^{-1} ({det}\circ
G(i) (S))\subset S^1(L^*)=P
  \end{equation}
(not almost geodesic {\it a priori}).
Now consider the pair of isogenies
 \begin{equation}
\mu_{k}\colon S^1(L^2)\to S^1(L^{2 k})
 \end{equation}
 \[
\mu_2\colon S^1(L^{k})\to S^1(L^{2 k})
 \]
and the lifting $ l\colon S\to S^1(L^{k}) $ given by a covariant
constant section over a $BS_{k}$ cycle $S$.
 \begin{dfn} The lift $l$ is  almost geodesic if
 \[
 \mu_{k} \circ {det}\circ G(i) (S)=\mu_2 \circ l (S).
 \]
 \end{dfn}

The number of geodesic lifts is obviously $\leq 2(k+d)$.

In summary, let $\La_k$ be the space of Legendrian subcycles of
$P$ the images of whose projection to $M_{I}$ is the $k$th root of
unity cover of a $BS_{k}$ Lagrangian cycle on $M_{I}$ (that is
Planckian cycles). Then the natural projection
 $ p\colon \La_{k}\to BS_{k} $
which sends  Legendrian cycle to Lagrangian cycle is a principal
$U(1)$-bundle.

The geodesic lifting
 \begin{equation}
l\colon \widetilde{BS_{k}} \to \La_{k}
 \end{equation}
we have described is a multisection of this principal bundle and
 $ p\colon \widetilde{BS_{k}}\to BS_{k} $
is a finite cyclic cover.

Consider a real polarization $\pi\colon M_{I}\to \De$ (1.13). Then
we have a finite set of Bohr--Sommerfeld fibres
 \[
 \De \cap BS_{k}=\{S_i\}, \quad \text{for $i=1,\dots, \vert SNW_{k}$.}
 \]

 \begin{dfn} A choice of geodesic lifts
 \[
\{\tilde{S_i}\}\subset \La_{k}
 \]
 is called a choice of  theta structure of the real polarization
$\pi_{\Ga}$.
 \end{dfn}

Marking these Lagrangian cycles with the half-forms  we get a
finite set
 \[
 \{\tilde{S_i},\hF_i\}
 \]
of half-weight Legendrian cycles.

We know that
 \[
BPU_k (\{\tilde{S_i},\hF_i\})\subset \PP H^0(M_{I},L^{k})
 \]
is a linear independent system of vectors (states) if $k\gg0$ and
an accurate choice of half forms.  In particular, if
\begin{equation}
 \#(BS_{k} \cap B)=rk H^0(M_{I},L^{k}),
 \end{equation}
  we
 get a Bohr--Sommerfeld basis.

 \begin{rmk}  For other descriptions and
applications of the geodesic lift from Lagrangian to Legendrian
cycles see \cite{T1}, \cite{T2} and \cite{T3}.
 \end{rmk}

Now on every $BS$-fiber we can fix the canonical half-form. For
this return to the description of $BS$-fibers (4.36) - (4.38).
Remark that  groups (4.30) admit bi-invariant half-forms
$\hF_1$ on $U(1)$ and $\hF_2$ on $\SU(2)$. For every $w$ we can
normalize these form $\hF_1(w)$ and $\hF_2(w)$ so that the
half-form
\begin{equation}
 \hF_w=(\hF_1(w))^{t-s} \cdot (\hF_2(w))^{p+s}
 \end{equation}
is homogeneous of degree 1 on $\pi^{-1}(w)$ (see (4.37)) with
respect to scaling $\hF_i(w)\to t\cdot\hF_i(w)$. We say that such
half-form is homogeneously normalized.

It's easy to see (\cite{JW2}, 4.7) that a normalized half-form
for a nonsingular $BS_{k}$ fibre is given by a Hamiltonian vector
field with Hamiltonian in $\R^{3g-3}$ of volume 1.

Thus every $BS_{k}$ fibre is endowed with the covariant constant
half-form (8.37), and we can proceed to construct the
corresponding Legendrian distributions in $P$. Recall that
$M_{I}$ is almost homogeneous with respect to the Goldman torus
action (4.5). Thus the Schwartz kernel of coherent states does
not depend on points and, outside singular points of fibres, they
behave as in the homogeneous case (see \cite{BPU}, (10--13)). By
lifting to $P$ every $BS_{k}$ fibre $\pi^{-1}(w)$ defines
Legendrian subcycle $\La_w\subset P$ marked with the half-form
 \begin{equation}
 \bar{hF}_w=(\fie_4\circ {det}\circ G(i))^* hF_w,
 \end{equation}
and having monodromy a $k$th root of 1.

To apply the BPU construction the principal bundle $P=S^1(L^*)$
must be given a metaplectic structure. This can be done at once
using the equality $K_{M_{I}} = - 4 \Theta$.

Our long trip is finishing at the final statement:

{\it the collection of functions}
\begin{equation}
\{ \theta^{\Ga_j} \}, \quad \Ga_j \in SNW^k_\Ga
\end{equation}
(7.25) {\it coincides with the collection of sections}
\begin{equation}
\{\theta_w = BPU_{k-theta} (\widetilde{\pi^{}-1}_\Ga(w), \bar{hF}_{w}) \}
\end{equation}
(see (8.12) and (8.38)) {\it under the meromorphic map} (3.44).

Now one can prove this up to phases (6.18). May be it is productive to input
these phase ambiguities into the non-abelian theta-structure notion?

\section{The main weapon}

\subsection{Large limit curves}

To see how all the constructions could be exploited  in the main game we have to show  behavior on the moduli space "pointwise" way. For this we are using the following strategy:
\begin{enumerate}
\item we choose special complex curves which are convinient for applications
of  constructions as algebro-geometrical so analytical  constructions and
\item extend obtained structures to general curves using the projective connections.
\end{enumerate}
Recall that the analytic geometry formulation of two-dimensional conformal field theory
(see \cite{FS})  predicts the existence of an extension of all geometrical objects like vector
 bundles (1.29) to the compactification $\ov{\sM_g}$ of the moduli space of curves. So our "special" curves may  be non smooth. Moreover our curves are
"maximally" non smooth and reducible. They look like
 maximally unipotent boundary points of moduli spaces  of 
Calabi-Yau threefolds.
Now  we define such type points for moduli spaces of algebraic
 curves  and investigate the geometry of their deformations. This
 geometry reflects
 mathematical relations between (1,1)-conformal field theories
and string theories observed  by physisists in a number of
papers. Recall that maximally unipotent boundary points  of
CY-moduli spaces correspond to large radius limit point according to 
mirror symmetry. We have not time to discuss all parallel to mirror symmetry
for algebraic curves here. However we  call our curves {\it  large limit
curves} (ll-curves for short).
 
Let $\Ga$ be any 3-valent graph of genus $g$.
Topologically it is equivalent to
3-dimensional handlebody $H_{\Ga}$ with boundary
$\p H_{\Ga} = \Si_{\Ga}$ (6.9)
where $\Si_{\Ga}$ is the Riemann surface given by the pumping up
trick (see Fig. 9) : we pump every edge of $\Ga$ to a tube (a {\it channel} in
terms of CQFT) and every vertex to a trinion (=  2-sphere with 3
holes). By the construction our Riemann surface $\Si_{\Ga}$ has a
trinions (or  "pair of pants") decomposition given by removing all
tubes.

Of course we  get a surface without any complex structure.
Thus it is not an algebraic curve, but we can certainly define uniquely   a
reducible algebraic curve $P_{\Ga}$. Namely, send every vertex $v \in
V(\Ga)$ to the complex Riemann sphere $P_{v} = \C\PP^{1}$ with a
triple points $p_{e_1}, p_{e_2}, p_{e_3}$ corresponding to edges of
 the star $S(v)$. Now for every edge $e$ with vertices $\p e = v, v'$
we identify points $ p_{e}$ of components $P_{v}$ and $P_{v'}$.
As the result we obtain a reducible algebraic curve with
properties:

\begin{enumerate}
\item  arithmetical genus of the connected reducible curve
 $P_{\Ga_{g}}$ is equal to $g$;
 \item for any 3-valent graph $\Ga_{g}$ the curve $P_{\Ga_{g}}$ is
  Deligne - Mumford stable and called  the {\it large limit
  curve};
  \item the curve $P_{\Ga_{g}}$ defines a point of the Deligne-Mumford
   compactification $\ov{\sM_{g}}$ (with the same notation);
  \item   this
point is contained by the boundary divisor $D = \ov{\sM_{g}} -
\sM_{g}$;
\item moreover,
\begin{equation}
\{ P_{\Ga}  \} \in D_{0} \cap ( \cap_{i=1}^{g-1} D_{i})
\end{equation}
where  general point of divisor $D_{0}$ corresponds to an
irreducible curve with one node and general point of $D_{i}$ is a
bouquet of smooth curves of genus $i$ and $g-i$ and
\begin{equation}
D = D_{0} \cup ( \cup_{i=1}^{g-1} D_{i})
\end{equation}
\end{enumerate}

So in the Deligne-Mumford compactification $\ov{\sM_{g}}$ of the
moduli space $\sM_{g}$ of smooth curves of genus $g$ we get the finite
configuration of points
\begin{equation}
\sP \subset \ov{\sM_{g}}
\end{equation}
enumerated by the set of 3-valent
 graphs $TG_{g}$ - the set of large limit curves now.

 \begin{rmk}
 These curves make a problem of the necessary part of two-dimensional conformal
 field theory formulation as analytic geometry (see \cite{FS}).
  Recall that the correlation functions of all local quantum
  fields can be recovered from the partition function when all
  channels (tubes)  of the surface are constricted down to nodes.
  For the surface $\Si_{\Ga}$ the constriction of all 3g-3 tubes
  produces our reducuble curve $P_{\Ga}$ with uniquely defined
  complex structure.   On the other hand from the point of view of  dynamic
  triangulation theory  these curves are the
  polymer
  phases of asymptotic  triangulations.
 \end{rmk}

\subsection{Canonical classes and canonical maps}

For every curve $P_{\Ga}$ the canonical class $K_{\Ga}$ is a line
bundle: a restriction of it to every component $P_{v}, v \in
V(\Ga)$ is the sheaf of meromorphic differentials $\om$ with
simple poles at $p_{e_1}, p_{e_2}, p_{e_3}$ where $e_1 \bigcup e_2 \bigcup e_3 =
S(v)$. Thus
\begin{equation}
K_{\Ga} \vert_{P_{v}} = K_{P_{v}}(p_{e_1} + p_{e_2} + p_{e_3} ) =
\Oh_{P_{\Ga}}(1).
\end{equation}
and in the Neron-Severi lattice $NS(P_\Ga) = \prod_{v \in V(\Ga) Pic(P_v)}$
one has
\begin{equation}
c (K_{\Ga}) = (1,1,...,1) \in NS(P_\Ga).
\end{equation}
Every holomorphic section $s$ of $K_{\Ga}$ is a collection of
meromorphic differentials $\{ \om_{v}\}$ on the components $P_{v}$
with poles at $\bigcup_{e \in E()\Ga} p_{e}$ with constraints: for every $e$ such that $\p e = v, v'$
\begin{equation}
res_{p_{e}}\om_{v} +  res_{p_{v'}} \om_{v'} = 0.
\end{equation}
On the other hand we have 2g-2 linear relations: for every $v \in
V(\Ga) $ with $(e_1,e_2, e_3 )= S(v)$
\begin{equation}
res_{p_{e_1}} \om_{v} + res_{p_{e_2}} \om_{v} + res_{p_{e_3}}
\om_{v} = 0
\end{equation}
The {\it thickness} $th(\Ga)$ of the graph is the minimal number of edges that may be removed to make the graph disconnected. In \cite{Ar} Artamkin proved the following

\begin{prop} The canonical linear system $\vert K_{P_\Ga} \vert$ is
\begin{enumerate}
\item base points free iff $th(\Ga) \geq 2$;
\item very ample iff  $th(\Ga) \geq 3$.
\end{enumerate}
\end{prop}

For simplicity below we are working with 3-valent graphs of thickness 3 only.
Then we have the canonical embedding
\begin{equation}
\phi_\Ga
 \colon P_\Ga \to \PP^{g-1} 
\end{equation}
and it is easy to check

\begin{prop} In $\PP^{g-1}$  dimension of  complete linear system
\begin{equation}
\vert 2h - \phi_{K_{\Ga}}(P_{\Ga}) \vert
\end{equation}
of quadrics through the canonical curve is equal to
\begin{equation}
dim \vert 2h \vert - 3g + 3.
\end{equation}
\end{prop}
The proof is quite classical: we just repeat here the classical argument of Castelnuovo.

\begin{cor}
For  3-valent graph $\Ga_{g}$ the large limit curve $P_{\Ga}
\in \ov{\sM_{g}}$ is a smooth point of $\ov{\sM_{g}}$ as orbifold
and the fiber of the cotangent bundle at this point is given by the following 
equality
\begin{equation}
T^{*}_{P_{\Ga}} \ov{\sM_{g}} = H^{0}(P_{\Ga}, \Oh_{P_\Ga} (2K_{P_\Ga}).
\end{equation}
\end{cor}
Now consider  the double  canonical map of $P_\Ga$ given  by the
complete linear system $\vert 2K_{\Ga}\vert$
\begin{equation}
\phi_{2K_{\Ga}} \colon P_{\Ga} \to \PP^{3g-4}.
\end{equation}
Remark that the target projective space has the interpretation
\begin{equation}
\PP^{3g-4} = \vert 2K_{\Ga}\vert^{*} = \PP
T_{P_{\Ga}} \ov{\sM_{g}}
\end{equation}
is the projectivization of the tangent space to the $\ov{\sM_{g}}$
at $P_{\Ga}$.

As any double canonical curve our configuration $\phi_{2K_{\Ga}} (P_\Ga)$ of conics and
double points is determined up to a projective automorphism of
$\PP^{3g-4} = \PP T_{P_{\Ga}} \ov{\sM_{g}} $ only. But in the next
subsection we show that  {\it the double canonical
model of a large limit curve is determined absolutely canonically 
without the projective transformation ambiguity} (unlikely of
canonical model of general algebraic curve).

Indeed, the images  of nodes 
\begin{equation}
\{ \phi_{2K_{\Ga}}(\bigcup_{e \in E(\Ga)} p_{e}) \}
\end{equation}
define the configuration of 
$$3g-3 = rk H^0(P_\Ga, \Oh(2K_{P_\Ga}))$$ 
linear independent  points ( since components are coming to irreducible conics
and every conic is defined by any triple of points on it).  Thus this configuration of points gives the decomposition of the tangent space
\begin{equation}
T_{P_{\Ga}} \ov{\sM_{g}}= \bigoplus_{e \in E(\Ga)} \C_{e}
\end{equation}
where $\PP \C_{e} = \phi_{2K_{\Ga}}(p_{e})$.

In the next
subsection we show that every tangent direction $\C_{e}$ corresponds
to a rational curve in $\ov{\sM_{g}}$. This identification of images of nodes
under the double canonical embedding and the directions of special deformations cancels  the projective transformation ambiguity.

\subsection{Special 1-parameter deformation of large limit curves}

Let $ e \subset \Ga$ be  a 3-valent graph $\Ga$ with
fixed edge  $e \in E(\Ga)$ such that $v, v' = \p e$ are two
different points with stars $S(v)= e, e_{1}, e_{2}$ and $S(v') =
e, e'_{1}, e'_{2}$. That is $e$ is not a loop. We call such pair
a simple {\it flag}.  Let us blow down this edge in our graph. We
get  new graph with 4-valent vertex $v_{new} = v = v'$ with the
star $S(v_{new}) = e_{1}, e_{2}, e'_{1}, e'_{2}$. There are 3
partitions this edges to pairs: the old one $(e_{1}, e_{2}) \vert
(e'_{1}, e'_{2})$ and couple of new ones:
\begin{equation}
(e_{1}, e'_{2} \vert e_{2}, e'_{1}) \text{ and } (e_{1}, e'_{1}
\vert e_{2}, e'_{2}).
\end{equation}
Now we can blow up the vertex  $v_{new}$ to the edge $e_{new}$
with vertices $\p e_{new} = v_{new}, v'_{new}$ and stars
\begin{enumerate}
\item $S(v_{new}) = e_{new}, e_{1}, e_{2}$ and $S(v'_{new}) = e_{new},
e'_{1}, e'_{2}$. This is our starting graph $\Ga$.
\item $S(v_{new}) = e_{new}, e_{1}, e'_{2}$ and $S(v'_{new}) = e_{new},
e'_{1}, e_{2}$. This is the first new graph $\Ga'$.
\item $S(v_{new}) = e_{new}, e_{1}, e'_{1}$ and $S(v'_{new}) = e_{new},
e'_{2}, e'_{2}$. This is the second new graph $\Ga''$.
\end{enumerate}
Remark that as a result of this construction every obtained graph
 has distinguished edge that is we have the triple of flags
\begin{equation}
(e \subset \Ga), ( e_{new} \subset \Ga'), ( e_{new} \subset
\Ga'').
\end{equation}
Such triple we call a {\it nest} of flags. Every nest is defined
uniquely by any flag $(e \subset) \Ga$ of the triple. Moreover,
let $S(e)$ be the  star  of an edge $e$ that is  the union
of stars of the boundary $v, v' = \p e$
\begin{equation}
S(e) = S(v)\bigcup S(v').
\end{equation}
Then we have the canonical identification of graphs
\begin{equation}
\Ga - S(e) = \Ga' - S(e_{new}) = \Ga'' - S(e_{new}).
\end{equation}

\begin{rmk}
It is easy to see that if $e$ is a loop this construction gives
the same graph with the same loop $e$ again.
\end{rmk}
Consider two componets $P_{v}, P_{v'}$ of a curve $P_{\Ga}$ with
 common point $p_{e}$ where $v, v' = \p e$. Remove this point and
 glue $P_{v}$ and $P_{v'}$ by a tube that is, consider the
 connected sum
 \begin{equation}
 P_{v}\#_{p(e)} P_{v'} = P_{v,v'} = S^{2}.
 \end{equation}
This is a 2-sphere with two pairs of points $(p_{e_{1}}, p_{e_{2}})$
and $(p_{e'_{1}}, p_{e'_{2}})$ where $(e, e_{1}, e_{2}) = S(v)$
and $(e, e'_{1}, e'_{2}) = S(v')$. If we fix a complex structure
on $S^{2}$ and consider the double cover
\begin{equation}
\phi \colon E \to \PP^{1}
\end{equation}
with ramification points
\begin{equation}
W = p_{e_{1}} \bigcup p_{e_{2}} \bigcup p_{e'_{1}} \bigcup
p_{e'_{1}}
\end{equation}
we obtain an elliptic curve $E$ with a point of second order
\begin{equation}
\si = p_{e_{1}} + p_{e_{2}} - p_{e'_{1}} - p_{e'_{2}}.
\end{equation}
So the moduli space of complex structures on $S^{2}$ is equal to
$\sM_{1}^{2}$ - the moduli space of smooth elliptic curves with fixed
point of order 2.

Every such complex structure $\tau \in \sM_{1}^{2}$ on $S^{2}$ and the
standard complex structures on all others components define a
stable algebraic reducuble curve $P_{e \subset\Ga, \tau}$. Thus
we obtain an embedding
\begin{equation}
\psi_{e \subset \Ga} \colon \sM_{1}^{2} \to \ov{\sM_{g}}.
\end{equation}
The moduli space
\begin{equation}
\sM_{1}^{2} = \PP^{1} - (\si, \si', \si'')
\end{equation}
is the projective line without 3 points. These 3 points
correspond to 3 possibility to divide 4 points $p_{e_{1}},
p_{e_{2}}, p_{e'_{1}}, p_{e'_{2}}$ in 2 pairs that is a choice of
a point of order 2 on an elliptic curve.

It is easy to see that

\begin{prop}
\begin{enumerate}
\item The map $\psi_{e \subset \Ga}$ (9.24) can be extended to the map
\begin{equation}
\psi_{e \subset \Ga} \colon \PP^{1} \to \ov{\sM_{g}};
\end{equation}
\item 
\begin{equation}
\psi_{e \subset \Ga} (\si) = p_{e} \in  P_{\Ga}
\end{equation}
where $\si$ is given by (9.23) and $P_{\Ga}$ is a large limit curve with
fixed node corresponding to the edge $e$,
\item
\begin{equation}
\psi_{e \subset \Ga} (\si') = P_{\Ga'};
\end{equation}
and
\item
\begin{equation}
\psi_{e \subset \Ga} (\si'') = P_{\Ga''};
\end{equation}
where the  triple
\begin{equation}
(e \subset \Ga), (e_{new} \subset \Ga'), ( e_{new} \subset \Ga'')
\end{equation}
is a nest of flags (9.17).
\item Now we can identify the tangent direction
\begin{equation}
T\psi_{e \subset \Ga} (\PP^{1})_{P_{\Ga}} = \C_{e}
\end{equation}
from the decomposition (9.15) with the corresponding node of the
double canonical curve. We obtain the canonical pluricanonical  model.
\end{enumerate}
\end{prop}
\begin{rmk} If  edge $e$ is a loop with a vertex $v$ such that
the star $S(v) = e, e'$ then $P_{v}$ is a rational curve with one
double point $p_{e}$ and the smooth point $p_{e'}$. The operation
of connected summing  around double point $p_{e}$ gives a smooth
2-torus $T^{2}$ with fixed point $p_{e'}$ and the isotopy class
$a \in H_{1}(T^{2}, \Z)$ which is the class of the neck of gluing
tube. The class $a$ mod 2 gives a point of order 2 on $T^{2}$.
The space $\sM_{1}^{2}$ of complex structures on $T^{2}$  with
such additional data is $\C^{*} - 1$ that is $\PP^{1}$ without 3
points again. In this case the rational curve $\psi_{e \subset
\Ga} (\sM_{1}^{2})$ admits the compactification by the double
point corresponding to $P_{\Ga}$.
\end{rmk}

So, the Deligne-Mumford  compactification $\ov{\sM_{g}}$ contains
the configuration $\sC$ of rational curves
\begin{equation}
\{ C \} = \sC \text{ every } C = \psi_{e \subset \Ga} (\PP^{1})
\end{equation}
for some flag $e \subset \Ga$.  We can consider this configuration of rational
curves as a reducible curve
\begin{equation}
\sC = \bigcup C
\end{equation}
that is the union of all components. It is easy to see that

\begin{prop}
\begin{enumerate}
\item for every  component $C$ of $\sC$
\begin{equation}
  C \bigcap \sP  = P_{\Ga} + P_{\Ga'} + P_{\Ga''}
\end{equation}
where $(\Ga, \Ga', \Ga'')$ is a nest of graphs (9.17 );
\item for a pair $C, C'$ of componets  either the intersection $C \bigcap C'$
 is empty or it is transversal and
 \begin{equation}
 C \bigcap C' \in \sP
 \end{equation}
 \item the set $S(\Ga)$ of components through every point
  $P_{\Ga} \in \sP$ is enumerated by the set $E(\Ga)$;
\end{enumerate}
\end{prop}

We saw that if an edge $e \in E(\Ga)$ is not a loop then the map
$\psi_{e \subset \Ga}$ (9.26) is an embedding and the component $C$
(9.24) of the reducible curve $\sS \subset \ov{\sM_{g}}$ is a
smooth rational curve with fixed 3 points $ ( P_{\Ga}, P_{\Ga'},
P_{\Ga''})$ such that the corresponding 3 graphs  form a nest
(9.17). Over every of such points the tangent spaces admit
decompositions (9.15). Geometry of the other curves from the family
of curves parametrized by $C$ is very near to the geometry of
large limit curves: let $e \in E(\Ga)$, $v, v' = \p e$ and $P_{E,
\si}$ is the point of $C$ corresponding to the elliptic curve
(9.21) with a point of order 2 (9.23). Then $P_{E, \si}$ has 2g-4
 old components
 \begin{equation}
 \bigcup_{v'' \neq v, v'} P_{v''}
 \end{equation}
with triples of points and 3g-1 nodes $p_{e'}, e \neq e'$ and one
new component $P_{v, v'}$  with the quadruple of points (see
(9.16) and 1), 2) and 3) after this formula). Then

\begin{enumerate} \item
the canonical class $K_{P_{E, \si}}$ is a line bundle: a
restriction of it to every component $P_{v''}, v'' \neq v, v' $
is the sheaf of meromorphic differentials $\om$ with simple poles
at $p_{e}, p_{e'}, p_{e''}$  where $\{e,e', e''\} = S(v)$;
\item the
restriction of the canonical class to the component $P_{v, v'}$
is the sheaf of meromorphic differentials $\om$ with simple poles
at $p_{e_{1}}, p_{e_{2}}, p_{e'_{1}}, p_{e''_{2}}$ (see (9.16)).
\item
Thus 
\begin{equation}
c (K_{P_{E, \si}}) = (2,1,...,1) \in NS_{P_{E, \si}}
\end{equation}
where the first coordinate corresponds to $P_{v, v'}$.
\item Again
every holomorphic section $s$ of the canonical class  is a
collection of meromorphic differentials $\{ \om_{v''}\}$ on the
components $P_{v''}$ with poles at $p_{e}, p_{e'}, p_{e''}$ where
$\{e,e', e''\} = S(v)$ and a meromorphic differential $\om_{v,
v'}$ on the component $P_{v, v'}$ with poles at the quadruple
with the same constraints (1.29) and (1.30).
\item
The canonical map defined by the complete linear system $\vert K_{P_{E,
\si }}\vert$
\begin{equation}
\phi_{K} \colon P_{E, \si} \to \PP^{g-1}
\end{equation}
sends $P_{v''}$ to a configuration of lines and $P_{v, v'}$ to a
conic  in $\PP^{g-1}$.
\item Again the dimension of the space of quadratic differentials
on $P_{E, \si}$ is equal to 3g-3 and this curve is an orbifold
smooth point of $\ov{\sM_{g}}$.
\item
The double  canonical map of $P_{E, \si}$  given  by the complete
linear system $\vert 2K_{\Ga}\vert$ is an embedding
\begin{equation}
\phi_{2K_{P_{E, \si}}} \colon P_{E, \si} \to \PP^{3g-4} = \PP
T\ov{\sM_{g}}_{P_{E, \si}}.
\end{equation}
(see (9.11)).
\item The images  of nodes
\begin{equation}
\{ \phi_{2K_{P_{E, \si}}}(p_{e'}) \}, \quad e \neq e' \in
E(\Ga)
\end{equation}
define the decomposition of the restriction of the tangent bundle
\begin{equation}
T\ov{\sM_{g}} \vert_{C} = TC \bigoplus (\bigoplus_{e \neq e' \in
E(\Ga)} L_{e'})
\end{equation}
where the fiber of the line bundle $L_{e'}$ over a point is the
component of the decomposition (9.15).
\item Thus every line bundle $L_{e'}$ from the previous decomposition
is the tautological line bundle
\begin{equation}
L_{e'} = \Oh_{C} (-1),
\end{equation}
\item hence the previous decomposition is
\begin{equation}
T\ov{\sM_{g}} \vert_{C} = \Oh_{C}(2)\bigoplus (3g-4) \Oh_{C}(-1).
\end{equation}
\item The restriction to $C$ of the canonical class of $\ov{\sM_{g}}$
is
\begin{equation}
K_{\ov{\sM_{g}}} \vert_C = \Oh_{C}(3(g-2)).
\end{equation}
\end{enumerate}

\subsection{ Special 2-parameter deformations of large limit curves}

Let $v \subset \Ga$ be a pair of a 3-valent graph $\Ga$ and fixed
vertex $v \in V(\Ga)$ and in the star  $ S(v) = e \bigcup e_{1}
\bigcup e_{2}$ let us choose a pair $e_{1}$ and $e_{2}$ of edges.
Then we have the pair of vertices $v_{1}$ and $v_{2}$ defined by
equalities $\p e_{1} = v, v_{1}$ and $\p e_{2} = v, v_{2}$.

Thus such choices are given by 3-flag on a graph $\Ga$
\begin{equation}
 v \subset e \subset \Ga
\end{equation}
where
\begin{equation}
e_{1} \bigcup e_{2} = S(v) - e
\end{equation}
 Now let us blow down
edges $e_{1}$ and $e_{2}$. We obtain a new graph $\Ga_{new}$
with a new 5-valent vertex $v_{new} \in V(\Ga_{new})$ and fixed
additional choices on the star of $v_{new}$:
\begin{enumerate}
\item fixed edge $e \in S(v_{new}) $ through $v_{new}$ and
\item a partition of other quadruple edges from the star $S(v_{new})$
 in two pairs $S(v_{1}) - e_{1}$ and $S(v_{2} - e_{1})$.
\end{enumerate}

There are exist 15 such choices divided into 5 triples corresponding
to the first choice of $e' \in S(v_{new}) $. For every such a choice
 let us blow up from $v_{new}$ two new edges $e'_{1}$ and $e'_{2}$
 such a way that
\begin{enumerate}
\item  $v_{new}$ becomes a new 3-valent vertex $v'$
 and $S(v') = e', e'_{1}, e'_{2}$,
 \item two vertices appear $v'_{1}$ and $v'_{2}$ such that
 $\p e'_{1} = v', v'_{1}$ and $\p e'_{2} = v', v'_{2}$.
 \end{enumerate}
Thus we get a new 3-flag
\begin{equation}
v' \subset e' \subset \Ga'
\end{equation}
and full collection of choices gives a collection of 3-flags
\begin{equation}
(v \subset e \subset \Ga),   (v' \subset e' \subset \Ga'), (v''
\subset e'' \subset \Ga''), ........
\end{equation}
We call this collection of 15 3-flags  a {\it nest of 3-flags}.

A partial case of the choice is the following:
\begin{enumerate}
\item let us blow down just one edge, for example, $e_{1}$,
\item then we get 4-valent vertex $v_{new}$ with the star
$S(v') = e, e_{2}, e_{3}, e_{4}$  where $e_{3} \bigcup e_{4} =
S(v_{1}) - e_{1}$.
\item Now let us change $e$ by any edge from the pair $e_{3},
e_{4}$, for example, by $e_{3}$ and blow up the vertex $v_{new}$
to a new edge $e'$ with vertices $\p e' = v', v''$ such that
\begin{equation}
S(v') = e' \bigcup e \bigcup e_{3} \text{ and } S(v'') = e'
\bigcup e_{2} \bigcup e_{4}.
\end{equation}
Then we have a pair of 3-flags
\begin{equation}
((v \subset e \subset \Ga)) \text{ and } (v' \subset e' \subset
\Ga')
\end{equation}
from the same nest of 3-flags such that the flags
\begin{equation}
(e_{1} \subset \Ga) \text{ and } (e' \subset \Ga')
\end{equation}
belong to the same nest (9.17).
\item if in this construction we replace $e$ by $e_{4}$ we get the
nest triple
\begin{equation}
(e_{1} \subset \Ga), (e' \subset \Ga'), (e'' \subset \Ga'')
\end{equation}
(see (9.17)).
\end{enumerate}

Any flag $(v \subset e \subset \Ga)$ defines
\begin{enumerate}
\item the reducible curve $P_{\Ga}$,
\item its irreducible component $P_{v}$,
\item a pair of its points $p_{e_{1}}, p_{e_{2}} \subset P_{v}$ and
\item
two components $P_{v_{1}}$ and $P_{v_{2}}$.
\end{enumerate}
Remove nodes $p_{e_{i}}, i = 1, 2$ and glue components by tubes
that is consider the connected sum
\begin{equation}
 P_{v_{1}} \#_{p_{e_{1}}} P_{v} \#_{p_{e_{2}}} P_{v_{2}} = S^{2}.
\end{equation}
This is 2-sphere with fixed 5 points with fixed the point $p_{e}$
and two pairs of points corresponding to the following two pairs of edges
\begin{equation}
(S(v_{1}) - e_{1}), (S(v_{2}) - e_{2}).
\end{equation}
The moduli space of complex structures on $S^{2}$ with fixed 5
points subjecting such partition conditions is
\begin{equation}
S((v \subset e \subset \Ga)) = \PP^{1} \times S^{2}(S^{2}
\PP^{1}) / PGL(2, \C),
\end{equation}
so it is rational affine surface. Every complex line with 5 points can
be uniquely extended to a reducible curve adding the corresponding
complex lines with 3 points. Thus the surface $S((v \subset e
\subset \Ga))$ is embedded
\begin{equation}
\psi_{(v \subset e \subset \Ga)} \colon S((v \subset e \subset
\Ga)) \to \ov{\sM_{g}}
\end{equation}
and its compactification belongs to
\begin{equation}
\ov{\psi_{(v \subset e \subset \Ga)}S((v \subset e \subset \Ga))}
\subset \ov{\sM_{g}}.
\end{equation}
The compactification divisor is
\begin{equation}
\ov{\psi_{(v \subset e \subset \Ga)}S((v \subset e \subset \Ga))}
- \psi_{(v \subset e \subset \Ga)}S((v \subset e \subset \Ga)) =
\bigcup C
\end{equation}
where 15 components of the curve $\sC$ corresponding to the nests
(9.48) and (9.52).

Again the Deligne-Mumford  compactification $\ov{\sM_{g}}$
contains the configuration $\sS$ of rational surfaces
\begin{equation}
\{ S \} = \sS \text{ every } S = \psi_{v \subset e \subset \Ga}
(S((v \subset e \subset \Ga)))
\end{equation}
for some 3-flag $v \subset e \subset \Ga$.  We can  consider 
this configuration of
rational surfaces as a reducible surface
\begin{equation}
\sS = \bigcup S
\end{equation}
that is the union of all components. It is easy to see that

\begin{prop}
\begin{enumerate}
\item for every  component $S$ of $\sS$
\begin{equation}
  S \bigcap \sP  = \bigcup C
\end{equation}
where the curves $ C $ belongs the collection  (9.32);
\item for a pair $S, S'$ of componets the intersection $C \bigcap C'$
 is either empty or it is transversal and
 \begin{equation}
 C \bigcap C' \subset  \sC
 \end{equation}
\end{enumerate}
\end{prop}

We saw that if an edge $e \in E(\Ga)$ from 3-flag $v \subset e
\subset \Ga$ is not a loop then the map $\psi_{v \subset e \subset
\Ga}$ (3.12) is an embedding and the component $S$ (3.13) of the
reducible surface $\sS \subset \ov{\sM_{g}}$ is a  rational
surface with fixed 15 components of the reducible curve $\sC$.
Over any of such irreducuble curve the restriction of the
tangent bundle of $\ov{\sM_{g}}$ admits decompositions (9.41) and
(9.43). The structure of a general stable curve $P_{s}, s \subset S$
from the 2-dimensional
 family parametrized by $S$
 is a slightly different from  the geometry of the large limit curves:
  $P_{s}$ has 2g-5
 old components
 \begin{equation}
 \bigcup_{v'' \neq v, v_{1}, v_{2}} P_{v''}
 \end{equation}
 and one
new component $P_{v, v_{1}, v_{2}}$  with 5 points. Then

\begin{enumerate} \item
the canonical class $K_{s}$ is a line bundle: the restriction of it
to every component $P_{v''}, v'' \neq v, v_{1}, v_{2} $ is the
sheaf of meromorphic differentials $\om$ with simple poles at
$p_{e}, p_{e'}, p_{e''}$  where $e \bigcup e' \bigcup e'' = S(v)$;
\item the
restriction of the canonical class to the component $P_{v, v_{1},
v_{2} }$ is the sheaf of meromorphic differentials $\om$ with
simple poles at $p_{e}, p_{e_{1}}, p_{e_{2}}, p_{e'_{1}},
p_{e'_{2}}$.
\item
Thus in terms of subsection 1.1
\begin{equation}
c (K_{P_{s}}) = (3,1,...,1) \in NS_{P_{s}}.
\end{equation}
where the first coordinate corresponds to $P_{v, v'_{1}, v_{2}}$.
\item
The canonical map by the complete linear system $\vert
K_{P_{s}}\vert$
\begin{equation}
\phi_{K} \colon P_{s} \to \PP^{g-1}
\end{equation}
sends $P_{v''}$ to a configuration of lines and $P_{v, v_{1},
v{2}}$ to a cubic  in $\PP^{g-1}$.
\item Again the dimension of the space of quadratic differentials
on $P_{s}$ is equal to 3g-3 and this curve is an orbifold smooth
point of $\ov{\sM_{g}}$.
\item
The double  canonical map of $P_{E, \si}$  given  by the complete
linear system $\vert 2K_{P_{s}}\vert$ is an embedding
\begin{equation}
\phi_{2K_{P_{s}}} \colon P_{s} \to \PP^{3g-4} = \PP
T\ov{\sM_{g}}_{s}.
\end{equation}
\item The images  of nodes
\begin{equation}
\{ \phi_{2K_{P_{s}}}(p_{e'}) \}, \quad e \neq e_{1}, e_{2} \in
E(\Ga)
\end{equation}
define the decomposition of the restriction of the tangent bundle
\begin{equation}
T\ov{\sM_{g}} \vert_{S} = TS \bigoplus (\bigoplus_{e \neq e_{1},
e_{2} \in E(\Ga)} L_{e})
\end{equation}
where the fiber of the line bundle $L_{e}$ over a point is the
component of the decomposition (9.43).
\end{enumerate}

\subsection{Modular configuration}

We understands the moduli subspaces (9.33) and (9.60) as partial cases of the general construction.
For  any graph $\Ga$ with vertices of any valence $ > 2$ consider the corresponding reducible curve
\begin{equation}
P_{\Ga} = \bigcup_{v \in V(\Ga)} P_v
\end{equation}
where each component equals to $\PP^1$  with $\vert S(v) \vert$
fixed points of transversal intersections with other components along the incidence correspondence  of the graph. If our graph
 isn't 3-valent such curve has moduli. These moduli are defined by  the positions of  $\vert S(v) \vert$
 points on the component $P_v$ and by the combinatoric of the set of expansions $\sE(\Ga)$ (see subsection 6.3). For example, if
$\Ga$ is the result of one contraction of 3-valent graph $\Ga_{max}$ the corresponding moduli space $\sM^2_1$ (9.24) is a 3-cover of the moduli space of smooth elliptic curves. 
For graphs with one 5-valent vertex the corresponding moduli space is the rational surface $S$ from the configuration
 (9.60) and so on. 

Every such moduli space $\sM_\Ga$ is embedded into the Deligne-
Mumford compactification and it admits itself the Deligne-Mumford compactification 
\begin{equation}
\ov{\sM_\Ga} \subset \ov{\sM_g}
\end{equation}
such that
\begin{equation}
\ov{\sM_\Ga} = \sM_v \bigcup (\bigcup_{\Ga'> \Ga}) \ov{\sM_{\Ga'}} 
\end{equation}
where  graph $\Ga'$ is an expansion of $\Ga$.

Obviously the maximal dimension of such moduli space
is $2g-3$. 

So in the compactification $\ov{\sM_g}$ of the moduli space of  genus $g$
curves  we have the configuration of the moduli spaces
\begin{equation}
\sM = \bigcup_{\Ga} \sM_\Ga \subset \ov{\sM_g}
\end{equation}
such that 
\begin{enumerate}
\item
\begin{equation}
\ov{\sM_\Ga'} \bigcap \ov{\sM_{\Ga''}}
= \bigcup_{\Ga > \Ga', \Ga'' } \sM_\Ga;
\end{equation}
\item the restriction of the tangent bundle of $\ov{\sM_g}$
to every component $\sM_\Ga$ admits the decomposition like (9.68).
\end{enumerate}

\subsection{ $Pic_0(P_\Ga)$ and moduli spaces $M_{vb}^{ss}(P_\Ga)$ of ll-curves}

For a singular curve the moduli spaces of  vector bundles are not compact. These moduli spaces admit the compactifications by torsion free sheaves. These sheaves are not local free  only over
nodes but the theory of compactification is very close to
such theory for algebraic surfaces. All details of this theory
can be found in Artamkin paper \cite{Ar}.

For every torsion free sheaf $F$ on $P_\Ga$ we have the standard exact surface
\begin{equation}
0 \to F \to F^{**} \to C(F) \to 0
\end{equation}
where 
\begin{enumerate}
\item the quotient sheaf has the form
$$
Supp C(F) \subset \bigcup_{e \in E(\Ga)} p_e;
$$
\item the double dual sheaf $F^{**}$ is a vector bundle on the normalization of the curve.
\end{enumerate}
The cardinality of $Supp F$ is the measure of deviation from the local freedom of
 $F$. Let this number be
\begin{equation}
\vert Supp C(F) \vert = c_2(F).
\end{equation}
So $c_2(F) = 0$ implies $F$ is local free thus it is a vector bundle.

A vector bundle $E$ on $P_\Ga$ is called topologicaly trivial if
the restrictions $E \vert_{P_v}$ are trivial for all $v \in V(\Ga)$.
Below we describe the moduli spaces of topologicaly trivial 
\begin{enumerate}
\item line bundles $Pic_0 (P_\Ga)$,
\item rk 2 semi-stable bundles $M_{vb}^{ss}$;
\item the compactification $J(P_\Ga) = \ov{Pic_0(P_\Ga)}$ of 
$Pic_0(P_\Ga)$ by torsion free sheaves $F$ (for which obviously
$F^{**} \in Pic_0(P_\Ga)$);
\item  the compactification $M^{ss}(P_\Ga) = \ov{M^{ss}_{vb}(P_\Ga)}$ of 
$Pic_0(P_\Ga)$ by torsion free sheaves $F$ (for which obviously
$F^{**} \in M^{ss}_{vb}(P_\Ga)$).
\end{enumerate}

Since each  component $P_{v}$ is  projective line the restriction of any topologicaly trivial line bundle equals
\begin{equation}
L \vert_{P_{v}} = \Oh_{P_{v}}.
\end{equation}

To describe a line bundle on $P_\Ga$ consider a collection of any trivializations of $L$ on all components $P_v$ of $P_\Ga$ and denote the line bundle with such addition structure by $L_0$.

Now to get a line bundle $L$ on $P_{\Ga}$ we have to concord every
 pair of line bundles $\Oh_{P_{v}}, \Oh_{P_{v'}}$
at common point $p_{e}$ if $v, v' = \p e$. Under our trivializations such concordance is
given by a multiplicative constant $ a (\vec e) \in \C^{*}$ such that under the orientation reversing involution $i_e$
$$
a(i_e(\vec e)) = a(\vec e)^{-1}
$$
(compare with (6.86)). Thus $L_0$ defines a map
$$
a \colon \vec E(\Ga) \to \C^*
$$
(compare with (6.85)). 

The changing of trivialization is given by a function
$$
\tilde{\la} \colon V(\Ga) \to \C^*
$$
(compare with (6.90)) which acts on the functions $a$  by the formula
$$
\tilde{\la}(a(\vec e)) = \tilde{\la}(v_s) \cdot a(\vec e) \cdot \tilde{\la}(v_t)^{-1}
$$
where $\p \vec e = v_s \bigcup v_t$ (compare with (6.91)).

Let $\sA_{\C}$ be the space of trivialized topologicaly trivial line bundles that is the space of functions $a$ 
 and $\sG_{\C}$ be the group of trivializations.
Then in terms of subsection 6.6 the space  $\sA_{\C}$ is the space of flat $\C^*$-connections,  $\sG_{\C}$ is the complexification of the unitar 
abelian gauge group and
\begin{equation}
Pic_0(P_\Ga) = \sA_{\C} / \sG_{\C}.
\end{equation}
Thus 
\begin{equation}
Pic_0(P_\Ga) = Hom(\pi_1(\Ga), \C^*) = S^a_g
\end{equation}
is the abelian Schottky space (2.68).

To describe a rk 2 bundle $E$ on $P_\Ga$ consider a collection of its trivializations  on all components $P_v$ of $P_\Ga$ and denote the  bundle with such addition structure by $E_0$.

Now to get a rk 2 vector  bundle $E$ on $P_{\Ga}$ we have to glue every
 pair of line bundles $\Oh_{P_{v}} \oplus \Oh_{P_{v}}, \Oh_{P_{v'}} \oplus \Oh_{P_{v'}}$
at common point $p_{e}$ if $v, v' = \p e$. Under our trivializations such concordance is
given by a matrix $ a (\vec e) \in SL(2, \C)$ with the standard
 property (6.57).

 Thus $E_0$ defines a flat $SL(2, \C)$-connection
$$
a \colon \vec E(\Ga) \to SL(2, \C)
$$
(compare with (6.56)). 

The changing of the trivialization is given by a function
$$
\tilde{g} \colon V(\Ga) \to SL(2, \C)
$$
(compare with (6.63)) which acts on the functions $a$  by the formula
$$
\tilde{g}(a(\vec e)) = \tilde{u}(v_s) \cdot a(\vec e) \cdot \tilde{g}(v_t)^{-1}
$$
where $\p \vec e = v_s \bigcup v_t$ (compare with (6.65)).

Thus the space  $\sA_{ \C}$  of trivialized topologicaly trivial rk 2 bundles coincides with the space of  flat  $SL(2, \C)$-connections
 and  the group of trivializations $\sG_{ \C}$ coincides with the complex gauge group.
Then in terms of subsection 6.5
\begin{equation}
M^{ss}_{vb}(P_\Ga) = \sA_{\C} / \sG_{\C}.
\end{equation}
Thus 
\begin{equation}
M^{ss}_{vb}(P_\Ga) = \CLRep(\pi_1(\Ga), SL(2, \C)) = S_g
\end{equation}
is the  Schottky space (see subsection 3.2).

Let us return to the subsection 3.2 and the handlebodies (6.9). Every such
 handlebody defines the exact sequence
\begin{equation}
1 \to  F_g \to \pi_1 (\Si_\Ga) \to \pi_1 (\Ga) \to 1
\end{equation}
where the normal subgroup is a free group with $g$ generators. Using this sequence we can choose the standard generators of $\pi_1 (\Si_\Ga)$ 
$$
\pi_1 = <a_1, ... ,a_g, b_1, ... , b_g \vert \prod_{i=1}^g [a_i, b_i] = 1>
$$
such that
\begin{equation}
F_g = <a_1, ... , a_g>
\end{equation}
$$
\pi_1 (\Ga) =  <b_1, ... , b_g>.
$$
Then we can consider the space of classes of $SL(2, \C)$-representations
\begin{equation}
\CLRep(\pi_1(\Ga), SL(2, \C)) = S_g
\end{equation}
as a family of flat vector bundles on $\Si_\Ga$. If we fix a complex structure
$I$ on $\Si_\Ga$ we get the forgetful map
\begin{equation}
f_I \colon S_g \to M^{ss}(\Si_I)
\end{equation}
(3.44).
The space of complex structures on $\Si_\Ga$ under fixing of the basis  is
the  Teichmuller space $\tau_g$ (5.31) and we can consider the universal family 
of the forgetful maps
\begin{equation}
F \colon S_g \times \tau_g \to \sM
\end{equation}
where $u \colon \sM \to \tau_g$ is the universal family of moduli spaces of vector bundles such that $u^{-1}(I) =  M^{ss}(\Si_I)$.

Now it is easy to see that 
\begin{enumerate}
\item the ll-curve $P_\Ga$ is contained by the Deligne-Mumford 
compactification of $\tau_g$,
\item the moduli space $ \ov{M^{ss}(P_\Ga)}$ is contained by the compactification of
$\sM$ corresponding to  the Deligne-Mumford 
compactification of $\tau_g$,
\item the specialization of the universal forgetful map F ( . ) on $P_\Ga$ coincides with the standard embedding
\begin{equation} 
F_{P_\Ga} = i \colon S_g =  M_{vb}^{ss}(P_\Ga)  \to  \ov{M^{ss}(P_\Ga)}
\end{equation}
described before.
\end{enumerate}

From this it easy to prove  the following partial solution of the Narasimhan problem

\begin{prop} For a general smooth algebraic curve $\Si_I$ of genus $g$ the image of the forgetful
map $f_I : S_g \to M^{ss}(\Si)$ (3.44) is Zariski dense in this moduli space.
\end{prop}
Indeed, this fact is true for $P_\Ga \in \sM_g$ and this property is preserved
by small deformations in $\sM_g$.  (I thank Gerd Detloff for the simplification of my arguments.)

This fact is quite productuve even in the classical abelian case (see subsection 2.4): the standard embedding
\begin{equation} 
i \colon (\C^*)^g = S^a_g \to \ov{Pic_0(P_\Ga)}
\end{equation}
is a limit of the standard forgetful covers $f_I \colon (\C^*)^g \to J(\Si_I)$
(see (2.68)- (2.73)). We can see what happens with the classical theta functions (2.73) under this limit procedure. In particular we can see why 
$\delta$-function (2.46) appears. 

Moreover, the theory of classical theta functions can be constructed by CST-method from subsection 2.3 on $Pic_0(P_\Ga)$ and can be extended to any smooth curve using
 the classical holomorphic flat connection perfectly described by Welters in
\cite{W}. We use the same strategy as well as in the non-abelian case:
CST-constructions from section 6 can be interpreted as the theory of non-abelian theta functions on $ M_{vb}^{ss}(P_\Ga)$. After this  theory
can be extended to any smooth curve using the Hitchin projective flat connection. In this situation to avoid the meromorphicy of the forgetful map
we apply the geometry of Lagrangian and Legandrian cycles and Bohr-Sommerfeld
quantization. 

Here we use  existence of an algebraic-geometrical compactification
of the spaces of vector bundles on ll-curves. Now we describe these compactifications precisely.

\subsection{Compactifications and AG-theory of  theta functions
 for ll-curves}

Both of our spaces $Pic_0(P_\Ga)$ and $M^{ss}_{vb}(P_\Ga)$ are
affine varieties admitting canonical compactifications.
The compactification in the abelian case  will be called {\it jacobian}:
 \begin{equation}
J(P_{\Ga}) = \ov{Pic_{0}(P_{\Ga})}.
\end{equation}

 Every rk 1
torsion free sheaf $F$ is given by a function
$a$ as before  but  we have to
compactify the target space  by two points $0$ and $\infty$. So the new target
space is the projective line $\PP^{1}_{\C^{*}}$ with fixed points
$0$ and $\infty$ or fixed standard $\C^{*}$-action
(distinguishing fixed points $0$ and $\infty$).

So the compactification is
\begin{equation}
J(P_{\Ga}) = (\PP^{1}_{\C^{*}})^{3g-3} / \sG_{\C^{*}}.
\end{equation}
The boundary is the union of $6g-6$ divisors
\begin{equation}
D_{\Ga} = J(P_{\Ga}) - Pic_{0}(\Ga) = \bigcup_{e \in E(\Ga)}
(D_{+}(e) \bigcup D_{-}(e))
\end{equation}
where
\begin{equation}
D_{+}(e) = \{ a \vert a (\vec e) = 0\},
\end{equation}
\begin{equation}
D_{-}(e) = \{ a \vert a (\vec e) = \infty\}
\end{equation}
and
\begin{equation}
D_{+}(e) \bigcup D_{-} (e)= \{F \vert  p_{e} \subset Sing F \}.
\end{equation}

Now consider an abelian spin network $\Ga_w$ (6.99)  on our graph $\Ga$
given by a function $w$ (6.99) subjecting to  the conditions (6.100).
Then we have  divisor
\begin{equation}
D_w = \sum w(e)\cdot D_{sign (w(e))}
\end{equation}
with properties
\begin{enumerate}
\item this divisor is principal that is it is divisor of rational
function $\theta_w$;
\item this function is regular on $Pic_0(P_\Ga)$ being the  limit of the ordinary theta
functions for smooth curves;
\item if $\Ga_w \in SNW^k_a$ that is $w $ is of level $ k$ (6.104)-(6.107) Then
we can fix a pole to make all $w(e)$ non negative and get the canonical basis
of theta functions with characteristics. 
\end{enumerate}

 Every rk 2
torsion free sheaf $F$ is given by the same type function
$a$  but  we have to consider as the target space  the compactification
\begin{equation}
SL(2, \C) \subset \PP^3 = \PP End(\C^2).
\end{equation}
with fixed quadric $Q$.

So 
\begin{equation}
M^{ss}(P_{\Ga}) = \ov{M^{ss}_{vb}}(P_\Ga) = (\PP^{3})^{E(\Ga)} / \sG_{\C}
\end{equation}
where the gauge group is the same (9.79).

 The boundary $D = M^{ss}(P_{\Ga}) - M^{ss}_{vb}(P_{\Ga}) $ is the union of $3g-3$ divisors
\begin{equation}
D =   \bigcup_{e \in E(\Ga)}
D(e). 
\end{equation}
Recall (see the text after formula (6.33) and Fig.12) that we can remove the edge $e$ from $\Ga$ and forget  the vertices $v_1, v_2 = \p e$. We get a new graph
$\Ga_e$ of genus $g-1$. Then the corresponding divisor can be expressed
\begin{equation}
D(e) = Q_e \times M^{ss}(P_{\Ga_e}).
\end{equation}

For every spin network $\Ga_w$ of level 1 the divisor
\begin{equation}
\Theta_w = \sum_{e \in E(\Ga)} w(e) \cdot D(e)
\end{equation}
is the limit of effective theta divisors on smooth curves and the collection of these divisors
\begin{equation}
\{ \Theta_w\} \quad , w \in SNW_\Ga^1 = W^1(\Ga)
\end{equation}
gives the Mumford basis (2.13) which coincides with the Bohr-Sommerfeld basis (4.39) in the wave function space (see Definition 4) for
ll-curve $P_\Ga$. 

Let us fix a spin network of level 1 $w_0 \in W^1(\Ga)$. Then one has
\begin{enumerate}
\item  the divisor $k \Theta_{w_0}$;
\item the collection of divisors
\begin{equation}
\{ \l.c.m.(\Theta_w, k \Theta_{w_0})\} \quad , w \in SNW_\Ga^k = W^k(\Ga)
\end{equation}
\item this collection of effective divisors from the complete linear system $| k \Theta_{w_0}|$ gives the  Bohr-Sommerfeld basis (4.39) in the wave function space  for
ll-curve $P_\Ga$ and the polarization of level $k$.
\end{enumerate}
By the restriction to $uS_g$ we can check that the extension of this basis to smooth curves due to the Hitchin connection coincides with the Bohr-Sommerfeld basis given by BPU-construction.

 \end{document}